\documentclass[letterpaper,11pt,oneside,reqno]{article}
\usepackage{bm}
\usepackage{mathrsfs}
\usepackage{wrapfig}
\usepackage{amsfonts,amsmath, amssymb,amsthm,amscd,stmaryrd,mathtools}
\usepackage[height=9.6in,width=5.95in]{geometry}
\usepackage[mpexclude,DIV13]{typearea}
\usepackage{verbatim}
\usepackage{hyperref}
\usepackage{graphicx}
\usepackage[latin1]{inputenc}
\usepackage{latexsym}
\usepackage{lscape}
\usepackage{epsfig}
\include{bibtex}

\usepackage{setspace}

\usepackage{parskip}
\newcommand{\e}[0]{\epsilon}

\newcommand{\PP}{\ensuremath{\mathbb{P}}}

\newcommand{\N}{\ensuremath{\mathbb{N}}}
\newcommand{\R}{\ensuremath{\mathbb{R}}}

\newcommand{\Z}{\ensuremath{\mathbb{Z}}}
\newcommand{\Q}{\ensuremath{\mathbb{Q}}}

\newcommand{\E}[0]{\mathbb{E}}

\newtheorem{theorem}{Theorem}[section]

\newtheorem{lemma}[theorem]{Lemma}
\newtheorem{proposition}[theorem]{Proposition}
\newtheorem{corollary}[theorem]{Corollary}

\theoremstyle{definition}

\theoremstyle{definition}

\theoremstyle{definition}
\newtheorem{definition}[theorem]{Definition}

\theoremstyle{definition}

\theoremstyle{definition}

\newcommand{\intint}[1]{\llbracket 1,#1 \rrbracket}

\newcommand{\mc}{\mathcal}

\newcommand{\dist}{\vert\vert}

\newcommand{\macthin}{C}
\newcommand{\constfine}{L}

\newcommand{\downward}{L}

\newcommand{\base}{\mathrm{base}}
\newcommand{\head}{\mathrm{head}}

\newcommand{\upwell}{\mathrm{UpClust}}
\newcommand{\upwellu}{\mathrm{UpClust}^u}
\newcommand{\upwellv}{\mathrm{UpClust}^v}
\newcommand{\upwellx}{\mathrm{UpClust}^x}

\newcommand{\downwell}{\mathrm{DownClust}}
\newcommand{\midwell}{\mathrm{MidClust}}
\newcommand{\midcluster}{\mathrm{MidClust}}

\newcommand{\midfragment}{\mathrm{MidFragment}}
\newcommand{\residualwell}{\mathrm{ResidualClust}}
\newcommand{\northzone}{\mathrm{NorthZone}}

\newcommand{\southzone}{\mathrm{SouthZone}}

\newcommand{\midzone}{\mathrm{Tropics}}
\newcommand{\tropics}{\mathrm{Tropics}}
\newcommand{\equator}{\mathrm{Equator}}

\newcommand{\downcluster}{\mathrm{DownClust}}
\newcommand{\upcluster}{\mathrm{UpClust}}

\newcommand{\hmax}{h_{\rm max}}
\newcommand{\hadj}{h_{\rm adj}}

\newcommand{\Base}{\mathrm{Bed}}
\newcommand{\Top}{\mathrm{Top}}

\newcommand{\fine}{\mathsf{Fine}}

\newcommand{\success}{\mathsf{SimpleJoin}}
\newcommand{\catastrophe}{\mathsf{Catastrophe}}
\newcommand{\cmni}{\mathsf{CrossMeagNearIsland}}

\newcommand{\veryfine}{\mathsf{VeryFine}}

\newcommand{\tenmac}{3}

\newcommand{\sting}{\mathsf{Str}}
\newcommand{\culdesac}{\mathsf{CulDeSac}}
\newcommand{\culdesacset}{\mathrm{CdS}}
\newcommand{\cds}{\mathsf{CulDeSac}}
\newcommand{\cdsset}{\mathrm{CdS}}
\newcommand{\icds}{\mathsf{InvCulDeSac}}
\newcommand{\icdsset}{\mathrm{ICdS}}

\newcommand{\stingset}{S}

\newcommand{\cluster}{\mathsf{Cluster}}

\newcommand{\clusteruthree}{\mathsf{Cluster}_{[u,u+3r]}}
\newcommand{\clusteruthreeopen}{\mathsf{Cluster}_{[u,u+3r)}}

\newcommand{\cset}{\mc{C}}

\newcommand{\csetuthree}{\mc{C}_{u,u+3r}}
\newcommand{\starcsetuthree}{\cluster^*_{[u,u+3r]}}
\newcommand{\overlinecsetuthree}{{\overline\cset}_{u,u+3r}}
\newcommand{\ucs}{\mathrm{UCS}}

\newcommand{\boxbig}{{\rm BigBox}}

\newcommand{\Long}{\mathsf{Long}}
\newcommand{\Short}{\mathsf{Short}}

\newcommand{\const}{K}

\newcommand{\surfacek}{\mathrm{Surface}_{\const}}

\newcommand{\slide}{{\rm Slide}}
\newcommand{\iter}{{\rm IterateSlide}}

\newcommand{\slab}{{\rm Slab}}

\newcommand{\renmin}{r_{\rm min}}
\newcommand{\renmax}{r_{\rm max}}

\newcommand{\sausageset}{\mathscr{S}}

\newcommand{\badcoll}{\mathsf{Bad}}

\newcommand{\kappasting}{\kappa_{\rm string}}

\newcommand{\kappacluster}{\kappa_{\rm clust}}

\newcommand{\maxrenewalgap}{{\rm MaxRenGap}}

\newcommand{\height}{{\rm height}}
\newcommand{\macb}{R}

\newcommand{\macd}{D}
\newcommand{\kay}{2(2d-1)p^{-1}K \macb}
\newcommand{\kaypbytwo}{K \macb}
\newcommand{\kaymac}{\macd}
\newcommand{\kaypbytwomac}{\macd \tfrac{p}{2(2d-1)}}

\newcommand{\fv}{J_v}

\newcommand{\rlevel}{{\rm RenLev}}

\newcommand{\rsetmac}{R}

\newcommand{\vecell}{{\pmb{\ell}}}
\newcommand{\cvl}{\cdot \vecell}
\newcommand{\uclust}{[\mathfrak{u}]}
\newcommand{\dclust}{\mathfrak{d}}

\newcommand{\simvecell}{\stackrel{\vecell}{\sim}}
\newcommand{\thickapproxvecell}{\stackrel{\vecell}{\thickapprox}}
\newcommand{\openint}[1]{(#1)}

\newcommand{\fubuthree}{\Z^d \cap F_u \cap B_{u+3r}}

\newcommand{\vfsOne}{{\bf S}_1(\macb)}

\newcommand{\ftOne}{{\bf T}'_1(R)}
\newcommand{\ftTwo}{{\bf T}'_2(R)}
\newcommand{\fsOne}{{\bf S}'_1(\macb)}
\newcommand{\fsTwo}{{\bf S}_2(\macb)}
\newcommand{\fg}{{\bf G}(\constfine)}
\newcommand{\dproj}{d_{{\rm proj}}}

\newcommand{\vertbdry}{\mathrm{VertBdry}}
\newcommand{\edgebdry}{\mathrm{EdgeBdry}}

\newcommand{\rmax}{r_{\rm max}}
\newcommand{\rmin}{r_{\rm min}}

\newcommand{\cover}{{\rm Cover}}
\newcommand{\Nmac}{N_0}
\newcommand{\lat}{\mathrm{Lat}}
\newcommand{\undisc}{\mathrm{UnDisc}}
\newcommand{\termin}{T}
\newcommand{\seed}{\mathrm{Seed}}

\newcommand{\nkk}{\overline{N}_{k,k+1}}
\newcommand{\nkell}{\overline{N}_{k,k+\ell}}

\makeatletter
\newcommand\xxleftrightarrow[2][]{%
  \ext@arrow 9999{\longleftrightarrowfill@}{#1}{#2}}
\newcommand\longleftrightarrowfill@{%
  \arrowfill@\leftarrow\relbar\rightarrow}
\makeatother

\usepackage{titlefoot}

\begin{document}


\title{Sharp asymptotics for finite point-to-plane connections in supercritical bond percolation in dimension at least three}

\author{Alexander Fribergh and Alan Hammond}
\AtEndDocument{
  \bigskip
  \small
  \par
  \textsc{Alexander Fribergh} \\
  \textsc{Montreal} \\
   \textsc{Canada} \\
  \textit{Email:} \texttt{alex@montreal.canada} \\
  
  }

\AtEndDocument{
  \bigskip
  \small
  \par
  \textsc{Alan Hammond} \\
  \textsc{Departments of Mathematics and Statistics, U.C. Berkeley} \\
   \textsc{899 Evans Hall, Berkeley, CA, 94720-3840, U.S.A.} \\
  \textit{Email:} \texttt{alanmh@berkeley.edu} \\

}


\maketitle

\begin{abstract}
We consider supercritical bond percolation in $\Z^d$ for $d \geq 3$. The origin lies in a finite open cluster  with positive probability, and, when it does, the diameter of this cluster has an exponentially decaying tail. For each unit vector $\vecell$, we prove sharp asymptotics for the probability that this cluster contains a vertex $x \in \Z^d$ that satisfies $x \cvl \geq u$.
For an axially aligned $\vecell$, we find this probability  to be of the form $\kappa \exp \{ - \zeta u \}(1+ {\rm err})$ for $u \in \N$, where $\vert {\rm err} \vert$ is at most $C \exp \{ - c u^{1/2} \big\}$;
for general $\vecell$, the form of the asymptotic depends on whether $\vecell$ satisfies a natural lattice condition.
To obtain these results, we prove that renewal points in long clusters are abundant, with a renewal block length whose tail is shown to decay as fast as $C \exp \big\{ - c u^{1/2} \big\}$.

\end{abstract}


\vspace{-2mm}

\setcounter{tocdepth}{2}
\tableofcontents


\section{Introduction}

We investigate the geometry and probability of finite connections in supercritical bond percolation on~$\Z^d$ for dimension $d$ at least three.
 We will present strongly quantified conclusions regarding the presence of renewal points in finite connections. These lead to precise estimates in Theorems~\ref{t.stingwelltrap.cluster} and~\ref{t.stingwelltrap.sting} for the probability of long point-to-hyperplane connections running in a given direction~$\vecell$. 



\subsection{The first main result, on the asymptotic probability of long finite clusters}

Write 
$S^{d-1} = \big\{ u \in \R^d: \vert \vert u \vert\vert_2 = 1 \big\}$ for the Euclidean unit ball.
\begin{definition}
Let $\vecell \in S^{d-1}$. 
This vector is called non-lattice if the $\Z$-linear span of its Cartesian coordinate projections $\vecell \cdot e_i$, $i \in \intint{d}$, is dense in~$\R$.
In the opposing case, where $\vecell$ is lattice, we define the span of $\ell$ to be the greatest positive $s$
for which $\big\vert \vecell \cdot e_i \big\vert \in s \cdot \N$ for $i \in \intint{d}$. (For definiteness,~$0 \in \N$.)

When $\vecell$ is lattice, the cover $\cover(\vecell)$ of $\vecell$ is $s \cdot \N$, with $s$ the span of $\vecell$. When $\vecell$ is non-lattice,  $\cover(\vecell)$  is simply $(0,\infty)$.
\end{definition}

\begin{definition}\label{d.asymptotic}
Let $f,g:\N \to (0,\infty)$.
We will use the convention that 
 $f \simvecell g$ means that $f(u)/g(u) \to 1$ as $u \to \infty$ 
 along the cover of $\vecell$. 
 We will say that $f$ and $g$ are {\em weakly asymptotic} in this case.
 (When $\vecell$ is non-lattice, this notion is simply the usual asymptotic~$\sim$; we will write  $f \simvecell g$ only in the lattice case.)
 Suppose that the pair of sequences $(f,g)$ meets the stronger condition that there exist constants $C$ and~$c$ in $(0,\infty)$ such that 
 $\big\vert f(u)/g(u) - 1 \big\vert \leq C \exp \{ -c u^{1/2} \}$ for $u$ in the cover of~$\vecell$. Then we say that $f$ and $g$ are {\em strongly} asymptotic, denoting this relation by~$f \thickapproxvecell g$.
\end{definition}


Let $d \geq 2$. By $(\Z^d,\sim)$, we denote the lattice $\Z^d$ with the nearest-neighbour notion of adjacency. Writing $p_c = p_c(d)$ for the critical value for bond percolation on $\Z^d$, we let $p \in (p_c,1)$ be a given supercritical value. Let $\PP$ denote the law of bond percolation on  $(\Z^d,\sim)$ with parameter~$p$.

The height of $x \in \Z^d$ equals $x \cdot \vecell$.

\begin{definition}
A cluster is a finite open connected component under $\PP$. 
The base of this cluster is a vertex therein attaining the minimum height.
(If there is more than one minimizer, one is selected as the base in some definite way, as we will specify in Section~\ref{s.maincharacters}.)  
When it exists, the cluster with base at $x \in \Z^d$ will be denoted $\mc{C}(x)$.

For $x \in \Z^d$ and $u \in (0,\infty)$, let $\cluster_u(x)$ denote the event that $x$ is the base of a cluster for which the maximum height attained among its vertices equals~$u$.
\end{definition}
The connection formed when $\cluster_u(x)$ occurs is {\em finite}, because $\vert \mc{C}(x) \vert < \infty$, and of {\em point-to-plane} type, in the sense that an open path travels from $x$ to the hyperplane of height~$u$.

Here is one of the principal conclusions of this article.
\begin{theorem}\label{t.stingwelltrap.cluster}
Let $d \geq 3$, $p \in (p_c,1)$ and $\vecell \in S^{d-1}$. 
Then there exist positive and finite constants $\kappacluster = \kappacluster(d,p,\vecell)$ and $\zeta = \zeta(p,d,\vecell)$ such that, if $\vecell$ is lattice, then 
$$ 
 \PP \big( \cluster_u \big) \thickapproxvecell \kappacluster \exp \big\{ - \zeta u \big\} \, ;
  $$
  while, if $\vecell$ is non-lattice, then
$$ 
 \PP \big( \cluster_{[u,\infty)} \big) \sim \kappacluster \exp \big\{ - \zeta u \big\} \, ,
  $$
  where $\cluster_{[u,\infty)} = \cup_{u' \geq u} \cluster_{u'}$.
\end{theorem}
The value of  $\zeta(p,d,\vecell)$ will be furnished by Proposition~\ref{p.zeta}.



Fundamental to the derivation of Theorem~\ref{t.stingwelltrap.cluster} is the notation of a renewal point and a cousin of the cluster called a string. 
The introduction continues with suitable notation; Theorem~\ref{t.stingwelltrap.sting}, which is a string counterpart of Theorem~\ref{t.stingwelltrap.cluster}; and
 Theorem~\ref{t.renewal}, on the abundance of renewal points. 

\subsection{The basics}\label{s.basics}

Let  $x \cdot y$ denote the Cartesian scalar product of $x,y \in \R^d$.

 We regard $\vecell$ as pointing vertically in the positive direction; so $x \in \R^d$
is upwards from $z \in \R^d$ if $x \cvl > z \cvl$.

For $h \in \R$, the {\em forward} and {\em backward} half-spaces of $h$ are given by
$$
 F_h = \big\{ u \in \R^d: u \cdot \vecell \geq h \big\}
 \, \, \, \textrm{and} \, \, \, B_h =  \big\{ u \in \R^d: u \cdot \vecell < h \big\} \, . 
$$

For $e = (x_1,x_2) \in E(\Z^d)$, write $\openint{e} = \big\{ \lambda x_1 + (1-\lambda) x_2: \lambda \in (0,1) \big\} \subset \R^d$
for the interior of the endpoint-interpolating interval. 
For $h \in \R$, we define the upper and lower half-spaces $U_h$ and $L_h$ of~$h$ to be subgraphs of $(\Z^d,\sim)$ whose edge sets are given by
$$
 E(U_h) = \big\{ e \in E(\Z^d): (e) \cap F_h \not= \emptyset \big\} \, \, \, \textrm{and} \, \, \,  E(L_h) = \big\{ e \in E(\Z^d): (e) \subset B_h   \big\} \, ,
$$ 
and whose vertex sets $V(U_h)$ and $V(L_h)$
are the sets of endpoints of elements of the respective edge sets.

\begin{figure}[htbp]
\centering
\includegraphics[width=0.75\textwidth]{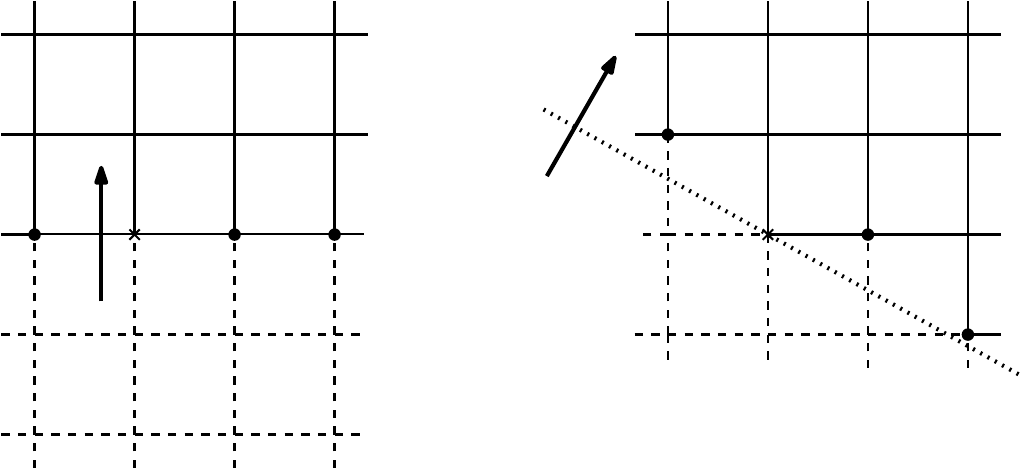}
\caption{The partition $E(L_h) \cup E(L_h)$ of edges in $\Z^2$ for an axial and a more general choice of $\vecell$. The arrow represents $\vecell$, while the cross marks the origin, which is an element of $\vertbdry_h$ for $h=0$. The discs mark the other elements of $\vertbdry_0$. Elements of $E(U_0)$ are drawn with solid lines; those of $E(L_0)$ are dashed.}\label{f.fblusketch}
\end{figure}

The {\em vertex boundary} $\vertbdry_h$ of $h$ equals the intersection of vertex sets $V(U_h) \cap V(L_h)$. The various characters appear in Figure~\ref{f.fblusketch}.


For $m,n \in \Z$ with $m \leq n$, we denote by $\llbracket m,n \rrbracket$ the integer interval $\big\{ k \in \Z: m \leq k \leq n \big\}$.

Let $G_i = \big( V(G_i), E(G_i) \big)$, $i \in \intint{2}$, be two subgraphs of $(\Z^d, \sim)$. Their intersection $G_1 \cap G_2$
is the subgraph $\big( V(G_1) \cap V(G_2), E(G_1) \cap E(G_2) \big)$. 
For $y,z \in \R$ with $z \geq y$,
we define the slab $\slab_{y,z}$  from $y$ to $z$ to be the subgraph  $U_y \cap L_z$.

An element $v \in V(\slab_{y,z})$ is called lower if it belongs to $\vertbdry_y$, and upper if it belongs to~$\vertbdry_z$.

\subsection{Notation, strings, and renewal}\label{s.maincharacters}

\subsubsection{The base and the head}\label{s.baseandhead}
For certain choices of $\vecell$, the hyperplane $H_s = \{ z \in \R^d: z \cvl = s \}$, with $s \in \R$, may intersect a finite set $A \subset \Z^d$ at more than one point. To distinguish a point in the intersection, we consider the projection $H_s \to H_0$ orthogonal to $\vecell$. We fix an orthonormal basis for~$H_0$.
This basis induces a lexicographical ordering on $H_0$ and this ordering gives rise to one on~$H_s$ when we identify $H_0$ and $H_s$ by translation by any given element in~$H_s$.
The latter ordering is a total order on~$H_s$.

Let $G$ be a finite subgraph of $(\Z^d,\sim)$. The base of $G$, denoted by $\base(G)$, is the lexicographically minimal element of $V(G)$ among those of minimum height.
For the head $\head(G)$, we replace `minimal' by `maximal' and `minimum' by `maximum'.
The usage of the term `height' is extended so that the height of $G$ is given by $\max \{ (v-u) \cdot \vecell : u,v \in V(G) \}$.

\subsubsection{Renewal}
Next is the important definition of a renewal point.

\begin{definition}\label{d.renewal}
Let $U$ be a finite connected subgraph of $(\Z^d,\sim)$
and let $x \in V(U)$.
The element $x$ is called a {\em renewal point} of~$U$ if $V(U) \cap \vertbdry_{x \cvl}$
equals the singleton set $\{ x \}$.
\end{definition}

\subsubsection{The string}\label{s.tubestingcluster}
\begin{figure}[htbp]
\centering
\includegraphics[width=0.35\textwidth]{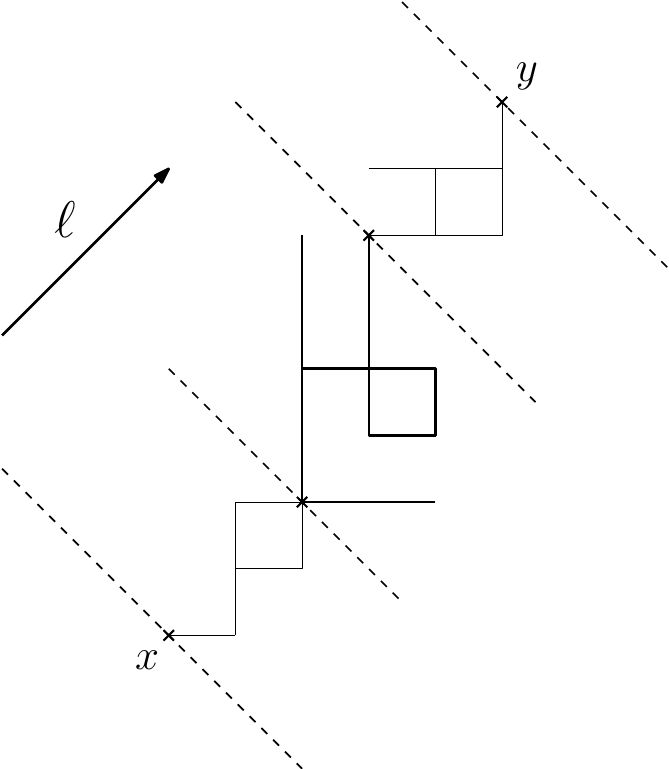}
\caption{Four crosses mark the renewal points in this string from $x$ to $y$. The middle of its three sausages is drawn in bold.}\label{f.sting}
\end{figure}

Let $x,y \in \Z^d$, $x \cvl < y \cvl$.
A string from $x$ to $y$ 
is a finite open connected component in $\slab_{x \cvl,y \cvl}$ 
such that $x,y$ are renewal points of~$S$. A string crosses the slab $\slab_{x \cvl,y \cvl}$ between locations $x$ and~$y$ at which it has renewal points.
The term `string' is shorthand for `string of sausages': as we will see, a string can be broken into sausages delimited by consecutive renewal points: see Figure~\ref{f.sting}.

For $u > 0$ and $x \in \Z^d$, let $\sting_u(x)$ denote the event 
there exists a string from $x$ to $z$ for some $z \in \Z^d$ with $(z - x) \cvl = u$.
(Whatever the value of $\vecell$, this event is empty for  some $u$.)
When $\sting_u(x)$ occurs, it is straightforward, and implied by the upcoming Lemma~\ref{l.stingbasics}, that there is a unique such~$z$, and a unique string from $x$ to $z$. 
We denote this string by $\stingset_u(x)$.
When $x = 0 \in \Z^d$, write $\sting_u = \sting_u(0)$ and $\stingset_u = \stingset_u(0)$.


\subsubsection{Notation}\label{s.notation}

The next definition develops the usage $\cluster_{[u,\infty)}$ seen in Theorem~\ref{t.stingwelltrap.cluster}.
\begin{definition}\label{d.enotation}
Here the symbol $E$ represents either $\sting$ or $\cluster$.
For $A \subseteq [0,\infty)$, we write $E_A = \cup_{u \in A} E_u$.
This notation will be applied when $A$ is a finite or semi-infinite interval.
\end{definition}

Set $r = \max \big\{ \vert e_i \cvl \vert : i \in \intint{d} \big\}$.

\subsection{The asymptotic probability of long strings}

The next  result asserts that
the string probability $\PP \big( \sting_{[u,u+r)} \big)$ decays in the form $\exp \{ - (\zeta + o(1))u \}$. 
In essence, the assertion is due to~\cite[Lemma~$8.5$]{FriberghHammond}, and we provide the straightforward derivation given this input at the end of Section~\ref{s.asp}. 
\begin{proposition}\label{p.zeta}
There exists $\zeta = \zeta(p,d,\vecell) \in (0,\infty)$ such that 
$$
 \lim_u u^{-1} \log \PP \big( \sting_{[u,u+r)} \big) = - \zeta \, .
$$
\end{proposition}

Here is the string counterpart to Theorem~\ref{t.stingwelltrap.cluster}.
\begin{theorem}\label{t.stingwelltrap.sting}
Let $d \geq 3$, $p \in (p_c,1)$ and $\vecell \in S^{d-1}$. 
 Write $\zeta$ for the constant $\zeta(p,d,\vecell) \in (0,\infty)$ in Proposition~\ref{p.zeta}. 
Then there exists a  constant $\kappasting = \kappasting(d,p,\vecell) \in (0,\infty)$
such that, if $\vecell$ is lattice, then 
$$ 
 \PP \big( \sting_u \big) \thickapproxvecell \kappasting \exp \big\{ - \zeta u \big\} \, ;
  $$
  while, if $\vecell$ is non-lattice, then
$$ 
 \PP \big( \sting_{[u,\infty)} \big) \sim \kappasting \exp \big\{ - \zeta u \big\} \, .
  $$
\end{theorem}

\subsection{Renewal results}

Let $U$ be a subgraph of $(\Z^d,\sim)$. A renewal level of $U$ is the height of a renewal point of~$U$. Let $\rsetmac(U)$, and $\rlevel(U)$, denote the set of renewal points, and of renewal levels, of~$U$.

We will write $\mc{C} = \mc{C}(0)$ denote the cluster with base at~$0$.
Recall from Definition~\ref{d.enotation} that $\clusteruthree$ denotes the event $\cup_{u \leq h \leq u+3r}\cluster_h$, under which $\mc{C}$  has height in $[u,u+3r]$.  
\begin{theorem}\label{t.renewal}
\leavevmode
\begin{enumerate}
\item 
There exist $c > 0$ and $u_0 \in \N$ such that, for $u \in [u_0,\infty)$,
$$
 \PP \Big(  \big\vert R(\cset) \big\vert \geq cu  \, \Big\vert \, \clusteruthree  \Big) \geq c   \, .
 $$ 
\item 
For some $c > 0$,
the condition $u^{1/2} \leq k \leq cu (\log u)^{-1}$ on $u \in (0,\infty)$ and $k \in \N$ implies that
$$
\PP \Big( \big\vert R(\cset) \big\vert \leq k \, \Big\vert \, \clusteruthree  \Big) \, \leq \, \exp \big\{ -cu/k \big\} \, .
$$
\end{enumerate}
\end{theorem}
In this theorem, and throughout, the positive constant $c > 0$ may depend on $d \geq 3$ and $p \in (p_c,1)$.   
\begin{corollary}\label{c.renewal}
Consider Theorem~\ref{t.renewal} with the event $\clusteruthree$ in the conditioning 
 replaced by $\clusteruthree \cap \big\{ \{0,h\} \subseteq  \rlevel(\mc{C}) \big\}$ where $h$ denotes the height of $\head(\mc{C})$.
 In other words, condition not only on  $\clusteruthree$ but also on the resulting cluster having renewal points at the highest and lowest possible locations, namely $\base(\mc{C}) = 0$  and $\head(\mc{C})$. Alternatively, we may condition 
\begin{enumerate}
\item Thus altered, Theorem~\ref{t.renewal}(1) remains valid. 
\item And likewise Theorem~\ref{t.renewal}(2).
\item Theorem~\ref{t.renewal}(1) remains valid if the conditioning is instead replaced by  $\clusteruthreeopen \cap \big\{ \{0,h\} \subseteq  \rlevel(\mc{C}) \big\}$, wherein the interval is open on the right.
\end{enumerate}
\end{corollary}


\subsection{Review of and comparison with prior work}

Ornstein and Zernike~\cite{OrnsteinZernike} presented a theory that offers a sharp asymptotic description of the density correlations in classical fluids by assuming the validity of a form of locality expressed via a renewal equation.
Finding rigorous manifestations of their theory amounts to giving probabilistic descriptions of objects such as long clusters in terms of essentially one-dimensional chains of irreducible components.
Early rigorous progress was made with perturbative arguments, such as~\cite{MinlosZhizhina} for the Ising model at very high temperature.
Self-avoiding walk with an exponential penalty for length was addressed for axial directions in~\cite{ChayesChayes}
and it offered the setting for an early non-perturbative analysis in~\cite{Ioffe}. 
 Campanino and Ioffe~\cite{CampaninoIoffe}
  developed a non-perturbative approach for  subcritical percolation on Euclidean lattices in  work that provided  sharp asymptotic estimates for the probability that distant pairs of points are connected. 
The analogous task for finite range Ising models  above the critical temperature  was undertaken by Campanino, Ioffe and Velenik in~\cite{CampaninoIoffeVelenik03} and for  the subcritical random cluster model  by these same authors in~\cite{CampaninoIoffeVelenik}.
The papers~\cite{CampaninoIoffe,CampaninoIoffeVelenik03,CampaninoIoffeVelenik} employ a ball-and-point search algorithm that marks points on a lengthy connection and tracks where the connection leaves a ball of fixed radius (in a suitable norm) centred at the marked point. The algorithm's output offers a coarse-grained description of the long connection. In these analyses, the rare object, a long subcritical connection, is being detected and explored by a local procedure. Here a basic difference arises with our problem, where a long but finite open connection forms in supercritical percolation. For us, it is not the presence of a long open path that is rare, but rather that this path is disjoint from the infinite open component. Indeed, the long finite path is enclosed by a surface formed of dual plaquettes that separates it from other open structure. This enclosing surface roughly plays the role analogous  to the long subcritical connection in~\cite{CampaninoIoffe} and~\cite{CampaninoIoffeVelenik}. But the surface in principle may lie far from a given point that it encloses, and there is no clear counterpart to locally operating ball-and-point algorithms that might roughly represent its form. This basic difference led us to consider very different techniques for proving asymptotics, and the resulting proofs that we present here offer new tools for deriving sharp Ornstein-Zernike asymptotics. 

Campanino, Ioffe and Louidor~\cite{CampaninoIoffeLouidor} find sharp asymptotics for finite point-to-point connections in supercritical bond percolation in dimension two. This counterpart to our work with $d=2$
has a different flavour to ours, because the trapping surface whose interior contains two given distant points is a subcritical loop in the dual model, which is itself an instance of bond percolation in this dimension: thus~\cite{CampaninoIoffeLouidor} may be viewed as a development of the Ornstein-Zernike theory of subcritical bond percolation.   

The need for new tools is also illustrated by the way in which the coarse-grained algorithmic output description of a long connection is, in the percolation case, analysed in~\cite[(2.2)]{CampaninoIoffe} by means of the BK inequality: disjoint connection probabilities are there bounded above by a product of marginal probabilities. In our context of finite connection, a comparable monotonicity is not available:  if two points are finitely connected, then there is an open path between them which is surrounded by a surface of dual plaquettes, and these structures may or may not facilitate finite connection with a third point. 
In the random cluster case of~\cite{CampaninoIoffeVelenik}, an exponential decay of correlations is used to gain independence between modestly separated balls. A roughly counterpart assertion in our context might claim that, if a certain collection of points is finitely connected, then a dual plaquette surface will enclose all these points whose cardinality exceeds that of the minimum one necessary by a random amount whose law has a rapidly decaying tail.  
This however is a delicate assertion, for which the renewal abundance Theorem~\ref{t.renewal} may offer technical support.

Rigorous implementations of Ornstein-Zernike theory such as~\cite{CampaninoIoffe,CampaninoIoffeVelenik03,CampaninoIoffeVelenik}  reach conclusions including sharp asymptotics for the probability of point-to-point connections $0 \longleftrightarrow x$, in which a local central limit theorem style factor $\dist x \dist^{-(d-1)/2}$ multiplies an exponential decay rate, by establishing a {\em mass gap} condition. This condition asserts that a moment generating function associated to the connections decays at an exponential rate that is strictly lower than the counterpart function where connection occurs without intervening renewal points. The condition typically entails other favourable features such as the analyticity and strict convexity of the equidecay profile. For mass gap to hold, it is necessary that the conditional probability given a lengthy connection that renewal is absent decay exponentially in the distance travelled. 

The principal conceptual progress made in the present article is the introduction and analysis of a slide resampling procedure to yield the renewal abundance Theorem~\ref{t.renewal}, with the sharp asymptotic point-to-plane Theorems~\ref{t.stingwelltrap.cluster} and~\ref{t.stingwelltrap.sting}  resulting. We emphasise however that the just mentioned condition would correspond to a strengthening of Theorem~\ref{t.renewal}(2)
in which the decay rate $C e^{-c u^{1/2}}$ is improved to $C e^{-cu}$. We anticipate that the point-to-plane asymptotic will be sufficient for accompanying work, which we will mention next, concerning biased walk on the supercritical cluster. We present our results and the technical arguments that demonstrate them in the hope of attracting attention to a new technique in the Ornstein-Zernike vein as well as to offer a fundamental ingredient in a broader study of biased walk in disorder. The problem of proving mass gap is open, with technical refinements of the slide resample offering one possible route to a solution.

\subsubsection{Motivation from random walk in random environment}
A random walk on the infinite open cluster of supercritical percolation on $(\Z^d,\sim)$ that has a bias towards a given preferred direction manifests a phase transition, as the strength of the bias rises, from a regime of linear progress in the long-run to one of sublinear progress in which the walk becomes waylaid for lengthy periods in traps: building on \cite{BGP,Sznitman}, \cite{FriberghHammond} proved the existence of a critical point for this transition. In the sublinear regime, traps resemble lengthy finite percolation connections which have been opened up at one end, through which the walker may enter to become waylaid near the opposing end. Arguments in~\cite{FriberghHammond} were facilitated by percolation results of the like of Proposition~\ref{p.zeta}, with asymptotic estimates on trap probabilities of the form~$\e^{-\zeta u (1+o(1))}$.
The far stronger inferences made in the present article speak to the aim of analysing the sublinear regime at a finer level of detail, and we expect this article to facilitate an analysis of this regime and the scaling limit of the walk therein.

\subsection{Structure of the paper}

Strings are easier objects than clusters to work with in deriving asymptotic finite connection probabilities, so that the string Theorem~\ref{t.stingwelltrap.sting} will be derived before the cluster Theorem~\ref{t.stingwelltrap.cluster}, which is largely a consequence of the former result and its proof. In Section~\ref{s.asp}, we will decompose strings into irreducible building blocks called sausages, and use renewal theory to deliver Theorem~\ref{t.stingwelltrap.sting} from the renewal abundance Theorem~\ref{t.renewal}. 
The proof of the asymptotic cluster probability Theorem~\ref{t.stingwelltrap.cluster} finishes this section.  
Theorem~\ref{t.renewal} is conceptually central to this article, and it is proved by a novel {\em slide} resampling procedure, which is applied once to obtain Theorem~\ref{t.renewal}(1)  and in an iterative scheme to yield Theorem~\ref{t.renewal}(2). A coarse-grained description of percolation involving a lattice of boxes of given width 
will facilitate these resampling arguments. 
In Section~\ref{s.heuristics} we offer an overview of the slide resample and introduce basic aspects of the coarse-graining. 
In the more substantial Section~\ref{s.renewalslide}, 
Theorem~\ref{t.renewal} is proved. 

\subsubsection{Acknowledgments}
The second author is supported by the National Science Foundation under DMS grants~$1855550$ and~$2153359$ and by the Simons Foundation as a $2021$ Simons Fellow.
He thanks Hugo Duminil-Copin and Yvan Velenik for helpful discussions in Geneva around 2012 about finite connections in dimension at least three.

\section{Asymptotic string and cluster probability via renewal plenitude}\label{s.asp}

There are three subsections.
In the first and principal Section~\ref{s.provingsting}, we use renewal theory to prove Theorem~\ref{t.stingwelltrap.sting}.
This result has a weaker counterpart in which the strong asymptotic is replaced by the weak in the lattice case. We will first derive this weaker assertion, obtaining Proposition~\ref{p.hnprob} (in the lattice case). The renewal abundance assertions
 Theorem~\ref{t.renewal}(1) and  Theorem~\ref{t.renewal}(2) will be harnessed in our proofs: the former to obtain the weaker counterpart just mentioned; the latter,  in deriving Theorem~\ref{t.stingwelltrap.sting} itself. In fact, we will use the following consequence  of  Theorem~\ref{t.renewal}(1), which is proved in Section~\ref{s.iterate}.
  \begin{proposition}\label{p.stingrenewal}
  There exists $c > 0$ such that, for $u$ sufficiently high,
  $$
   \PP \, \Big( \, \big\vert R(S) \big\vert \geq \delta u \, \Big\vert \, \sting_{[u,u+r)} \Big) \, \geq \, c \, . 
  $$
  \end{proposition}
  This statement in fact raises a point about notation. We will shortly see, in Lemma~\ref{l.geometrybasics}(2), that the union event $\sting_{[u,u+r)}$
  is characterized by the occurrence of $\sting_{u'}$ for a unique value of $u' \in [u,u+r)$. In the above, and later, we write $S$ on $\sting_{[u,u+r)}$ for the string that realizes $\sting_{u'}$
  for this value of $u'$.
 
 Section~\ref{s.stingconseq} is devoted to deriving several consequences of our analysis of strings and their renewal that will be applied in the companion article.
In Section~\ref{s.asympcluster}, we derive the cluster probability asymptotic Theorem~\ref{t.stingwelltrap.cluster} by showing how a typical lengthy cluster splits into a long middle string that abuts on each side a shorter endpiece.
  
  \subsection{Deriving Theorem~\ref{t.stingwelltrap.sting}}\label{s.provingsting}


Two basic lemmas concerning vertex boundary and renewal point sets will be useful. Here and later, it will be convenient to 
\begin{equation}\label{e.aoc}
\textrm{suppose that $e_i  \cvl \geq 0$  for each $i \in \intint{d}$ and that $e_d \cvl$ attains $\max \big\{  e_i \cvl : i \in \intint{d} \big\}$} \, .
\end{equation}
Reflecting in coordinate axes permits this assumption, which we will call the {\em axes ordering convention}.
Note that, consequently, the quantity~$r$ specified in Subsection~\ref{s.notation} equals  $e_d \cvl$.

\begin{lemma}\label{l.geometrybasics}
Let $h \in \R$ and $v \in \Z^d$.
\begin{enumerate}
\item The set $\vertbdry_h$ equals $\big\{ v \in \Z^d: v \cvl \in (h-r,h]\big\}$.
\item Let $S$ be a finite connected subgraph of $(\Z^d,\sim)$. The heights of distinct elements of $R(S)$ differ by at least $r$.
\end{enumerate}
\end{lemma}
{\bf Proof: (1).} Let $v \in \vertbdry_h$.
 The interior of some edge abutting $v$ lies in $B_h$, so that $v \cvl \leq h$; and the interior of another such edge intersects $F_h$, so that
 $(v+e_d)\cvl > h$
since $e_d \cvl$ attains $r  = \max \vert e_i \cvl \vert$ by the axes ordering convention.
 Hence, 
 $v \cvl > h - r$. 

Now let $v \in \Z^d$ satisfy $v \cvl \in (h-r,h]$. The edge interior $(v,v+e_d)$ intersects $F_h$, so $v \in V(U_h)$. And $(v-e_d,v)$
lies in $B_h$, so $v \in V(L_h)$.

{\bf (2).} If $y \in R(S)$, then $V(S) \cap \vertbdry_y = \{ y\}$. Thus the preceding part implies that no element of $R(S)$ may adopt a height in $(y\cvl -r,y\cvl)$. \qed

\begin{lemma}\label{l.stingbasics}
A string~$S$ contains a unique element of lowest height and a unique element of greatest height. These respective vertices are $\base(S)$ and $\head(S)$.
They are renewal points of $S$. The renewal point set $R(S)$ is strictly and totally ordered according to the elements' height.
\end{lemma}
{\bf Proof.} Recall that any given string $S$ runs from $y$ and to $z$, for certain $y,z \in \Z^d$, $y \cvl < z \cvl$, with $y,z \in R(S)$ and $S \subset \slab_{x\cvl, y \cvl}$.
Since $y \in R(S)$, $V(S) \cap \vertbdry_{y\cvl}$ equals $\{ y\}$. Any element of $V(S)$ lies in $V(U_{y\cvl})$; and it is easy to check that any element of $V(U_{y\cvl}) \setminus \vertbdry_{y\cvl}$
has height strictly exceeding $y \cvl$. Thus, $y$ is the unique element of $V(S)$ of lowest height, which forces $y = \base(S)$. Similar reasoning tells us that $z = \head(S)$ is the unique element of $S$ of greatest height.
Lemma~\ref{l.geometrybasics}(2) implies the strict ordering of renewal heights. \qed 

Lemma~\ref{l.stingbasics} prompts the next definition.
\begin{definition}\label{d.sausage}
A sausage is a string whose vertex set contains no renewal point other than its head and base. Any string~$S$ is uniquely decomposed into sausages via its renewal point set~$R(S)$. 
Indeed, any pair $(u,v)$, $u \cvl < v \cvl$, of vertices of consecutive height in $R(S)$ specifies a slab $\slab_{u \cvl, v \cvl}$, and the subgraph of $S$ induced by $\slab_{u \cvl, v \cvl}$
is a sausage in   the decomposition of~$S$. 

Suppose that the list of sausages of a given string takes the form $\big\{ S_i: i \in \intint{\ell} \big\}$ in increasing order of height. We will denote the decomposition $S = S_1 \circ S_2 \circ \cdots \circ S_\ell$.
\end{definition}
\begin{definition}\label{d.os}
The {\em slab} of a string~$S$ is $\slab_{\base(S) \cvl,\head(S) \cvl}$.
A string $S$ is said to occur if, under~$\PP$, every edge in $S$ is open, and every edge that is incident to $S$ and lies in the slab of $S$ is closed. 
We write $\mathsf{O}(S)$ for the event that $S$ occurs.
\end{definition} 
\begin{definition}\label{d.sausagemu}
Let $Y$ denote the set of strings whose base equals $0 \in \Z^d$. Let $\mu:Y \to (0,\infty)$ be given by $\mu(S) = \exp \{ \zeta \cdot {\rm hgt}(S) \} \PP \big( \mathsf{O}(S) \big)$, where ${\rm hgt}(S)$ is the height of $S$ as this notion is specified in Subsection~\ref{s.baseandhead}, and where Proposition~\ref{p.zeta} provides the constant $\zeta > 0$.
\end{definition}
\begin{lemma}\label{l.percoprod}
Let $S = S_1 \circ S_2 \circ \cdots \circ S_\ell$ be a string. Then
$$
\mu(S) = \prod_{i = 1}^\ell \mu(S_i) \, \, \, \, \textrm{and} \, \, \, \, \PP \big( \mathsf{O}(S) \big) = \prod_{i = 1}^\ell \PP \big( \mathsf{O}(S_i)  \big) \, .
$$
\end{lemma}
{\bf Proof.} The event $\mathsf{O}(S)$ occurs precisely when all of the events $\mathsf{O}(S_i)$, $i \in \intint{\ell}$, do. The latter events are specified by the open-or-closed status of edges where a given edge is involved in determining the occurrence of at most one of the events $\mathsf{O}(S_i)$. Thus, these events are independent, and the latter formula in Lemma~\ref{l.percoprod} is obtained. If we multiply this formula by 
$$
\exp \big\{ \zeta \cdot {\rm hgt}(S) \big\} \, = \, \prod_{i = 1}^\ell \exp \big\{ \zeta \cdot {\rm hgt}(S_i) \big\} \, ,
$$
 then we obtain the former formula  by distributing the terms in the product in the evident way. \qed 
\begin{definition}\label{d.sausageset}
Let $\sausageset$ denote the set of sausages whose base equals $0 \in \Z^d$.
\end{definition}
\begin{proposition}\label{p.muprob}
The sum $\sum_{S \in \sausageset} \mu(S)$ equals one. 
\end{proposition}
A lemma will help to prove this result. 

Let $u > 0$.
Under $\PP$, let $C(u)$ denote the open connected component of $0 \in \Z^d$ 
in $\slab_{0,u}$.
  Let $\mc{S}_u$ denote the set of strings~$S$ with $\base(S) = 0$ and $\head(S) \cdot \vecell = u$.
\begin{lemma}\label{l.essn}
To each $A \subseteq \mc{S}_u$, we may associate the event $\mathsf{A} = \big\{ C(u) \in A \big\}$. For $u >0$ and any such~$A$,
$$
 \PP \big( \mathsf{A} \big) =   \exp \big\{ - \zeta u \big\}   \sum_{S \in A} \mu(S) \, .
$$
\end{lemma}
{\bf Proof.} For $S \in \mc{S}_u$, $\mathsf{O}(S) = \big\{ C(u) = S \big\}$. 
Thus,
$$
 \PP \big( \mathsf{A} \big) =  \sum_{S \in A} \PP \big( \mathsf{O}(S) \big) =  \sum_{S \in A}  \prod_{i = 1}^\ell \PP \big( \mathsf{O}(S_i)  \big) \, ,
$$
where the generic $S \in A$ has been decomposed into sausages according to our convention. The quantity $\prod_{i = 1}^\ell \exp \big\{ \zeta \cdot {\rm hgt}(S_i) \big\}$
equals $\exp \big\{ \zeta u \big\}$ for any string $S \in \mc{S}_u$. Thus, we may multiply the last display by this shared value and use the definition of $\mu$ to obtain
$$
   \exp \big\{ \zeta u \big\} \PP \big( \mathsf{A} \big) \, = \,   \sum_{S \in A}  \prod_{i = 1}^\ell \mu (S_i)   \, .
$$
Since the right-hand summand equals $\mu(S)$ by Lemma~\ref{l.percoprod}, we obtain Lemma~\ref{l.essn}. \qed


{\bf Proof of Proposition~\ref{p.muprob}.}  For $k \geq 2$, let $\mc{S}_u(k)$ denote the set of $S \in \mc{S}_u$ for which $\vert R(S) \vert \geq k$.
Recall Definition~\ref{d.renewal} and note that, by Lemma~\ref{l.essn}, 
$$
 \PP \Big( \sting_u , \big\vert R (\stingset_u) \big\vert \geq cu \Big) \, = \, \exp \big\{ - \zeta u \big\}   \sum_{S \in \mc{S}_u(cu)} \mu(S) \, .
$$
Write $\mc{S}_{[u,u+r)}(k)$ for the union $\bigcup_{h \in [u,u+r)}\mc{S}_h(k)$. Summing the preceding equality over values of $u$ in an interval of the form $[h,h+r)$,
and relabelling $h$ as $u$, we find that
$$
 \PP \Big( \sting_{[u,u+r)} , \big\vert R (\stingset) \big\vert \geq cu \Big) \, \leq \, \exp \big\{ - \zeta u \big\}   \sum_{S \in \mc{S}_{[u,u+r)}(cu)} \mu(S) \, .
$$
We have that
$$
 \sum_{S \in \mc{S}_{[u,u+r)}(k)} \mu(S) \leq   \sum_{j=k-2}^\infty \mu(\sausageset)^j \, ,
$$
because Lemma~\ref{l.percoprod} implies that the right-hand side equals the $\mu$-value of the set of strings $S$ with $\base(S) = 0$ and at least $k-2$ non-trivial renewal points (or, equivalently, with a decomposition having at least $k-1$ sausages).
Suppose that $\mu(\sausageset) < 1$.
Choosing $k = cu$, we see from the last display that  
$$
\PP \big( \sting_{[u,u+r)} , \vert R(\stingset) \vert  \geq cu \big) \leq C \exp \big\{ - (\zeta + \e) u \big\} \, .
$$
 for some $\e > 0$ and a positive constant $C$. 
 By Proposition~\ref{p.stingrenewal}, this bound is equally valid for 
$\PP \big(  \sting_{[u,u+r)} \big)$. 
The bound for 
$\PP \big(  \sting_{[u,u+r)} \big)$ contradicts Proposition~\ref{p.zeta}: thus, we learn that $\mu(\sausageset) \geq 1$.

Suppose now that $\mu(\sausageset) > 1$. For $\ell \in \N$, let $\sausageset_\ell$ denote the set of sausages whose base is $0 \in \Z^d$ and whose head has Euclidean distance~$\vert\vert \cdot \vert\vert$ at most $\ell$ from $0$. When the condition on the head distance is suspended, the resulting set may naturally be denoted $\sausageset_\infty$, and then we have  $\sausageset = \sausageset_\infty$. Choose $\ell \in \N$ such that $\mu(\sausageset_\ell) > 1$. For $n \in \N$, let $\sausageset_\ell^n$ denote the set of strings whose base is $0$ and whose decomposition into sausages has $n$ elements, each of which lies in $\sausageset_\ell$. Lemma~\ref{l.percoprod} implies that
$$
 \sum_{S \in \sausageset_\ell^n} \mu(S) = \mu(\sausageset_\ell)^n \, .
$$
For $S \in \sausageset_\ell$, it follows readily from Lemma~\ref{l.geometrybasics}(2) that  $\vert\vert \head(S) \vert\vert \in [ r,\ell ]$, where $r = e_d \cvl > 0$. Thus $\vert\vert \head(S) \vert\vert \in [ nr,n \ell ]$ whenever $S \in \sausageset_\ell^n$. Choose $z \in \Z^d$
with $\vert\vert z \vert\vert \leq n\ell$
such that
$$
 \sum_{\begin{subarray}{c}  S \in \sausageset_\ell^n :  \\
                            \head(S) = z \end{subarray}}
  \mu(S) \, \geq \, c (n\ell)^{-d} \mu(\sausageset_\ell)^n \, ,
$$
where $c \in (0,\infty)$ satisfies $\limsup_{u \to \infty} u^{-d} \cdot \big\vert \{ z \in \Z^d: \dist u \dist \leq u \} \big\vert < c^{-1}$.
Writing $m = \vert\vert z \vert\vert$, this left-hand side is at most $\exp \{ \zeta m \}   \PP (\sting_m)$ by Lemma~\ref{l.essn}. Since $rn \leq m \leq n\ell$ and   $\mu(\sausageset_\ell) > 1$,
$$
 \PP \big(\sting_m \big) \geq  c r^d (m\ell)^{-d} \Big(  \exp \{ - \zeta  \}  \mu(\sausageset_\ell)^{1/\ell} \Big)^m \, .
$$
From  $\mu(\sausageset_\ell) > 1$, we see that $\limsup_{m \to \infty} m^{-1} \log \PP \big(\sting_m \big) > -\zeta$. 
But this bound 
contradicts  Proposition~\ref{p.zeta}. 
We learn then that $\mu(\sausageset) \leq 1$ and thus complete the proof of Proposition~\ref{p.muprob}. \qed

Proposition~\ref{p.muprob} permits an abuse of notation whereby $\mu$ denotes a probability measure supported on $\mc{S}$.
Let $\mu^{\otimes \N}$ denote the law on infinite strings~$S$ with base at $0 \in \Z^d$ given by the concatenation $S = S_1 \circ S_2 \circ \cdots$ of an independent sequence of sausages drawn from the law~$\mu$. Renewal notions naturally extend to infinite strings. Let $\rlevel$ denote the set of renewal levels in the random infinite string~$S$ specified under the law $\mu^{\otimes \N}$.
In a notational abuse, we will permit $\mu^{\otimes \N} [a,b]$ to denote the mean value of $\rlevel \cap [a,b]$ under $\mu^{\otimes \N}$, for a real interval $[a,b]$.
\begin{lemma}\label{l.zetarenewal}
Let $u \in (0,\infty)$.
\begin{enumerate}
\item
We have that
 $$
  \exp \big\{ \zeta u \big\} \PP \big( \sting_u \big)  = \mu^{\otimes \N} \big(  u \in \rlevel\big) \, .
 $$
 \item Let $A$ denote a measurable collection of finite subsets of $(0,u)$. Then
 $$
 \exp \big\{ \zeta u \big\} \PP \big( \sting_u , \rlevel(\stingset_u) \cap (0,u) \in A  \big)  = \mu^{\otimes \N} \big( u \in \rlevel \, , \, \rlevel \cap ( 0,u )  \in A \big) \, .
 $$
 (Elements of $A$ lie in the union of the spaces $S_k = (0,u)^k$ for $k\in \N_+$. The union is a measurable space when it is embued with the $\sigma$-algebra associated to the Borel $\sigma$-algebras on each $S_k$.)
 \item 
 For $k \in \N$ at least two,
 $$
 \exp \big\{ \zeta u \big\} \PP \big( \sting_u , \big\vert R(\stingset_u) \big\vert = k \big)  = \mu^{\otimes \N} \big( u \in \rlevel \, , \, \big\vert \rlevel \cap ( 0,u ) \big\vert = k-2 \big) \, .
 $$
 \end{enumerate}
Now  suppose that $\vecell$ is non-lattice. 
\begin{enumerate}
  \setcounter{enumi}{3}
 \item For $\phi \in (0,\infty)$,  
 $\lim_{a \to \infty} \mu^{\otimes \N}[a,a+\phi) = \phi r^{-1}\mu_\infty$, 
 where $\mu_\infty = \lim_{u \to \infty} \mu^{\otimes \N} [u,u+r)$.
 \item We have that
 $$
  e^{\zeta u} \PP \big( \sting_{[u,u+r)} \big) \, \to \,  \mu_\infty r \zeta^{-1} (1 - e^{-\zeta r}) \, \, \, \textrm{and} \, \, \,  e^{\zeta u} \PP \big( \sting_{[u,\infty)} \big) \, \to \,  \mu_\infty r \zeta^{-1}  \, .
 $$
 \end{enumerate}
\end{lemma}
{\bf Proof: (1).} For a string~$S$, let $d(S) = \vert R(S) \vert -1$ denote the number of sausages into which $S$ decomposes. 
Note then that 
\begin{eqnarray*}
  \exp \big\{ \zeta u \big\} \PP \big( \sting_u \big) & = & \sum_{S \in \mc{S}_u} \mu(S) \, = \, \sum_{\ell =1}^\infty \sum_{\begin{subarray}{c}  S \in \sausageset_u :  \\
                            d(S) = \ell \end{subarray}} \mu(S) \\
    &
   = &  \mu^{\otimes \N} \bigg( \exists \, \ell \in \N:  \sum_{i=1}^\ell\ {\rm hgt} (S_i) = u \bigg) \, = \, \mu^{\otimes \N} \big( u \in \rlevel \big) \, ,
\end{eqnarray*}
where the first equality is due to Lemma~\ref{l.essn}.

{\bf (2).} Similarly, $\exp \big\{ \zeta u \big\} \PP \big( \sting_u , \rlevel(\stingset_u)  \in A \big)$ equals 
$$
\mu^{\otimes \N} \bigg( \exists \, \ell \in \N_+ :  \sum_{i=1}^\ell h_i  = u \, , \, \big( r_i: i \in \intint{\ell-1} \big) \in A \bigg) \, ,
$$ 
where here we write $h_i = {\rm hgt} (S_i)$ and $r_i = \sum_{j=1}^i h_i$ under $\mu^{\otimes \N}$;
whence the result. 

{\bf (3).} 
This result arises from the preceding part by taking $A$ to equal the collection of all subsets of $(0,u)$ of cardinality $k-2$. 

{\bf (4).} This is due to the non-lattice renewal theorem~\cite[Chapter~$XI$]{Feller}.

{\bf (5).}
Recall from Definition~\ref{d.enotation} that  $\sting_{[u,u+r)}$ equals $\cup_{h \in [u,u+r)} \sting_h$. Lemma~\ref{l.geometrybasics}(2)
implies that at most one of the events $\sting_u$ in this union may occur; thus, $\PP \big( \sting_{[u,u+r]} \big)$ equals the sum of $\PP \big( \sting_h \big)$  over those $h \in [u,u+r)$
for which this probability is non-zero. From Lemma~\ref{l.zetarenewal}(1), we find then that
$$
  e^{\zeta u} \PP \big( \sting_{[u,u+r)} \big) \, = \,  \mu^{\otimes \N} [u,u+r) \, \sum_h e^{-\zeta(h-u)} \tfrac{\mu^{\otimes \N} (h \in \rlevel)}{\mu^{\otimes \N}[u,u+r)} \, ,
$$
where the sum is over those $h \in [u,u+r)$ for which the summand is non-zero.   From Lemma~\ref{l.zetarenewal}(4) follows the claimed convergence for   $e^{\zeta u} \PP \big( \sting_{[u,u+r)} \big)$.
To finish the proof, note that
$$
e^{\zeta u} \PP \big( \sting_{[u,\infty)} \big) = \sum_{j=0}^\infty e^{-\zeta j r} \cdot e^{\zeta (u+jr)} \PP \big( \sting_{ [u+jr,u+(j+1)r)}\big) \to \mu_\infty r \zeta^{-1} \big( 1 - e^{-\zeta r}\big) \sum_{j=0}^\infty e^{-\zeta j r} = \mu_\infty r \zeta^{-1} \, .
$$
This completes the proof of Lemma~\ref{l.zetarenewal}(5). \qed

\begin{proposition}\label{p.finitemean}
The mean value of the head height  $\head(S) \cdot \vecell$ under the sausage law~$\mu$ is finite. 
\end{proposition}
{\bf Proof.} 
Suppose instead that $\E_\mu \, \head(S)  \cdot \vecell = \infty$. Let $\delta > 0$ be given. We may choose $k \in \N$ so that
$\E_\mu \big[ \min \{ k, \head(S)   \cdot \vecell  \} \big] \geq 2/\delta$.

  With $S$ denoting a sample of $\mu^{\otimes \N}$, we write $S_i$ for the $i$\textsuperscript{th} sausage in the decomposition of~$S$, translated in order that the base of $S_i$ is the origin.  Set
$$
X_j = X_j(S) : =  \sum_{i=1}^j  \min \{ k , \head(S_i) \cdot \vecell \}  \, .
$$
Note then that, whatever the choice of $k \in \N$,
\begin{equation}\label{e.mux}
\mu^{\otimes \N} \Big(  \big\vert \rlevel \cap [0,u] \big\vert \geq \delta u \Big) \leq 
\mu^{\otimes \N}  \big( X_{\delta u} \leq u \big) \, .
\end{equation}
Note however that, for $j \in \N$, $X_j$ is a sum of $j$ independent and identically distributed random variables whose shared mean is at least $2/\delta$. By the exponential Markov inequality,
there exists $\e > 0$ such that, for $u \in \N$,
\begin{equation}\label{e.muxbound}
\mu^{\otimes \N}  \big( X_{\delta u} \leq u \big) \leq \exp \big\{ - \e u \big\} \, .
\end{equation}
By Lemma~\ref{l.zetarenewal}(3),
$$ 
   \exp \big\{ \zeta u \big\}  \PP \big( \sting_u , \big\vert R(\stingset_u) \big\vert \geq \delta u +2 \big) =  
\mu^{\otimes \N} \Big(  \big\vert \rlevel \cap (0,u) \big\vert \geq \delta u  \,  , \,  u \in \rlevel \Big)  \, .
$$
Recalling that   $\sting_{[u,u+r)} = \cup_{h \in [u,u+r)} \sting_h$, we now sum over those heights $u$ lying in a given range $[h,h+r)$, and relabel $h$ as $u$, to obtain
$$
  \exp \big\{ \zeta u \big\}  \PP \big( \sting_{[u,u+r)} , \big\vert R(\stingset) \big\vert \geq \delta u +2 \big) \leq  
\mu^{\otimes \N} \Big(  \big\vert \rlevel \cap (0,u) \big\vert \geq \delta u  \,  , \,  [u,u+r) \cap  \rlevel \not= \emptyset \Big)  \, .
$$
From this display alongside~(\ref{e.mux}) and~(\ref{e.muxbound}), we find that
$$
  \PP \big( \sting_{[u,u+r)} , \big\vert R(\stingset) \big\vert  \geq \delta u + 2 \big) \leq \exp \big\{ -(\zeta +\e)u \big\} 
$$ 
when $u$ is high enough.  Since $\delta u +2 \leq 2\delta u$ for high $u$, we may select $\delta$
to be one-half of the value of this constant in
Proposition~\ref{p.stingrenewal}
and learn from this result  that
$$
  \PP \big( \sting_{[u,u+r)} \big) \leq  c^{-1} \exp \big\{ -(\zeta +\e)u \big\} \, .
$$ 
This conclusion is in contradiction with Proposition~\ref{p.zeta}.
Thus, $\E_\mu \, \head(S) \cdot \vecell < \infty$, as we sought to show in demonstrating Proposition~\ref{p.finitemean}. \qed 
 
 We are ready to prove the weaker counterpart of Theorem~\ref{t.stingwelltrap.sting} in the lattice case.
\begin{proposition}\label{p.hnprob}
Suppose that $\vecell \in S^{d-1}$ is lattice.
As $u \to \infty$ through the cover of $\vecell$, the limit $\mu^{\otimes \N}  \big( u \in \rlevel \big)$ exists and lies in $(0,1)$. Denoting it by $g$, we have that
 $$
   \PP \big( \sting_u \big) \, \simvecell \, g     \exp \big\{ -\zeta u \big\} \, .
 $$
\end{proposition}
{\bf Proof.} The height $\head(S) \cvl$ under $\mu$ has finite mean by Proposition~\ref{p.finitemean}, and the support of the distribution of $\head(S) \cvl$ under $\mu$ is readily seen in the lattice case to have greatest common denominator equal to the lattice span $s$. 
The renewal theorem~\cite[Chapter~$XI$]{Feller} thus implies the existence of the concerned limit.  The latter statement of Proposition~\ref{p.hnprob} then follows from Lemma~\ref{l.zetarenewal}(1). \qed

 \begin{lemma}\label{l.renewalrate}
 Let $\vecell \in S^{d-1}$ be lattice with span $s$.
 There exist positive choices for $C$ and $c$ such that, for $u \in s\cdot \N$ (and $g$ specified by Proposition~\ref{p.hnprob}),
 $\big\vert \mu^{\otimes \N}  \big( u \in \rlevel \big) - g \big\vert \leq C \exp \big\{ - c u^{1/2} \big\}$.
 \end{lemma}
 The next definition and result will facilitate the proof of the one just stated.
\begin{definition}\label{d.maxgap} 
Let $G$ be a finite subgraph of $(\Z^d,\sim)$. Let $h_{\rm min}$ and $h_{\rm max}$ denote the minimum and maximum heights of elements of $V(G)$.
Consider the sequence that begins with $h_{\rm min}$, ends with $h_{\rm max}$, and, in between, lists in increasing order the heights of renewal points of $G$ that lie strictly between $h_{\rm min}$ and $h_{\rm max}$. Let 
 $\maxrenewalgap(G)$ denote the maximum difference of consecutive terms in this sequence. 
 
 For use a little later, take $k \in \N$ and let 
 $\maxrenewalgap(G,k)$ denote the maximum difference among the final $k$ pairs of consecutive terms, where we take $\maxrenewalgap(G,k) = \maxrenewalgap(G)$ in the case that there are at most $k$ terms. 
 \end{definition}
  \begin{lemma}\label{l.clustersting}
 Let $\vecell \in S^{d-1}$ be lattice with span $s$, and 
 let $u \in s \cdot \N$.
 \begin{enumerate}
 \item We have that 
 $$
 \PP \big(\rlevel(S) \cap (0,u)= \emptyset , \sting_u \big)
   \leq  (1-p)^{-d}  \, \PP \big(\rlevel(\cset_u) \cap (0,u)= \emptyset , \cluster_u \big) \, .
   $$
 \item And that 
 $$
 \PP \big( \cluster_{(u,u+r]} , \maxrenewalgap(\mc{C} ) \geq \ell   \big) \leq   p^{-2}(1-p)^{5-4d}  \PP \big( \sting_{(u+2r,u+3r]} , \maxrenewalgap(S) \geq \ell \big) \, ,
 $$  
 where recall that the event $\cluster_{(u,u+r]}$  is specified
 by Definition~\ref{d.enotation} and $\mc{C}$ denotes the cluster containing~$0$. 
 \item  And also that $\PP \big( \sting_{[u,u+r)}  \big)  \leq  p^{-1} (1-p)^{1-2d} \PP \big( \sting_{[u+r,u+2r)}  \big)$.
 \end{enumerate}
 \end{lemma}
 {\bf Proof: (1).} Set ${\bf h} = \head(\stingset_u)$ and $h = {\bf h} \cvl$. When $\sting_u$ occurs, and edges incident to ${\bf h}$  and to an element of $V(U_h) \setminus \vertbdry_h$
 are closed, then $\cluster_u$ occurs, with $\cset_u$ equal to $\stingset_u$. 
 The event $\sting_u$ is independent of the status of the just mentioned edges, of which there are at most $d$, so that the desired bound results.
 
{\bf (2).} The probability $\PP \big( \cluster_{(u,u+r]}  , \maxrenewalgap(\mc{C} \big)$ may be expressed as a sum of terms $p_D$, the probability of the event $E_D$ that the cluster of the origin equals $D$, 
as $D$ ranges over all connected components whose base equals $0$, whose head has height in $[u,u+r]$ and for which $\maxrenewalgap(\mc{C}) \geq \ell$. 
The occurrence of $E_D$ is characterized by the edges in $D$ being open, and the neighbouring edges being closed. Consider instead the event $F_D$ which is specified by the same set of open-and-closed conditions as is $E_D$, with the exception that we instead demand that
the edges $e := (-e_d,0)$ and 
$e' := \big( \head(D),\head(D) + e_d \big)$ be open, and edges that lie in $\slab_{-e_d \cvl , (\head(D) + e_d) \cvl}$ but not in $D$ and that border either of these two edges be closed.
We have that $\PP (F_D) = \alpha \beta \PP(E_D)$ where the factor $\alpha = p^2(1-p)^{-2}$ addresses the switch in status for $e$ and $e'$ and $\beta = (1-p)^\sigma$, with $\sigma = \sigma(D) \in \llbracket 0 , 4d-3 \rrbracket$, represents the closure of the further edges under $F_D$.
When $F_D$ occurs, the subgraph with  edge-set  $D \cup \{ e \} \cup \{ e' \}$ is a string with base $-e_d$
and head $\head(D) + e_d$. Indeed, $\head(D) + e_d$ has height at least $r$ more than any other element of this subgraph, while $-e_d$ has height less than that of any such element by at least $r$; thus, Lemma~\ref{l.geometrybasics}(1)
shows that these two vertices are renewal points in the subgraph and demonstrates that the subgraph is indeed a string. 

Let $F_D^\uparrow$ denote the equiprobable event that $F_D$ occurs after the configuration is shifted by $-e_d$. This collection of events, indexed by the above set of $D$, is disjoint, with union lying in
$\sting_{(u+2r,u+3r]} \cap \big\{ \maxrenewalgap(S) \geq \ell \big\}$. In summary,  
\begin{eqnarray*}
\PP \big(  \cluster_{(u,u+r]} , \maxrenewalgap(\mc{C}) \geq \ell \big) & = & \sum_D \PP ( E_D ) \, \leq \, p^{-2}(1-p)^{5-4d} \sum_D \PP ( F^\uparrow_D )\\
 &  \leq & p^{-2}(1-p)^{5-4d} \PP \big( \sting_{(u+2r,u+3r]} , \maxrenewalgap(S) \geq \ell \big) \, ,
\end{eqnarray*}
so that  Lemma~\ref{l.clustersting}(2) is obtained. 

{\bf (3).} This proof operates similarly to the preceding: open the $e_d$-directed edge bordering the string~$S$ that realizes $\sting_{(u,u+r]}$, and then close the neighbours of the edge's other endpoint that lie in $\slab_{0,(\head(S) + e_d) \cvl}$. \qed

 {\bf Proof of Lemma~\ref{l.renewalrate}.} We will argue that, for  some positive constants $C$ and $c$,
 \begin{equation}\label{e.muargue}
  \mu \big( {\rm hgt}(S) \geq u \big) \leq C \exp \big\{ - c u^{1/2} \big\} \, \, \,  \textrm{for $n \in s\N$} \, .
 \end{equation} 
  Indeed, admitting this, 
   we may invoke~\cite[Theorem~$1$]{StoneWainger} with the function $M(x)$ set equal to the expression $\exp \big\{ - c/2 \cdot u^{1/2} \big\}$, and relabel $c > 0$, in order to obtain  Lemma~\ref{l.renewalrate}.

   To derive~(\ref{e.muargue}), note that
   \begin{eqnarray*}
    \mu \big( {\rm hgt}(S) \geq u \big) & = & \sum_{v \geq u} \mu^{\otimes \N} \big( \rlevel \cap (0,v) = \emptyset \, , \, v \in \rlevel \big) \\
     & = & \sum_{v \geq u} \PP \Big( \rlevel(S) \cap (0,v) = \emptyset \Big\vert \sting_v \Big) \PP (v \in \rlevel) \, ,
   \end{eqnarray*}
  where $v$ is understood to run values for which the summand at stake is non-zero. Here we used that  $\mu^{\otimes \N} \big( \rlevel \cap (0,v) = \emptyset \, \big\vert \, v \in \rlevel \big) = \PP \big( \rlevel(S) \cap (0,v) = \emptyset \big\vert \sting_v \big)$, a fact due to $0$ being a renewal level in any of the concatenated sausage sequences sampled under~$\mu^{\otimes \N}$.
   To obtain~(\ref{e.muargue}), it is thus enough to verify that
   \begin{equation}\label{e.rs}
   \PP \Big( R(S) \cap (0,u) = \emptyset \, \Big\vert \, \sting_u \Big) \, \leq \, C \exp \big\{ - c u^{1/2} \big\} \, .
\end{equation}
 To this end, note that, by considering the possibility that the edges not in $\slab_{0,u}$ that are incident to either the base or the head of the string realizing $\sting_u$ are closed, we obtain the first bound as we write
 \begin{eqnarray*}
   \PP \Big( R(S) \cap (0,u) = \emptyset \, , \,  \sting_u \Big) & \leq & (1-p)^{-2d} \PP \Big( R(\mc{C}) = \{ 0,u\} \, , \, \cluster_u \Big) \\
   & \leq &  (1-p)^{-2d} \PP \Big( R(\mc{C}) = \{ 0,h(\mc{C})\} \, , \, \cluster_{[u-3r,u]} \Big) \\
   & \leq & C \exp \big\{ - c u^{1/2} \big\} \PP \Big(  \{ 0,h(\mc{C})\} \subseteq R(\mc{C}) \, , \, \cluster_{[u-3r,u]}  \Big) \, ,
 \end{eqnarray*}
where we write $h(\mc{C})$ for the height of $\head(\mc{C})$, with $\mc{C}$ the cluster of the origin specified when $\cluster_{[u-3r,u]}$ occurs; and where the latter bound  is due to Corollary~\ref{c.renewal}(2).

\begin{figure}[htbp]
\centering
\includegraphics[width=0.35\textwidth]{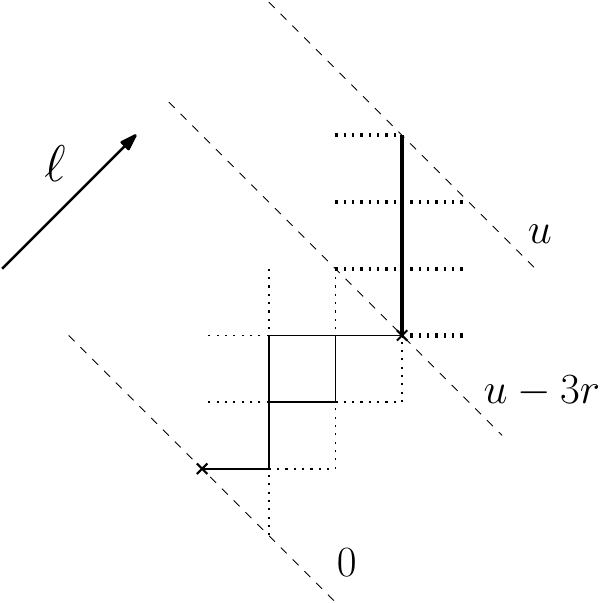}
\caption{Illustrating the bound in~(\ref{e.finiteenergybound}). 
In this two dimensional example, $\vecell$ makes a forty-five degree angle with the axes, and the quantities $s$ and $r$ in (\ref{e.finiteenergybound}) equal $2^{-1/.2}$.
The two crosses mark the base and head of a string that realizes $\sting_{u-3r}$, with closed edges bordering the string shown with dotted lines. 
This string may be extended to form a string realizing $\sting_u$: three consecutive vertical edges emanating from the head are marked in bold, with bordering edges shown bold dotted; if the former are realized open, and the latter closed, then the desired extension occurs, and it does so with positive probability.}\label{f.extend}
\end{figure}

By restriction of the configuration to a slab, we see that
$$
 \PP \Big(  \{ 0,h(\mc{C})\} \subseteq R(\mc{C}) \, , \, \cluster_{[u-3r,u]}  \Big) \leq \PP \Big( \bigcup_{v \in [u-3r,u]} \sting_v \Big) \, .
$$
Now, 
\begin{equation}\label{e.finiteenergybound}
\PP \Big( \bigcup_{v \in [u-3r,u]} \sting_v \Big) \, = \, \sum_{
\substack{
v \in s \N \\
v \in [u-3r,u]
}} \PP (  \sting_v )
\, \leq \, C \, \PP (\sting_u) \, .
\end{equation}
The equality invokes that $\vecell$ is lattice with span $s$.
The bound, illustrated in Figure~\ref{f.extend}, arises by considering the configuration outside $\slab_{0,v}$ and the event that an open path delimited by closed edges reaches to height $u$.
Since this path has finite length, the constant $C$, while dependent on $p$, is finite.

Assembling these bounds, we learn that 
$\PP \big( R(S) \cap (0,u) = \emptyset \, , \,  \sting_u \big) \leq C e^{-c u^{1/2}} \PP (\sting_u)$ for a $p$-dependent constant $C$.
From this conclusion,~(\ref{e.rs}) naturally follows; in this way, we obtain Lemma~\ref{l.renewalrate}. \qed

{\bf Proof of Theorem~\ref{t.stingwelltrap.sting}.} 
The non-lattice case has already been obtained, with 
$$
\kappasting = \zeta^{-1} r \lim_{u \to \infty} \mu^{\otimes \N}[u,u+r) \, :
$$
 it is Lemma~\ref{l.zetarenewal}(5). In the lattice case, the conclusion is due to Lemmas~\ref{l.zetarenewal}(1) and~\ref{l.renewalrate}, with $\kappasting = \lim_{u \to \infty} \mu^{\otimes \N}( \{u\} )$. \qed 

\subsection{Further results on the renewal structure of strings}\label{s.stingconseq}

In separate subsections, we state and prove several results that address the renewal structure of strings and that will be applied in the companion paper.
And then we prove Proposition~\ref{p.zeta}.

\subsubsection{Strings  with a small constant fraction of renewal points are rare} Our first result is an assertion to this effect.
\begin{proposition}\label{p.stingfewrenewal}
We have that
\begin{equation}\label{e.stingfewrenewal}
\PP \bigg(  \bigcup_{w \in [u,u+r)}  \sting_w \cap \big\{ \vert R(\stingset_w) \vert \leq cw \big\} \bigg) \, \leq \, C \exp \big\{ - \zeta u -c u^{1/2} \big\} \, .
\end{equation}
\end{proposition}
Just as after Theorem~\ref{t.renewal} we permitted generally $c > 0$ to depend on $d \geq 3$ and $p \in (p_c,1)$, now we point out that in Proposition~\ref{p.stingfewrenewal} and later, we also permit this dependence for $C > 0$.

For the proof of the proposition,  let $Z_i$ be specified under the law 
$\mu^{\otimes \N}$ to be the difference between the $(i+1)$\textsuperscript{st} and the $i$\textsuperscript{th} elements of  $\rlevel$ for $i \in \N$.
\begin{lemma}\label{l.z}
The sequence $\big\{Z_i : i \in \N \big\}$ is independent and identically distributed. For a $p$-dependent positive constant $C$, we have that $\PP (Z_1 \geq k) \leq C \exp \big\{ - c k^{1/2}\big\}$. 
\end{lemma}
{\bf Proof.} The first statement holds by the definition of~$\rlevel$. 
The second holds because
\begin{eqnarray*}
\PP (Z_1 \geq u ) & = & \mu^{\otimes \N} \Big( \, \exists \, v \geq u: v \in\rlevel \, , \,\rlevel \cap (0,u) = \emptyset  \Big) \\
 & = & \sum_{v \geq u} \mu^{\otimes \N} \Big( v \in\rlevel \, , \,\rlevel \cap (0,u) = \emptyset  \Big) \\
 &= & \sum_{v \geq u} e^{\zeta v} \PP \big( \sting_v \, , \, \vert R(S_v) \vert = 2 \big) \\
 & \leq & C \sum_{v \geq u} e^{\zeta v} \PP \big( \cluster_v \, , \,  R(\mc{C})  = \{ 0, v\} \big) \\
 & \leq & C \sum_{v \geq u} e^{\zeta v - c v^{1/2}} \PP \big( \cluster_v \, , \, \{ 0,v \} \subseteq R(\mc{C}) \big) \\
 & \leq & C \sum_{v \geq u} e^{\zeta v - c v^{1/2}} \PP \big( \sting_v  \big) \\
 & \leq & C \exp \big\{ - c u^{1/2} \big\} \, ,
\end{eqnarray*}
where each sum is taken over all $v \geq u$ for which the summand is non-zero. 
The first inequality is due to Lemma~\ref{l.clustersting}(1).
To derive the second, we split the sum over $v \geq u$ by division of $[u,\infty)$ into intervals of length three, and apply 
 Corollary~\ref{c.renewal}(2) to each of the resulting probabilities.
That the restriction of a configuration realizing $\cluster_v$ to $\slab_{0,v}$
realizes $\sting_v$ justifies the third bound. The fourth is due to Theorem~\ref{t.stingwelltrap.sting}. \qed

{\bf Proof of Proposition~\ref{p.stingfewrenewal}.} By Lemma~\ref{l.zetarenewal}(1,3), 
$$
\PP \big( \sting_w , \vert R(\stingset_w) \vert \leq cw  \big) \, \leq \,
  \exp \big\{- \zeta w \big\}  \mu^{\otimes \N} \big( w \in\rlevel \, , \,   \big\vert\rlevel \cap (0,w) \big\vert \leq \lfloor cw \rfloor -2 \big) \, .
$$ 
for $w >0$ with $cw \geq 2$.
The left-hand side in~(\ref{e.stingfewrenewal})
is bounded above by the sum of the left-hand of the preceding display over all those $w \in [u,u+r)$
for which the event $\sting_w$ is non-empty. Thus, the concerned probability is also bounded above by 
the sum of the right-hand terms over such indices, and thus by$$
  \exp \big\{- \zeta u \big\}  \cdot \mu^{\otimes \N} \Big(    \big\vert\rlevel \cap (0,u) \big\vert \leq \lfloor c(u+1) \rfloor -2 \Big) \cdot \kappa 
$$
where $\kappa$ equals the mean value of $\big\vert\rlevel[u,u+r) \big\vert$
under $\mu^{\otimes \N}$ given $\big\vert\rlevel \cap (0,u) \big\vert \leq \lfloor c(u+1) \rfloor -2$.
The value of $\kappa$ may be bounded above by instead conditioning on the maximum value of $\rlevel$ that is at most $u$; since the random variable $Z_1$ is almost surely positive, we see that $\kappa$ is bounded above, uniformly in $u$. Further, we have that
$$
 \mu^{\otimes \N} \Big(    \big\vert\rlevel \cap (0,u) \big\vert \leq \lfloor c(u+1) \rfloor -2 \Big) \, = \, \PP \bigg( \sum_{i=1}^{\lfloor c(u+1) \rfloor -1} Z_i \geq u \bigg) \, .
$$
We find from Lemma~\ref{l.z} and~\cite{Nagaev69} that $\PP \big( \sum_{i=1}^{\lfloor cu \rfloor} Z_i > u \big) \leq C \exp \big\{ - cu^{1/2}  \big\}$ if $c>0$ is chosen to be suitably small. In light of the bounds above, we obtain Proposition~\ref{p.stingfewrenewal}. \qed

\subsubsection{The start and end of a long string}
After presenting some notation, we state and prove Lemma~\ref{l.startend}, which addresses the joint form of the two ends of a long string.

\begin{definition}\label{d.muproduct}
For $k \in \N$, let $\mc{S}_k$ denote the set of strings that contain $k$ sausages. 
Write $\mu^{\otimes k}(S) = \prod_{i=1}^k \mu(S_i)$
for $S \in \mc{S}_k$ of the form $S = S_1 \circ \cdots \circ S_k$.
\end{definition}

For $k,\ell \in \N$, let $S_+ \in \mc{S}_k$ and  $S_- \in \mc{S}_\ell$.
Let $E(u,k,\ell,S_+,S_-)$ denote the event that $\sting_u$ occurs and that there exists a string $S$ such that $\stingset_u$ has the form $S_+ \circ S \circ S_-$: in other words, that $S_u$
contains at least $k+\ell$ sausages, with the initial string of length $k$ being $S_+$
and the final string of length $\ell$ being $S_-$.

\begin{lemma}\label{l.startend}
We have that
$$
\PP \big( E(u,k,\ell,S_+,S_-) \big\vert \sting_u \big) = \mu^{\otimes k}(S_+) \mu^{\otimes \ell}(S_-)  \frac{\mu^{\otimes \N} \big( u - \height(S_+) - \height(S_-) \in \rlevel \big)}{\mu^{\otimes \N} \big( u  \in \rlevel \big)}  \, .
$$
\end{lemma}
{\bf Proof.} Since $E(u,k,\ell,S_+,S_-) \subset \sting_u$ and 
$$
\PP \big( E(u,k,\ell,S_+,S_-)  \big) = \mu^{\otimes k}(S_+)  \PP \big( \sting_{u -  \height(S_+) - \height(S_-)}  \big) \mu^{\otimes \ell}(S_-) \, ,
$$
the result follows from Lemma~\ref{l.zetarenewal}(1). \qed

\subsubsection{The maximum gap between renewals}
Here we prove Proposition~\ref{p.maxrenewalgap}, which concerns the maximum displacement of consecutive renewal levels in a long string or cluster, and which is expressed in the terminology of Definition~\ref{d.maxgap}.
 \begin{proposition}\label{p.maxrenewalgap}
 Let $\ell \in \N$ and $(\mathsf{E},E) \in \big\{ (\sting,S), (\cluster,\mc{C})
  \big\}$. 
  \begin{enumerate}
  \item Let  $u \in (0,\infty)$ and
  suppose  that $\ell \geq C (\log u)^2$. If $\vecell$ is non-lattice, then 
  $$
  \PP \big( \maxrenewalgap( E_{(u,u+r]}) \geq \ell \, , \,  \mathsf{E}_{(u,u+r]} \big) \, \leq \,      C \exp \big\{-\zeta u - c \ell^{1/2 }\big\} \, .
  $$
 If $\vecell$ is lattice, this statement holds when $E_{(u,u+r]}$ is replaced by $E_u$.
 \item Let $k \in \N$ and suppose that  $\ell \geq C (\log k)^2$. If $\vecell$ is non-lattice, then 
  $$
  \PP \big( \maxrenewalgap( E_{(u,u+r]},k) \geq \ell \, , \,  \mathsf{E}_{(u,u+r]} \big) \, \leq \,      C \exp \big\{-\zeta u - c k^{1/2 }\big\} \, .
  $$
 If $\vecell$ is lattice, this statement holds when $E_{(u,u+r]}$ is replaced by $E_u$.
 \end{enumerate}
 \end{proposition}
 {\bf Proof: (1).} Call the  claims the {\em string} and {\em cluster} 
  bounds. 
 The string bound is needed to prove the cluster bound and we derive it first.
 Suppose that $\vecell$ is non-lattice.
 When $\sting_{(u,u+r]}$ occurs, Lemma~\ref{l.geometrybasics}(2) implies that there is a unique string with base $0$ whose height lies in $(u,u+r]$, which we will denote by $S$. Apply Lemma~\ref{l.zetarenewal}(2)
 with $u=h$ and $A$ equal to the collection of finite subsets of $(0,h)$ for which some pair of consecutive terms in $\{0\} \cup A \cup \{ h\}$
differ by at least $\ell$. Sum the resulting conclusion over $h \in (u,u+r]$ to find that
 $$
 \PP \big( \maxrenewalgap(S) \geq \ell , \sting_{(u,u+r]} \big) \, =  \, \exp \{ -\zeta u \} \mu^{\otimes \N} (u,u+r] \cdot \psi_{u,\ell} \, ,
 $$
  where $\psi_{n,\ell}$ is the probability that there exist consecutive elements in $\rlevel \cap [ 0, u +r]$
 of difference at least~$\ell$ under the law $\mu^{\otimes \N} \big( \cdot \big\vert \rlevel \cap ( u, u+r] \not= \emptyset \big)$. Since $\mu^{\otimes \N}(u,u+r]$ is uniformly positive by Lemma~\ref{l.zetarenewal}(3),
 we may apply
 Lemma~\ref{l.z} to find that $\psi_{\ell,u} \leq u \exp \big\{ - c \ell^{1/2} \big\}$. The condition  $\ell \geq C (\log u)^2$ then implies the string bound in the non-lattice case.
 When $\vecell$ is lattice, we instead have 
 $$
 \PP \big( \maxrenewalgap(S) \geq \ell , \sting_u \big) \, =  \, \exp \{ -\zeta u \} \mu^{\otimes \N} \{ u\} \cdot \psi_{u,\ell} \, ,
 $$
where we respecify $\psi_{n,\ell}$  by instead considering the conditional law  $\mu^{\otimes \N} \big( \cdot \big\vert u \in  \rlevel \big)$. In this lattice case, $\mu^{\otimes \N}(u \in\rlevel)$ is uniformly positive by Proposition~\ref{p.hnprob}. The lattice string bound arises with these variations to the non-lattice case.
 
 The cluster bound is due to 
 Lemma~\ref{l.clustersting}(2).
 
 {\bf (2).} This proof is similar to the preceding one. \qed
  
 \subsubsection{Deriving Proposition~\ref{p.zeta}} 
We conclude Section~\ref{s.asp} by deriving  Proposition~\ref{p.zeta} from \cite[Lemma~$8.5$]{FriberghHammond}.
 
{\bf Proof of Proposition~\ref{p.zeta}.} 
In the notation of~\cite{FriberghHammond},
$$
\mc{H}^+(k) = \big\{ x \in \Z^d: x \cvl \geq k \big\} \, \, \, \textrm{and} \, \, \, 
\mc{H}^-(k) = \big\{ x \in \Z^d: x \cvl \leq k \big\} \, .
$$
In \cite{FriberghHammond}, a notion of sausage connection is specified. In our terms, the event
$\big\{ 0 \xleftrightarrow{\text{s.c.}}  \mc{H}^-(-u) \big\}$ is characterized by the existence of a point $b \in \Z^d$ with $- u -r \leq b \cvl \leq -u$
such that there exists a string $S$ with base $b$ and head $0$ with
the set of edges in $S$ that are incident to either $b$ or $0$ equal to 
$\big\{ (b,b+e_1) , (-e_1,0) \big\}$. 

We {\em claim} that
\begin{equation}\label{e.urclaim}
\PP \Big( 0 \xleftrightarrow{\text{s.c.}}  \mc{H}^-(-u) \Big) \leq \PP \big( \sting_{[u,u+r]}(0) \big) \, .
\end{equation}
Indeed, by translating the string that characterizes the occurrence of the sausage connection event by the vector $-b$,
we see that the left-hand side of the above display is at most the probability that there exists a vertex $h \in \Z^d$
with $u \leq h \cvl \leq u+r$ such that there exists a string with base $0$ and head $h$; the latter event is simply $\sting_{[u,u+r]}$, so the claim follows.

Consider the union~$U$ over $b \in \Z^d$ satisfying $-u \leq b \cvl < -u + r$ of the event $S(b,-e_1) \cap O^-(b) \cap O^+(-e_1)$,
where $S(b,-e_1)$
is the event that there exists a string with base $b$ and head $-e_1$; $O^-(b)$ is the event that the edge $(b-e_i,b)$ is open for $i \in \intint{d}$
if and only if $i=1$; and $O^+(-e_1)$ is the event that the edge $(-e_1,-e_1 + e_i)$ is open for $i \in \intint{d}$
if and only if $i=1$.

Note that $U \subseteq \big\{ 0 \xleftrightarrow{\text{s.c.}}  \mc{H}^-(-u) \big\}$. By translating configurations so that the point $b$ is relocated to be the origin, we obtain the bound
$$
 \big( p(1-p)^{d-1} \big)^2 \PP \Big( \sting_{[u-r,u)}(0) \Big) \leq \PP \Big( 0 \xleftrightarrow{\text{s.c.}}  \mc{H}^+(u) \Big) \, .
$$
Using the same translation on~(\ref{e.urclaim}),
we find that
$$
\PP \Big( 0 \xleftrightarrow{\text{s.c.}}  \mc{H}^+(u) \Big) \leq \PP \Big( \sting_{[u-r,u]}(0) \Big) \, .
$$
Since $\sting_{[u-r,u]}(0) \subset \sting_{[u-r,u)}(0) \cup \sting_{[u,u+r)}(0)$, the two preceding displays imply that
$$
\PP \Big( 0 \xleftrightarrow{\text{s.c.}}  \mc{H}^+(u) \Big) \leq \PP \Big( \sting_{[u-r,u]}(0) \Big)
\leq \big( p(1-p)^{d-1} \big)^{-2} \bigg( \PP \Big( 0 \xleftrightarrow{\text{s.c.}}  \mc{H}^+(u) \Big) + \PP \Big( 0 \xleftrightarrow{\text{s.c.}}  \mc{H}^+(u+r) \Big) \bigg) 
\, .
$$
\cite[Lemma~$8.5$]{FriberghHammond} asserts that $\PP \big( 0 \xleftrightarrow{\text{s.c.}}  \mc{H}^+(u) \big) = \exp \big\{ -\zeta u (1+o(1)) \big\}$.
The just derived bounds shows that this asymptotic passes to $\PP \big( \sting_{[u-r,u]}(0) \big)$. 
We invoke Lemma~\ref{l.clustersting}(3) in the final inequality as we write
\begin{eqnarray*}
\PP \Big( \sting_{[u,u+r)}(0) \Big) \leq \PP \Big( \sting_{[u,u+r]}(0) \Big) & \leq & \PP \Big( \sting_{[u,u+r)}(0) \Big)+ \PP \Big( \sting_{[u+r,u+2r)}(0) \Big) \\
 & \leq &  \big( 1 + C \big)\PP \Big( \sting_{[u+r,u+2r)}(0) \Big) \, .
\end{eqnarray*}
From $\PP \big( \sting_{[u,u+r]}(0) \big) = e^{-\zeta u (1+o(1))}$, we learn that 
$\PP \big( \sting_{[u,u+r)}(0) \big)$ also satisfies this bound.
This completes the proof of Proposition~\ref{p.zeta}. \qed


\subsection{Deriving the asymptotic probability of clusters}\label{s.asympcluster}

Here we prove Theorem~\ref{t.stingwelltrap.cluster}. A cluster typically has at least two renewal points. When it does so, it may be split into three pieces at its lowest and highest renewal points. 
The middle piece is a string, which is bookended by a cul-de-sac and an inverted cul-de-sac, each with unique renewal points. First we specify these {\em culs-de-sac}.

{\bf The cul-de-sac.} 
Let $x \in U_0$. A cul-de-sac with entry $x$ is a finite open connected component of  the lower half-space $L_{x \cvl}$ of which $0$ is the base and $x$ is a renewal point. 
We denote this object by $\culdesacset[x]$ and write $\culdesac[x]$ for the event that it exists.
For $u > 0$, we write $\culdesac_u = \cup \, \culdesac[x]$ where the union is over $x \in \Z^d$ with $x \cvl = u$. 
Note that, if $u$ is such that $\culdesac_u$ is non-empty, then there is only one such value of $x$ for which $\culdesac[x]$ occurs.

Now let $x \in L_0$. An inverted cul-de-sac with entry $x$   is a finite open connected component of  the upper half-space $U_{x \cvl}$ of which $0$ is the head and $x$ is a renewal point. 
We write $\icds$ and $\icdsset$ in place of $\cds$ and $\cdsset$.
 
\begin{figure}[htbp]
\centering
\includegraphics[width=0.45\textwidth]{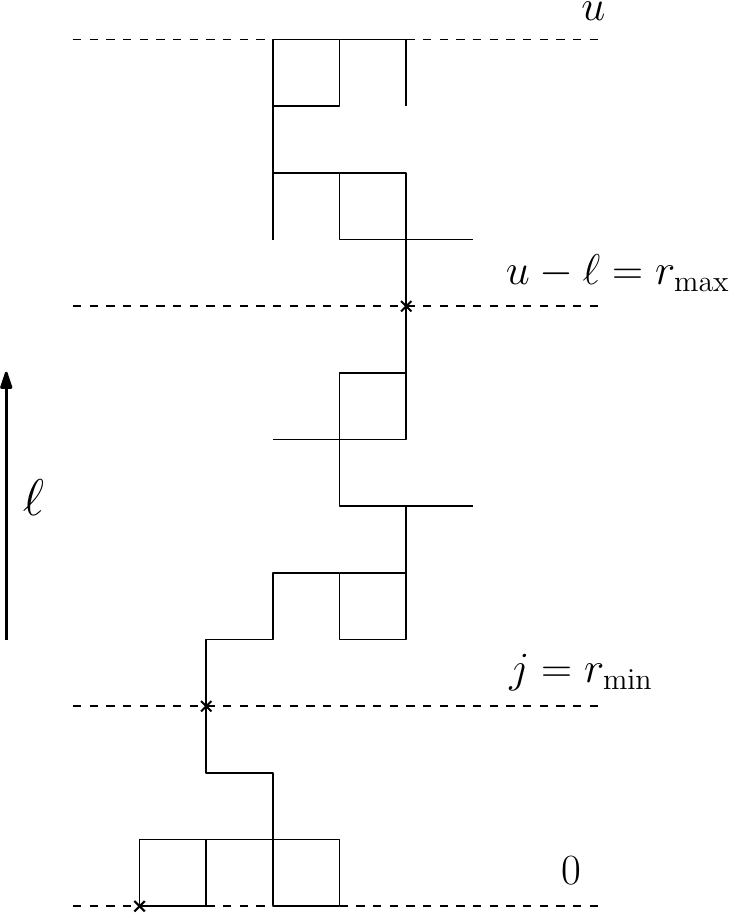}
\caption{Illustrating Lemma~\ref{l.jell.cds}. A cluster splits into an initial cul-de-sac, a middle string and final inverted cul-de-sac.}\label{f.decomposition}
\end{figure}

 The next definition and result speak to the decomposition of a cluster into three pieces.
 \begin{definition}
 On the event $\cluster_u$, with $u \in (0,\infty)$ given, we write $\rmax$ and $\rmin$ for the maximum and minimum heights of renewal points of $\mc{C}_u$.
 \end{definition}
  \begin{lemma}\label{l.jell.cds}
Let $u \in (0,\infty)$ and $j,\ell \in [0 , u ]$ with $j+\ell \leq u$. We have that
  \begin{eqnarray*}
 & &    \PP \Big( \cluster_u \, , \, \rmin = j \, , \, \rmax = u - \ell \Big) \\
 & = & \PP \big( \cds_j, \rlevel(\cdsset_j) = \{ j\} \big) \PP \big( \sting_{u-j-\ell} \big)  \PP \big( \icds_\ell ,\rlevel(\icdsset_\ell) = \{ 0 \} \big) \, . 
 \end{eqnarray*}
 \end{lemma}  
 {\bf Proof.} Consider the percolation configurations restricted to the sets $L_{\rmin}$, $\slab_{\rmin,\rmax}$ and $U_{\rmax}$ in the event 
  $\cluster_u \cap \big\{ \rmin = j ,  \rmax = u - \ell \big\}$. As Figure~\ref{f.decomposition} illustrates, the latter event is characterized by the occurrence of the 
  three right-hand events in the above display, where in the second case, the configuration is translated so that $\head(\cdsset_j)$ appears at the origin, and the third configuration is also shifted, so that 
  the head of the string in the second event is also at the origin. These three events are respectively measurable with respect to percolation configurations restricted to $L_{\rmin}$, $\slab_{\rmin,\rmax}$ and $U_{\rmax}$, so their probabilities multiply, as asserted in the lemma. \qed

Next we prove that unique-renewal culs-de-sac are rare when judged relative to the $e^{-\zeta u}$ decay in Theorems~\ref{t.stingwelltrap.cluster} and~\ref{t.stingwelltrap.sting}.
 \begin{theorem}\label{t.estimates.cds}
 \leavevmode
 \begin{enumerate}
 \item
When $\vecell$ is lattice,
$$
\PP \Big( \culdesac_u \, , \, \rlevel(\culdesacset) = \{ u \} \Big) \leq C \exp \big\{ - \zeta u - c u^{1/2} \big\} 
$$
for $u \in \cover(\vecell)$.
When $\vecell$ is non-lattice, this bound is valid for $u \in [0,\infty)$ when the left-hand expression is taken to be $\sum \PP \Big( \culdesac_h \, , \, \rlevel(\culdesacset) = \{ h \} \Big)$
where here the sum is performed over those $h \in [u,u+r)$ for which the summand is non-zero. 
\item A counterpart statement holds for inverted culs-de-sac. Indeed, the above assertion is equally valid after the changes $\cds \to \icds$ and $\cdsset \to \icdsset$ are made.
\end{enumerate}
\end{theorem}
 {\bf Proof: (1).}
 In the case that $\vecell$ is lattice, we must show that
there exist positive constants $C$ and $c$ such that, for $u \in (0,\infty)$,
$$
 \PP \Big( \culdesac_u ,R(\culdesacset_u) = \big\{ \head(\culdesacset_u) \big\} \Big)
  \leq C \exp \big\{ - \zeta u - c u^{1/2} \big\}  \, .
$$
 To the end of doing so, we first claim that 
 \begin{equation}\label{e.wellcluster}
  \culdesac_u \cap \big\{ R(\culdesacset_u) = \{ \head(\culdesacset_u) \} \big\} \cap \mathsf{Q} \, \subseteq \, \cluster_{u} \cap \{ \mc{C} = \culdesacset_u \} \, ,
 \end{equation}
 where $\mathsf{Q}$ is the event that every element of $U_u$ that is incident to $\head(\culdesacset_u)$ is closed under the percolation, and $\mc{C}$ is the open cluster containing~$0$. The inclusion is evident enough: an open path from $0$ in the cul-de-sac  must leave $L_u$ via the entry $x$, but it has no ongoing route from there if $\mathsf{Q}$ occurs.
 
  If the edge running out of $\head(\culdesacset_u)$ in the direction $e_i$ lies in $E(U_e)$, then the edge running out of this vertex in the direction $-e_i$ does not lie in this set. 
 Thus, the number of elements of $E(U_e)$ that are incident to $\head(\culdesacset_u)$ is at most $d$. We see then that
 \begin{equation}\label{e.qprob}
  \PP \Big( \mathsf{Q}
 \, \Big\vert \,  
\culdesac_u \cap \big\{ R(\culdesacset_u) = \{ \head(\culdesacset_u) \} \big\} \Big) \geq (1-p)^d \, .
   \end{equation}
 As a result of this and~(\ref{e.wellcluster}), we learn that
\begin{equation}\label{e.culdesaccluster}
 (1-p)^d \PP \Big( \culdesac_u , R(\culdesacset_u) = \big\{ \head(\culdesacset_u) \big\} \Big) \leq \PP \Big( \cluster_u(0) , R ( \cset ) = \big\{ \head (\cset) \big\} \Big) \, .
\end{equation}
Suppose that  $\cluster_u(0)$ and $R ( \cset) = \big\{ \head(\cset) \big\}$ occur. Consider a resampling of the edges incident to $-e_d$. Suppose this resampling results in $\{-e_d,0\}$
being open, and the remaining edges being closed. And suppose further that we translate the resulting configuration by $e_d$.
Then in the new configuration the event $\sting_{u+r}$ will occur with $R(S_{u+r}) = \{ 0 , \head(S_{u+r}) \}$.
Thus,
\begin{equation}\label{e.clusterub}
 p (1-p)^{2d-1}  \PP \Big( \cluster_u(0) , R \big( \cset \big) = \big\{ \head (\cset ) \big\} \Big) 
 \leq \PP \Big( \sting_{u+r} \cap  \big\{ \rlevel(S_{u+r}) = \{ 0, u + r \}  \big\}  \Big) \, .
\end{equation}
The right-hand probability is at most  $C \exp \big\{ - \zeta u - c u^{1/2} \big\}$ by Proposition~\ref{p.stingfewrenewal}.
What we have learnt is that
\begin{equation}\label{e.cdsbound}
  \PP \Big( \culdesac_u , R(\culdesacset_u) = \big\{ \head(\culdesacset_u) \big\} \Big) \ \leq p^{-1} (1-p)^{1-3d} C  \exp \big\{ - \zeta u - c u^{1/2} \big\} \, .
\end{equation}
This completes the proof of Theorem~\ref{t.estimates.cds}(1) in the case that $\vecell$ is lattice.

Suppose now that $\vecell$ is non-lattice.  Replace $u \to h$ in~(\ref{e.culdesaccluster}), and denote the right-hand side of the resulting bound summed over $h \in [u,u+r)$ by~$\rho$.  Analogously to~(\ref{e.clusterub}), 
$$
 p (1-p)^{2d-1} \rho \, \leq \, \PP \Big( \sting_{[u,u+r)} \, , \, R(S_{[u,u+r)}) = \{ 0, \head (S_{[u,u+r)}) \}    \Big) \, ,
$$
where $S_{[u,u+r)}$ denotes the string between $0$ and some vertex of height in $[u,u+r)$ whose uniqueness is ensured by Lemma~\ref{l.geometrybasics}(1).
Since the right-hand probability is at most $C \exp \big\{ - \zeta u - c u^{1/2} \big\}$ by Proposition~\ref{p.stingfewrenewal}, assembling the estimates yields that
the sum of 
$$
\PP \big( \culdesac_h , R(\culdesacset_h) = \big\{ \head(\culdesacset_h) \big\} \big)
$$
over those $h \in [u,u+r)$ for which the summand is non-zero is bounded above by the right-hand expression in~(\ref{e.cdsbound}). 
Thus we obtain Theorem~\ref{t.estimates.cds}(1) in the non-lattice case. 

{\bf (2).} Write $-e = (-x,-y)$ and $e \oplus u = (x+u,y+u)$ for $e = (x,y) \in E(\Z^d)$ and $u \in \Z^d$. For $x \in \Z^d$, define a bijection $\omega \to \omega'$ on configurations by $\omega'(e) = \omega(-e \oplus x)$. Let $x \in U_0$. An inverted cul-de-sac with entry $-x$ occurs under $\omega'$ precisely when a cul-de-sac with entry $x$ occurs under $\omega$ and any edge incident to $x$
both of whose endpoints has height $x \cvl$ is $\omega$-closed. (This $\omega$-closed condition arises because the strict and non-strict inequalities in the definitions of forward and backward half-spaces in Section~\ref{s.basics} change places under reflection through the origin.) The ratio of probabilities for an inverted and a standard cul-de-sac thus takes the form $(1-p)^\kappa$ for a $\vecell$-determined value of $\kappa$ lying in $\llbracket 0, 2(d-1) \rrbracket$. In this way, Theorem~\ref{t.estimates.cds}(2) is reduced to the preceding part.
\qed

Now we harness the just proved result to learn that the culs-de-sac that bookend a typical lengthy cluster are rather short.
 \begin{lemma}\label{l.cluster}
\leavevmode
\begin{enumerate}
\item Suppose that $\vecell$ is lattice with span $s$.
There exist positive $C$ and $c$ such that the quantities 
 $$
 \alpha_u =   \frac{\PP \big( \culdesac_u, \rlevel(\culdesacset_u) = \{ u \} \big)}{\PP \big( \sting_u \big)} 
  $$
 and  
 $$
 \beta_u  =  \frac{\PP \big( \icds_u , \rlevel(\icdsset_u) = \{ 0 \} \big)}{\PP \big( \sting_u \big)} 
 $$
 satisfy
 $$
  \max \{ \alpha_u , \beta_u  \} \leq C \exp \big\{ - c u^{1/2} \big\}
   $$
 whenever $u \in s \cdot \N$.
 \item Now suppose that $\vecell \in S^{d-1}$ is general, lattice or non-lattice.
 Write 
 $$
 \alpha_{[u,u+r)} =   \frac{\sum \PP \big( \culdesac_h, \rlevel(\culdesacset_h) = \{ h \} \big)}{\PP \big( \sting_{[u,u+r)} \big)} 
  $$
 and  
 $$
 \beta_{[u,u+r)}  =  \frac{\sum \PP \big( \icds_h , \rlevel(\icdsset_h) = \{ 0 \} \big)}{\PP \big( \sting_{[u,u+r)} \big)} 
 $$
where the sums are taken over those $h \in [u,u+r)$ for which the summand is non-zero.
Then  $\max \{ \alpha_{[u,u+r)}  , \beta_{[u,u+r)}   \} \leq C \exp \big\{ - c u^{1/2} \big\}$ for $u \in (0,\infty)$.

 \end{enumerate}
\end{lemma}
{\bf 
(1).} Theorem~\ref{t.stingwelltrap.sting} and the lattice aspect of Theorem~\ref{t.estimates.cds}(1) imply the bound on $\alpha_u$, while the bound on $\beta_u$
arises when we instead invoke  Theorem~\ref{t.estimates.cds}(2). 

{\bf (2).}
Note that in the non-lattice case Theorem~\ref{t.stingwelltrap.sting} implies that $\PP \big( \sting_{[u,u+r)} \big) \sim \kappasting (1 - e^{-\zeta r} )e^{-\zeta u}$.
 The sought result thus follows from the just given argument where we instead invoke the non-lattice assertions made in the two parts of Theorem~\ref{t.estimates.cds}. \qed

We are ready to give the promised derivation.

{\bf Proof of Theorem~\ref{t.stingwelltrap.cluster}.} 
Suppose first that $\vecell$ is lattice, with span $s \in (0,\infty)$.
Set  $\gamma_i = \exp \{ \zeta i \} \PP \big( \sting_i \big)$ for $i \in s \cdot \N$. Theorem~\ref{t.stingwelltrap.sting} implies the existence of 
 $\gamma_\infty \in (0,1)$ such that the strong asymptotic $\gamma_i \thickapproxvecell \gamma_\infty$ holds (as $i \to \infty$)  in the sense of Definition~\ref{d.asymptotic}.

For $u \in s \cdot \N$, and in the notation of Lemma~\ref{l.cluster}, set $\phi_u = \sum_{j,\ell \in [0,u]}  \alpha_j \gamma_j \beta_\ell \gamma_\ell \gamma_{u-j-\ell} {\bf 1}_{j+\ell \leq u}$ for $u \in s \cdot \N$; here, and throughout the proof, summations over $s \cdot \N$ are understood.
Also write $\phi_\infty = \gamma_\infty  \sum_{j \geq 0}  \alpha_j \gamma_j   \sum_{\ell \geq 0} \beta_\ell \gamma_\ell$.

We claim that
\begin{equation}\label{e.phiu.cds}
\PP \big( \cluster_u ,\rlevel(C_u) \not= \emptyset  \big) =   \phi_u \, .
\end{equation}
Indeed, this follows from
  \begin{eqnarray}
   & & \PP \big( \cluster_u ,\rlevel(C_u) \not= \emptyset  \big) \nonumber  \\
   & = &  \sum_{j,\ell \in [0,u]} \PP \Big( \cluster_u \, , \, \rmin = j \, , \, \rmax = u - \ell \Big) {\bf 1}_{j+\ell \leq u}   
   \nonumber \\
 & = &  \sum_{j,\ell \in [0,u]}  \PP \big( \cds_j, \rlevel(\cdsset_j) = \{ j\} \big) \PP \big( \sting_{u-j-\ell} \big) \nonumber \\
 & & \qquad \qquad \qquad \times \, \PP \big( \icds_\ell , \rlevel(\icdsset_\ell) = \{ 0 \} \big)  {\bf 1}_{j+\ell \leq u}\nonumber  \\
 & = &   \sum_{j,\ell \in [0,u]}  \alpha_j \beta_\ell \PP \big( \sting_j \big)  \PP \big( \sting_{u-j-\ell} \big)  \PP \big( \sting_\ell \big)   {\bf 1}_{j+\ell \leq u} \, , \nonumber
 \end{eqnarray}
 where the second equality is due to Lemma~\ref{l.jell.cds}.

We next argue in favour of the strong asymptotic relation $\phi_u \thickapproxvecell \phi_\infty$. To derive this, note that $\phi_\infty - \phi_u$ equals $A_1 + A_2$, where
$$
 A_1 =  \sum_{j,\ell \geq 0}  \alpha_j \gamma_j \beta_\ell \gamma_\ell \big( \gamma_\infty - \gamma_{u-j-\ell} \big)  {\bf 1}_{j+\ell \leq u} 
$$
and
$$
 A_2 = \gamma_\infty \sum_{j,\ell \geq 0}  \alpha_j \gamma_j \beta_\ell \gamma_\ell  {\bf 1}_{j+\ell > u}  \, .
$$
Note here that $A_2$, being at most $\gamma_\infty \sum_{j \geq u/2} \alpha_j \gamma_j \sum_{\ell \geq u/2} \beta_\ell \gamma_\ell$, satisfies $A_2 \thickapproxvecell 0$ in view of Lemma~\ref{l.cluster}(1).
We also have that $\vert A_1 \vert \leq A_{11} + A_{12} + A_{13}$, where the latter expressions are sums of 
$$
 \alpha_j \gamma_j \beta_\ell \gamma_\ell \cdot \big\vert \gamma_\infty - \gamma_{u-j-\ell} \big\vert \cdot {\bf 1}_{j+\ell \leq u} 
$$
over the respective index sets $0 \leq j,\ell \leq u/4$; $j \geq u/4$ with $\ell \geq 0$; and $\ell \geq u/4$ with $j \geq 0$.
Since $\sup \gamma_i$ is bounded, the relations $\beta_\ell \gamma_\ell \thickapproxvecell 0$ and $\alpha_j \gamma_j \thickapproxvecell 0$ available from $\gamma_i \thickapproxvecell \gamma_\infty$ and Lemma~\ref{l.cluster}(1) serve to show that $A_{1i} \thickapproxvecell 0$ for $i \in \intint{3}$. We may now return to $\phi_\infty - \phi_u = A_1 + A_2$ and find, as we sought to do, that $\phi_u \thickapproxvecell \phi_\infty$.
From~(\ref{e.phiu.cds}), 
$\PP \big( \cluster_u ,\rlevel(C_u) \not= \emptyset  \big)$ is seen to be strongly asymptotic to $\kappacluster \exp \{ - \zeta u \}$,
 where 
 $$
 \kappacluster \, =  \,   \sum_{j \in s \cdot \N} \alpha_j \gamma_j \, \sum_{\ell \in s \cdot \N}  \beta_\ell \gamma_\ell \, . 
 $$
But 
\begin{equation}\label{e.channelr}
\PP \big( \cluster_u   \big) \thickapproxvecell \PP \big( \cluster_u ,\rlevel(C_u) \not= \emptyset  \big)
\end{equation}
 by Theorem~\ref{t.renewal}(2) and $\PP (\cluster_u ) \geq c \PP (\cluster_{[u,u+3]})$.
 In this way, we obtain Theorem~\ref{t.stingwelltrap.cluster} in the case that $\vecell$ is lattice. 
 
 When $\vecell$ is non-lattice, we seek to show that 
 \begin{equation}\label{e.desasyrel.cds}
 \PP \big( \cluster_{[u,\infty)} \big) \sim \kappa \exp \big\{ - \zeta u \big\} \, ; 
 \end{equation}
 what has changed is that the event of height at least, rather than exactly, $u$ is considered, and the claimed asymptotic relation is $\sim$ rather than the stronger $\thickapprox$.

 Let $\cds[x]$ denote the event that $\cds_{x \cvl}$ occurs, with $x \in V \big( \cdsset_{x \cvl} \big)$; and similarly set $\icds[x] = \icdsset_{x \cvl} \cap \big\{ x \in V \big( \icdsset_{x \cvl} \big) \big\}$.
 On $\cds[x]$, set $\cdsset[x] =  \cdsset_{x \cvl}$ and on $\icds[x]$, set $\icdsset[x] = \icdsset_{x \cvl}$.
 Write $\sting[x]$ for the event that $\sting_{x \cvl}$ occurs alongside $x \in R(\stingset_x)$. 
 We write 
 $$
 \alpha(x) =   \frac{\PP \big( \cds[x],  R(\cdsset[x]) = \{ x\} \big)}{\PP \big( \sting[x] \big)} \, \, \, \, , \, \,  \, \, 
  \beta(y)  =  \frac{\PP \big( \icds[y] , R(\icdsset[y]) = \{ 0 \} \big)}{\PP \big( \sting[y] \big)} 
 $$
and $\gamma(x) = \PP (\sting[x]) \exp \{ \zeta x \cvl \}$.
 We now consider 
  \begin{eqnarray*}
   & & \PP \big( \cluster_{[u,\infty)} ,\rlevel(C_u) \not= \emptyset  \big)  \\
   & = &  \sum_{x,y \in F_0} \PP \Big(  \cluster_{[u,\infty)} \, , \, \renmin = x \, , \, \renmax = y \Big) \cdot {\bf 1}_{(x+y) \cvl \leq u}   \\
 & = &   \sum_{x,y \in F_0}   \PP \big( \cds[x],  R(\icdsset[x]) = \{ x \} \big) \PP \big( \sting_{[u-(x+y) \cvl,\infty)} \big)  \\ 
  & & \qquad \qquad \qquad \times \, \PP \big( \icds[y] , \rlevel(\icdsset[y]) = \{ 0 \} \big) \cdot  {\bf 1}_{(x+y) \cvl \leq u}  
  \\
 & = &    \sum_{x,y \in F_0}   \alpha(x) \beta(y) \PP \big( \sting[x] \big)  \PP \big( \sting_{[u-(x+y) \cvl,\infty)} \big)  \PP \big( \sting[y] \big) \cdot  {\bf 1}_{(x+y) \cvl \leq u}   \, . 
 \end{eqnarray*}
Writing $\psi:(0,\infty) \to (0,\infty)$, $\psi(u) = e^{\zeta u } \PP (\sting_{[u,u+r]})$, we find that 
$$
  \PP \big( \cluster_{[u,\infty)} ,\rlevel(C_u) \not= \emptyset  \big) \exp \{ \zeta u \}
  =   \sum_{x,y \in F_0}   \alpha(x) \beta(y) \gamma(x) \gamma(y) \psi \big( u - (x+y) \cvl \big) {\bf 1}_{(x+y) \cvl \leq u}   \, . 
$$ 
 Lemma~\ref{l.zetarenewal}(5) implies that $\psi(u) \sim \mu_\infty \zeta^{-1} (1 - e^{-\zeta})$ as $u \to \infty$, where $\mu_\infty = \lim_{u \to \infty} \mu^{\N} [u,u+r)$. 
 We thus obtain the desired asymptotic relation~(\ref{e.desasyrel.cds}) with $\kappa = \kappacluster$ given by
 $$
 \kappacluster \, =  \,  \mu_\infty   \zeta^{-1} (1 - e^{-\zeta}) \sum_{x \in F_0} \alpha(x) \gamma(x) \, \sum_{y \in F_0}  \beta(y) \gamma(y) \, . 
 $$
 Thus we obtain Theorem~\ref{t.stingwelltrap.cluster} in the non-lattice case. 
 
 {\bf (2).} For $x \in \Z^d$, consider the operation that acts on a percolation configuration by reflecting it through the origin.
The operation leads the percolation law invariant, and, neglecting a certain priviso, it maps
a cul-de-sac with entry $x$ to an inverted cul-de-sac with entry~$-x$.
 
 \qed

\section{Renewal abundance heuristics and basic tools}\label{s.heuristics}

It remains to prove our key technical contribution  on renewal abundance, Theorem~\ref{t.renewal}. We will do so in the next section. Here, we explain in a first subsection 
some of the main ideas and in a second present some basic tools involving coarse-grained percolation needed to implement them.

\subsection{Slide into renewal}\label{s.slideintorenewal}

How to prove that renewal points are abundant in a typical cluster of high diameter? 
Here we discuss an approach to this question; in so doing, we offer a heuristic overview of the proof of our result, Theorem~\ref{t.renewal}(1), on the typical positive density of renewal in a lengthy finite cluster.
In this presentation, we take $\vecell =e_d$ to be vertical, and consider the typical form of the cluster $\mc{C}$ under percolation $\PP$
given the event $\cluster_u$ that a finite cluster based at the origin runs to height equal to $u \in \N$. We will write $\mc{C}_u$ for $\mc{C}$ to emphasise this conditioning.

The cluster $\mc{C}_u$ is a finite open subgraph of vertical extension $u$ that is encased away from the infinite open cluster~$\mc{I}$ by dual plaquettes through whose midpoints run closed edges in the cluster boundary. In a crude but useful distortion, this dual plaquette surface may be viewed as a subcritical object, whose cardinality has an exponential tail under the unconditioned percolation law~$\PP$, because 
the opposite ends of the closed edge running through any constituent plaquette lie inside and outside of the cluster, so that the outside endpoint
is at a significant
(and possibly even infinite) remove from the inside endpoint in the chemical distance along open edges; and this  event is a rarity that may be diagnosed locally. 

\begin{figure}[htbp]
\centering
\includegraphics[width=0.75\textwidth]{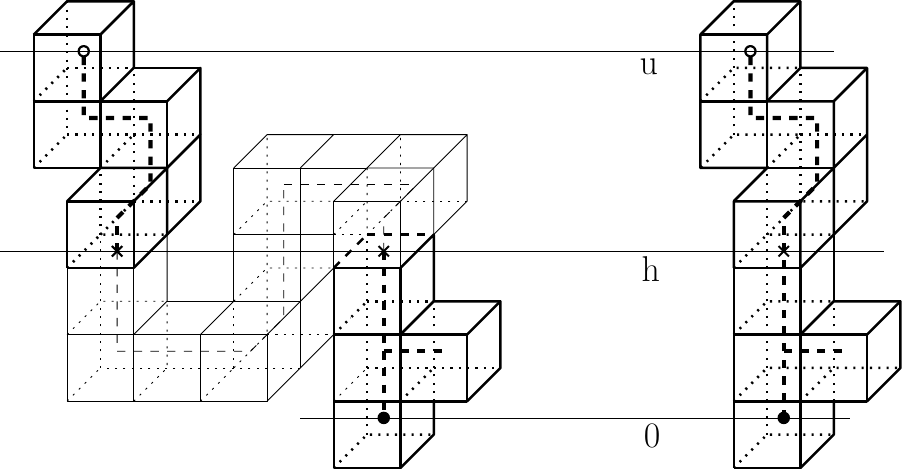}
\caption{A slide eliminates from a capital $N$ its midcluster.}\label{f.Nslide}
\end{figure}

In Proposition~\ref{p.nobigcluster}, we will develop this notion, by using a coarse-grained percolation and invoking results of Antal and Pisztora~\cite{AntalPisztora}. For now, we take the lesson of this result to be that, under $\PP$ given $\cluster_u$, $\mc{C}_u$ typically has volume of the order of the minimum possible,~$u$. The cluster $\mc{C}_u$ runs between heights zero and $u$. Suppose that we slice through the cluster at a typical height~$h$, around level $u/2$. The lesson indicates that the cross-section of $\mc{C}_u$ at height $h$ is likely to have a cardinality that is bounded above independently of $u$.
If this cross-section occurs at a unique vertex, that vertex would be a renewal point for the cluster; so in seeking to derive renewal abundance, we wish to argue that this circumstance is typical under the conditioned law. As it stands, though we know that the cross-section is typically bounded in cardinality, it may be high in diameter. The cluster could resemble a capital `$N$', as depicted in Figure~\ref{f.Nslide}(left).
Elements of $E(\mc{C}_u)$ are dashed, and have the rough shape of a path that runs upwards from the base, doubling back across level $h \approx u/2$, and then running upwards again to the head. Slicing through at height $h$ decomposes the cluster into three pieces: the downcluster and upcluster, emboldened in the depiction, which are the respective open connected components of the cluster's base and head within the half-spaces delimited by the slice; and a remaining piece, the midcluster, which runs between the two vertices marked with crosses that lie one apiece in the downcluster and the upcluster.

The midcluster is a costly item in the example, because it is sheathed by a subcritical surface  of dual plaquettes away from the two surfaces associated to the upcluster and the downcluster.
Not only expensive but also, it would seem, needless: the downcluster covers a height of $h$; the upcluster, one of $u-h$; so together, all the distance needed to realize the event $\mc{C}_u$ is covered. How then to gauge the resulting rarity of this picture? It is tempting to shift the upcluster horizontally so that its base coincides with the head of the downcluster (which is to say, so that the two crosses coincide). Figure~\ref{f.Nslide}(right) shows the shifted picture. There is no expensive surface enclosing a midcluster and so the right picture seems of the two to be much more probable. 
A probability resampling experiment that retains enough information to quantity this intuition works as follows. Sample the law $\PP$ given $\cluster_u$. Record the downcluster; and also record the upcluster, but only up to a horizontal shift, so that only the equivalence class of the upcluster under identification by translation in $\Z^d$ is retained. Then work to reconstruct the law $\PP$ given $\cluster_u$ and also given the two pieces of retained data. In the depicted case, this reconstruction job involves two tetris-like pieces, the lower one held in place, and the upper one resting on the level-$h$ hyperplane, and capable of gliding over it, with any given spot being chosen with
what we may call the {\em connection chance}, namely with
 a weight proportional to the probability that, were the rest of space to be completed according to the ordinary percolation law, a suitable (and perhaps trivial) finite midcluster would form to connect the two pieces and reconstitute a copy of the cluster~$\mc{C}_u$.  

When the upcluster is shifted so that the two crosses coincide, there is unit connection chance; whereas the connection chance for distant displacements may be expected to be of rapid decay.

These ideas may show some promise, but there are objections. One  problem is that, if the slice height is such that the downcluster has more than one vertex at the height, no horizontal shift of the upcluster will lead to a renewal point at that height. To overcome this difficulty, we may simply cut the downcluster off from above at height $h$, while cutting the upcluster from below at height $h+2$; then, if we shift the upcluster so that its base hovers two units above the head of the downcluster, a resampling in which the two intervening vertical edges are declared open will, alongside a suitable set of other edges being declared closed, result in a resampled cluster~$\mc{C}_u$ that has a renewal point at height at $h+1$; we say that a {\em simple join} occurs when the resampling acts in this way.

A second difficulty is represented by the `collapsed capital $N$': the midcluster could be lengthy, but if its surface is shared with those that enclose the downcluster and the upcluster, it is hard to argue that the connection chance is small, since the formation of the midcluster surface is not an independent event.

\begin{figure}[htbp]
\centering
\includegraphics[width=0.75\textwidth]{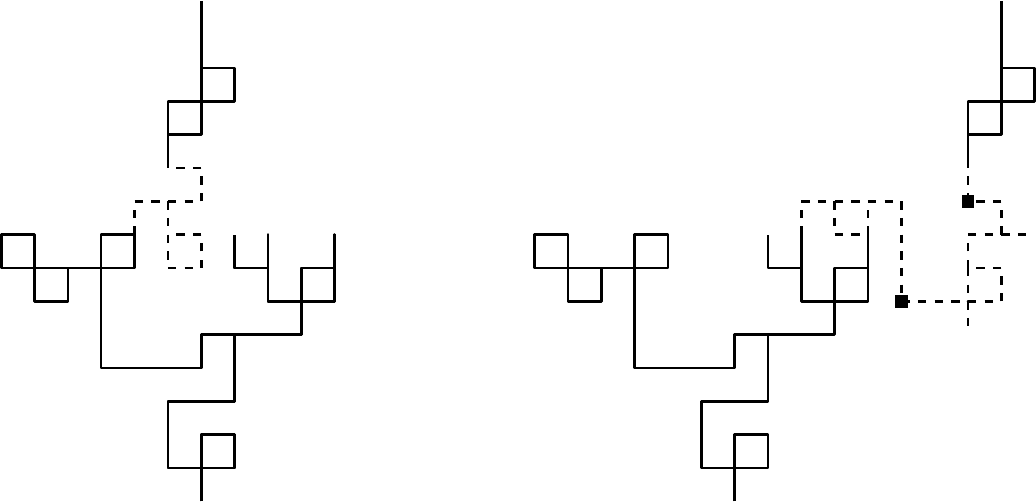}
\caption{The slide resample pushes the upcluster some way to the right in this depiction. The downcluster has some presence near the slide height and this permits the midcluster in the right sketch to travel for a while around the downcluster, but the slide is big enough to force the midcluster to make a journey in virgin territory between the two square vertices. The disjoint surface that encloses this part of the midcluster reflects the low connection chance (and the high connection cost!) of the proposed slide location.}\label{f.slideschematic}
\end{figure}

Figure~\ref{f.slideschematic} shows a slide resample with a gap of two at the slide height. The illustration is schematic, depicting a two-dimensional example for simplicity. 
Simple join would occur were the upcluster on the left resampled by a shift to the left by one unit, with two open vertical midcluster edges running down from the resulting upcluster base, and suitable neighbouring edges being closed.
The resampling on the right shows a dashed midcluster whose surface (which is not depicted)
shares several facets with that of the downcluster. That the downcluster has several routes up to around level~$h$ broadens the ways in which the midcluster and downcluster surfaces may coincide in part. Note however between the two square vertices runs a midcluster section whose surface is disjoint from those of the upcluster and the downcluster. In general, we may think of drawing a path in the midcluster that connects the downcluster head to the upcluster base and selecting the subpath---call it {\rm Connect}!---from the last visit to a neighbour of the downcluster until the first visit a neighbour of the upcluster.  {\rm Connect} is plausibly enclosed by a disjoint surface and roughly speaking the connection chance for a given proposed upcluster shift is of the order of $e^{-c\kappa}$ where $\kappa$ is the minimum length of ${\rm Connect}$ as we vary over viable choices of midcluster given the shift. 

If the downcluster and upcluster cardinalities are bounded above suitably in a neighbourhood of the slice, then the number of shift locations for which the $\kappa$-value is bounded above is limited. The $\kappa$-value grows linearly in the distance away from such special sites. This means that a given shift that may lead to simple join and a renewal point (such as one unit to the left in Figure~\ref{f.slideschematic}(left))
has a connection chance that is bounded below, while the sum over all shifts of the connection chance is bounded above, with a bounded population of sites making a unit-order contribution and others making one that is exponentially small as a function of distance between downcluster and shifted upcluster. 

Our plan is to present the fairly weak (and in essence already known) estimate that we have mentioned, namely Proposition~\ref{p.nobigcluster}, to the effect that cluster volume under $\PP$ given $\cluster_u$ is typically linear in $u$,
and to harness this to argue that downcluster and upcluster cardinalities at heights~$h'$ close to the slide height~$h$ typically grow at a suitably low rate as $\vert h' - h \vert$ rises. Then analysing a partition function given by the sum of connection chances will serve to show that the slide resample results with positive probability  in renewal near the slide height. Since the slide resample leaves the conditioned law at equilibrium, we will thus prove Theorem~\ref{t.renewal}(1), our assertion of positive density of renewal.

\subsection{Coarse-grained percolation}\label{s.coarse}

This subsection is principally devoted to stating and proving Proposition~\ref{p.nobigcluster} on the linear volume of the enclosing surface of a lengthy finite cluster.
We use a coarse-graining of percolation to present the result, dividing $\Z^d$ into boxes.

Two distinct vertices $v,w \in \Z^d$ are {\em $*$-neighbours} if  $w_i - v_i \in \{-1,0,1\}$ for $i \in \intint{d}$. A $*$-path is a $\Z^d$-valued path whose pairs of consecutive elements are $*$-neighbours. 
The $*$-distance of $v,w \in \Z^d$ equals $\ell - 1$, where $\ell$ is the number of elements in the shortest $*$-path that connects $v$ and $w$.

Let $K \in \N$ be a given {\em even} integer. 
For $v \in 2K\cdot \Z^d$, the $K$-box $B(v)$ with central vertex $v$ is $v + \llbracket -K,K \rrbracket^d$.
For $A$ a subgraph of $(\Z^d,\sim)$, $\overline{A}$ will denote the set of $K$-boxes that contain an edge in $A$.

The set $\big\{ B(v): v \in \Z^d \big\}$ of $K$-boxes is naturally indexed by the central vertex $v \in \Z^d$,
so that the notion of $*$-neighbour may be applied to such boxes just as it is to elements of $\Z^d$.

Recall 
from Section~\ref{s.basics} that a path makes nearest neighbour steps.
\begin{definition}\label{d.exterior}
Let $\partial_{\rm ext} \overline{\cset}$ denote the set of $K$ boxes $B$ 
such that $B$ is a $*$-neighbour (in the sense of $K$-boxes) of an element of $\overline{\cset}$
and  there exists an infinite path $Q$ in $K$ boxes emanating from $B$ that is disjoint from  $\overline{\cset}$. 
\end{definition}

Here is the promised result.

\begin{proposition}\label{p.nobigcluster}
For $d \geq 2$, there exist positive constants $C$, $c$ and $D_0$ such that, when $D \geq D_0$,
 $$
 \PP \Big(  \big\vert \partial_{\rm ext} \overline{\cset} \big\vert  \geq Du   \, \Big\vert \, \cluster_{[u,u+1]} \Big) \leq 
 C \exp \big\{ - c D u \big\} \, .
$$
\end{proposition}

\subsubsection{Box percolation}\label{s.box}
As tools we now present will serve to show, the above proposition will be derived by arguing that the $K$-box elements of  $\partial_{\rm ext} \overline{\cset}$ are unusual in a way that can be locally diagnosed.

The {\em boxlet} ${\rm boxlet}(v)$ with central vertex $v$ is $v + \llbracket -K/2,K/2 \rrbracket^d$ and the {\em big box} $\boxbig(v)$ with this central vertex is  $v + \llbracket -5K,5K \rrbracket^d$. 
 We will refer to the $K$-box $B(v)$ as a function of the big box $D = \boxbig(v)$ in the form Centre = Centre($D$).

For $G$ a subgraph of $(\Z^d,\sim)$,
we write $\partial G = \big\{ v \in V(G): \exists \, w \in \Z^d \setminus V(G) \, , \, v \sim w \big\}$.

Let $D = 2K v + \llbracket -5K,5K \rrbracket^d$, $v \in \Z^d$, denote a big box. 
The $K$-box  $B(v)$ will be called {\em bad}  if at least one of two conditions is satisfied:
\begin{enumerate}
\item there exists no open path from ${\rm boxlet}(v)$ to $\partial B(v)$---we say that `${\rm boxlet}(v)$ is isolated'; or
\item there exists an open path $P$ from  ${\rm boxlet}(v)$ to $\partial B(v)$ (so that $P$ lies in $B(v)$ and thus also in $\boxbig(v)$), and an open path $Q$, also contained in $\boxbig(v)$, of Euclidean diameter at least~$K/2$ that starts in a $K$-box that $*$-neighbours $B(v)$, such that $P$ and $Q$
lie in different open connected components of the percolation in~$\boxbig(v)$.
\end{enumerate} 
When $\cluster_u$ occurs, we write $\mc{C}_u$ for the cluster $\mc{C}$ containing $0$, for emphasis.
 \begin{lemma}\label{l.bdbad}
 Let $u \in (0,\infty)$ and $K \in \N$ satisfy $u \geq K >0$.
 Suppose $\cluster_u$ occurs.
 Let $B = B(v)$ be a $K$-box (with $v \in 2K \cdot \Z^d$) and let $D$ be a $K$-box $*$-neighbour of $B$. If $B \cap \cset_u = \emptyset$ and $D \cap \cset_u \not= \emptyset$, then $B$ is bad.  
 \end{lemma}
 {\bf Proof.} We will show that if the first condition that would qualify $B$ as bad is not met, then the second one is.  Suppose then that ${\rm boxlet}(v)$ is not isolated, and let $P$ denote an open path from this set to $\partial B$. Let $x$ be a vertex in both $\cset_u$ and $D$. 
  The diameter of $\cset_u$, being at least $u$, is at least $K$, so that $\max \big\{ d(x,v): v \in V(\mc{C}_u) \big\} \geq K/2$ (where $d$ is Euclidean distance). We may thus find 
 an open path~$Q$ in $\cset_u$ of diameter between $K/2$ and $K/2 + 1$ that emanates from $x$.  We see that
 the path $Q$ lies in $\boxbig(v)$, because $x$ lies in this big box at distance at least $2K$ from its boundary.  The path $P$ lies in $B$,  so that $P$ is vertex disjoint from $\cset_u$. Since $\mc{C}_u$ is an open cluster, there is no open path from $\cset_u$ to $P$; and specifically, none in $\boxbig(v)$.
 We have shown what we sought to show:  if ${\rm boxlet}(v)$ is not isolated, then $B$ is bad via the second condition specifying this notion.  This completes the proof of  Lemma~\ref{l.bdbad}. \qed

In the next result, we again rely on the indexing of $K$-boxes by $\Z^d$.
\begin{lemma}\label{l.percodominate}
Define a random process on $\Z^d$ by declaring that $v \in \Z^d$ is bad if the $K$-box $B(v)$ is bad. For any $\e >0$, there exists $K \in \N$ such that the process of bad indices---a random subset of $\Z^d$---is stochastically dominated by site percolation of parameter~$\e$. 
\end{lemma}
{\bf Proof.}
A given big box $D$ can be contained in the set $B'_{\bf i}(N)$ specified in \cite[Equation~$(2.7)$]{AntalPisztora}, where the parameter $N$ in this reference is set to be minimal to permit this containment.
It is readily verified that, if Centre$(D)$ is bad, then the complement of the event $R_{\bf i}^{(N)}$ in~\cite[Equation~$(2.9)$]{AntalPisztora} occurs. It is~\cite[Equation~(2.24)]{AntalPisztora} which establishes the rarity of this complementary event. Indeed, this bound shows that any given index in $\Z^d$ is bad with probability at most $\e$, where $\e > 0$ tends to zero in the limit of high $K$. The process of bad indices is $5$-dependent in the sense of~\cite{LSS}. Thus Lemma~\ref{l.percodominate} follows from \cite[Theorem~$0.0$(i)]{LSS}. \qed

We indicated in Section~\ref{s.slideintorenewal} how the slide resample will operate with a small gap, and next we present some related generalities of later use.
 
\begin{definition}\label{d.plustwominustwo}
Let $A$ denote a subgraph of $(\Z^d,\sim)$.
We denote by $\overline{A}[2]$  the set of $K$-boxes whose $*$-distance from $\overline{A}$ is most two.
And we write  $\overline{A}[-2]$ for the set of $K$-boxes in  $\overline{A}$ whose $*$-distance from the complement of $\overline{A}$ is at least three (so that at least {\em two} elements in any $*$-path in $K$-boxes between $\overline{A}[-2]$ and the complement of $\overline{A}$ lie outside of $\overline{A}[-2]$).
\end{definition}
\begin{lemma}\label{l.condbadperc}
Let $A$ be a subgraph of $(\Z^d,\sim)$. Let $\mathsf{E}$ denote a percolation event that is measurable with respect to the edges in $A$. 
Let $\badcoll$ denote the random set of bad $K$-boxes that do not lie in~$\overline{A}[2]$. Then the distribution of $\badcoll$ is equal under the laws $\PP$ and $\PP ( \cdot \vert \mathsf{E} )$.
\end{lemma}
{\bf Proof.} The value of the random process $\badcoll$ is determined by the open/closed status of the collection of edges lying in the union~$\mc{U}$ of the big boxes whose Centre is a $K$-box that does not lie in $\overline{A}[2]$. Although~$\mc{U}$ may contain edges that lie in the boundary of a $K$-box in $\overline{A}$, such edges do not belong to $A$, because, if they did, an adjacent $K$-box would also lie in $\overline{A}$ at $*$-distance two from Centre of a big box that lies in the union~$\mc{U}$. Thus the process $\badcoll$ is unperturbed by conditioning the law $\PP$ on the status of edges in $A$. \qed

\subsubsection{Bounding the cardinality of the cluster boundary}

Here we prove Proposition~\ref{p.nobigcluster}.

When $\cluster_u$ occurs, 
we write
 $\overline{\cset}_u$ in place of $\overline{\cset_u}$. 
\begin{lemma}\label{l.isoprep}
When $\cluster_u$ occurs,
every element of $\partial_{\rm ext} \overline{\cset}_u$ is a bad $K$-box. 
\end{lemma}
{\bf Proof.} Let $B \in \partial_{\rm ext} \overline{\cset}_u$. Let $D$ be a $K$-box $*$-neighbour of $B$ that lies in $\overline{\cset}_u$.
 Lemma~\ref{l.bdbad} then implies that $B$ is bad. \qed

\begin{lemma}\label{l.bigwell}
For $d \geq 2$, there exist positive constants $C$ and~$c$ such that, for $u \in (0,\infty)$ and $D \geq 1$,
$$
\PP \Big( \cluster_{[u,u+1]} , \big\vert \partial_{\rm ext} \overline{\cset} \big\vert \geq   D u \Big) \, \leq \, C \exp \big\{ - c  K^{-d} D u \big\} \, .
$$
\end{lemma}
{\bf Proof.} Suppose the occurrence of $\cluster_{[u,u+1]}$ and that $\big\vert \partial_{\rm ext} \overline{\cset} \big\vert = m$ for a given $m$ whose value is at least $D u$. Since $0 \in V(\cset)$,  we see that $\partial_{\rm ext} \overline{\cset}$ contains a $K$-box whose $\Z^d$-index has Euclidean norm at most $m/K$. 

By Lemmas~\ref{l.percodominate} and~\ref{l.isoprep}, we see that $\cluster_{[u,u+1]} \cap \big\{ \big\vert \partial_{\rm ext} \overline{\cset} \big\vert  = m \big\}$ entails that a site percolation  of parameter $\e$ in~$\Z^d$
contains a connected component  of size at least $m (2K+1)^{- d}$ to which a vertex belongs   whose distance from the origin is at most~$m$. Here, the parameter $\e > 0$ is a function of $K \in \N$ that may be chosen to be arbitrarily small by a high enough choice of~$K$.
The probability of the just mentioned eventuality is at most $C m^d \exp \big\{ - c  m K^{- d} \big\}$ for suitable positive constants $C$ and~$c$ in view of the exponential decay of the  tail of the size of the     
 connected component containing a given vertex in subcritical site percolation~\cite{AizenmanBarsky}. Summing this bound over $m \geq Du$, and using $D \geq 1$ to absorb the factor of $m^d$, we obtain Lemma~\ref{l.bigwell}. \qed

{\bf Proof of Proposition~\ref{p.nobigcluster}.} That there exists $c > 0$ such that $\PP(\cluster_{[u,u+1]}) > c^u$
is demonstrated by considering an explicit cluster, such as a concatenation of $\lceil u \cdot (e_d \cvl)^{-1} \rceil$ $e_d$-directed edges rising from~$0$.  Proposition~\ref{p.nobigcluster} then follows from Lemma~\ref{l.bigwell}
by making a suitably high choice for the lower bound $D_0$ on $D$ and by absorbing the term $K^{-d}$, which is constant for our given value of $K$, into the constant $c$ in the argument of the exponential. Note that the value of $D_0 \in (0,\infty)$ may be chosen independently of $\vecell \in S^{d-1}$ since $e_d \cvl \geq d^{-1/2}$ in view of the axes ordering convention~(\ref{e.aoc}). \qed

\subsubsection{Latitude and isoperimetry for $K$ boxes}\label{s.boxisoperimetry}

It is useful to know that when an object such as the downcluster is large in a slice, the enclosing surface is also large there.
Here we state and prove Proposition~\ref{p.dn.gen} to this effect, via Lemma~\ref{l.lift}.

\begin{definition}\label{d.boxlatitude}
Recall that a $K$-box takes the form $x + \llbracket -K,K \rrbracket^d$ with $x \in 2K \cdot \Z^d$. The latitude of $B$ equals $\bigl \lfloor \frac{x \cvl}{2K e_d \cvl} \bigr  \rfloor$.
 Latitude measures $K$-box location relative to the hyperplane $z \cvl = 0$ in integer units. 
\end{definition}
 For $k \in \Z$ and $\ell \in \N_+$, let $\overline{N}_{k,k+\ell}$ denote the set of $K$ boxes of latitude  in $\llbracket k,k+\ell \rrbracket$.

\begin{proposition}\label{p.dn.gen}
For some $k \in \Z$, suppose that $\overline{A}$ and $\overline{B}$ 
are finite subsets of $\nkell$
such that $\overline{B}$ separates $\overline{A}$
from infinity in $\nkell$. Then 
$$
 \vert \overline{A} \vert \, \leq \, (\ell+1) \vert \overline{B} \vert^{\tfrac{d-1}{d-2}} \, .
$$ 
\end{proposition}

Writing $\mc{B}$ for the set of $K$-boxes, we introduce the projection $P:\mc{B} \to \Z^{d-1}$
that sends $2Ku + \llbracket -K,K \rrbracket^d$ to $(u_1,\cdots,u_{d-1})$: namely, $P$ sends any $K$-box to the projection on the first $d-1$ coordinates of its natural index in $\Z^d$.
\begin{lemma}\label{l.lift}
Let $\overline{C} \subset \overline{N}_{k,k+\ell}$. Suppose that $x \in \Z^{d-1}$ is not separated from infinity by $P(C)$ in~$\Z^{d-1}$.
Then nor is any $B \in \overline{N}_{k,k+\ell}$ with $P(B) = x$
separated from infinity by~$\overline{C}$ in $\nkell$.
\end{lemma}
{\bf Proof.}
Let $y$ and $z$ be nearest neighbours in $\Z^{d-1}$, with $z = y + e_i$ for $i \in \intint{d-1}$.
Let $\mc{B}[y]$ denote the set of $K$-boxes in $\nkell$ that project to $y$ under $P$. 
Then $\vert \mc{B}[y] \vert$ equals $\ell +1$; label the lowest and highest elements of $\mc{B}[y]$ in the form 
 $B_y$ and $B'_y = B_y + 2K \ell e_d$.
Similarly for~$z$. We {\em claim} that some element of $\mc{B}[z]$ is obtained from some element of $\mc{B}[y]$ by translation by the vector $2K e_i$.
Indeed, if this were not the case, then either $B_z$ would be displaced from $B'_y + 2K e_i$ by a positive multiple of $2K e_d$; or $B'_z$ would be displaced from $B_y + 2K e_i$
by a negative multiple of $2K e_d$. In the former case, $B'_y + 2K e_i$ would have latitude less than $k$; in the latter, $B_y + 2K e_i$ would have latitude more than $k+\ell$. However, $0 \leq e_i \cvl \leq e_d \cvl$
by the axes ordering convention; and thus the difference in latitude between boxes $B$ and $\hat{B}$ satisfying $\hat{B} = B + 2K e_i$ is seen to be at most one in absolute value.  
Since $B'_y$ has latitude $k+\ell$ and $B_y$ has latitude $k$, we have found a contradiction and proved the claim.
 
 Whenever elements of $\nkell$ respectively project to nearest neighbours $y,z \in \Z^{d-1}$, we may invoke the claim to construct a  path of at most  $2(\ell +1)$ $K$-boxes that connects these elements, with the path valued in $\nkell$ among $K$-boxes that themselves project to one or other of $y$ and $z$. 
This local construction permits us to lift paths from $\Z^{d-1}$ to 
$\nkell$, and such a lifted path will prove Lemma~\ref{l.lift}.
Indeed, since $x \in \Z^{d-1}$ is not separated from infinity by $P(\overline{C})$ in $\Z^{d-1}$, there exists a path $Q$ in $\Z^{d-1}$ from $x$ to infinity that does not visit $P(\overline{C})$.
Arbitrarily select a sequence of elements of $\nkell$ that successively project to the consecutive elements of $Q$. Each consecutive pair of elements in this sequence may be interpolated by a path in $K$-boxes of length at most $2(\ell+1)$, by the local lift. When the resulting paths are concatenated, the path in $K$-boxes obtained reaches infinity without any visit to $\overline{C}$. This path validates the conclusion of Lemma~\ref{l.lift}. \qed

{\bf Proof of Proposition~\ref{p.dn.gen}.}
The first assertion of \cite[Lemma~$2.3$]{GLA} 
is a simple consequence of the Loomis-Whitney inequality~\cite{LoomisWhitney}: 
if $R \subseteq \Z^{d-1}$
denotes the collection of sites that a finite set $W \subseteq \Z^{d-1}$
separates from
infinity, then $\vert R \vert \leq \vert W \vert^{(d-1)/(d-2)}$.

Since $\overline{A}$ is separated from infinity in $\nkell$ by $\overline{B}$,
 Lemma~\ref{l.lift}
 implies that   
 $$
 \textrm{$P \big( \overline{A} \big)$ is separated from infinity in $\Z^{d-1}$ by $P ( \overline{B})$} \, .
 $$
 Taking $R = P( \overline{A})$ in the recalled bound, we  find that
 $\big\vert P ( \overline{A} ) \big\vert \leq \big\vert P ( \overline{B} ) \big\vert^{\tfrac{d-1}{d-2}}$.
 Using $\vert\overline{A} \vert \leq (\ell +1)    \big\vert P ( \overline{A} ) \big\vert$
 and $\vert P (\overline{B}) \vert \leq \vert \overline{B} \vert$, we find that
 $\vert \overline{A} \vert \, \leq \, (\ell +1) \big\vert \overline{B}  \big\vert^{\tfrac{d-1}{d-2}}$
 and thus obtain Proposition~\ref{p.dn.gen}. \qed

\section{Renewal levels via the slide resample}\label{s.renewalslide}

 Here we prove both parts of Theorem~\ref{t.renewal} by
 slide resampling  percolation conditional on the presence of the cluster.
 As we outlined in the preceding section,
 examining the slide resample when it operates at a single height is enough to prove Theorem~\ref{t.renewal}(1).
 The resample will be applied iteratively at many heights in order to obtain Theorem~\ref{t.renewal}(2). 
 
 There are four subsections. Critical to the utility of the slide for finding renewal levels is the {\em simple join} Proposition~\ref{p.success}, which shows that, when the slide resample operates at a given height on a cluster with typical geometry, it often results in a resampled cluster with a clean join between the downcluster and the upcluster. In Section~\ref{s.onetime}, we introduce the slide resample and demonstrate Theorem~\ref{t.renewal}(1) given  Proposition~\ref{p.success}, whose proof is the subject of Section~\ref{s.simplejoin}.
 As we turn to iterative use of the slide to reach the stronger Theorem~\ref{t.renewal}(2), we face the problem that the slide at a given height may fail to produce a simple join.
 In principle, failure may undo renewal levels far from the present slide height  that have been secured at earlier slide iterations. Section~\ref{s.damagecontrol}
 is devoted to showing that, when the slide fails, damage to the cluster  typically occurs only near the slide height.  
 With this support, we introduce the iterated slide and prove Theorem~\ref{t.renewal}(2) in Section~\ref{s.iterate}. 
 

The proof of Theorem~\ref{t.renewal} may be reprised to give Corollary~\ref{c.renewal}, the counterpart result in which the conditioned cluster is supposed to have renewal points at both extreme heights. The proof of Corollary~\ref{c.renewal} is given  at the end of Section~\ref{s.iterate}, alongside the derivation of 
Proposition~\ref{p.stingrenewal}, a string counterpart to Corollary~\ref{c.renewal}(1).

\subsection{The single slide: positive density of renewal}\label{s.onetime}

As we begin to prove Theorem~\ref{t.renewal}(1), we mention that we will develop the ideas of Section~\ref{s.slideintorenewal} closely enough, though since our context is now of a general rather than an axially aligned vector~$\vecell$, certain adjustments are needed. In the axial case, we could speak of proposed slides as horizontal, and indexed clearly by a copy of $\Z^{d-1}$. The vector $\vecell$ may however be non-lattice, and we have to accept that there no natural substitute for the lattice $\Z^{d-1}$ as an index of possible shifts. It is inevitable that the height of the cluster head will change, albeit perhaps only by a bounded quantity, under the actions of reasonable shifts. This consideration even finds expression in the form of Theorem~\ref{t.renewal}, where we condition on a cluster of height in a short interval $[u,u+3r]$, rather than of the specific value $u$; the latter choice is a luxury that is unavailable in the non-lattice case. We will remark at suitable moments on how the proof of Theorem~\ref{t.renewal}(1) varies from the template offered in Section~\ref{s.slideintorenewal}, including in regard to the need to accomodate general values of $\vecell$.

Let $u,h \in (0,\infty)$ satisfy $u > h$.
We introduce the upspace and the downspace, two half-spaces with boundaries near the slide height~$h$.
 The {\em downspace} is simply the lower-half space $L_h$.
 
 Recall from a specification made in Subsection~\ref{s.notation} and from the axes ordering convention~(\ref{e.aoc}) that $r = e_d \cvl$.
 
 The upspace is specified only under the event $\clusteruthree$, and then by means of a little more notation. 
 When this event  occurs, 
 there is a  cluster~$\mc{C}$, namely a finite open connected component with base $0$, whose  height lies in $[u,u+3r]$; for emphasis, we denote it by $\csetuthree$.
On the same event, we  write 
$$
\hmax = \max \{ z \cdot \vecell: z \in V(\csetuthree) \}
$$
for the maximum height among vertices in $\csetuthree$. In this way, 
  $\head(\csetuthree)$ is   
  the lexicographically minimal vertex in $\csetuthree$ of height~$\hmax$. 
  The upspace  is then defined to be the upper half-space $U_{\hadj}$, where the {\em adjusted height} $\hadj$ denotes $h + (\hmax -u)$.  When $\clusteruthree$ occurs, $\hmax \in [u,u+3r]$, so that the height $\hadj 
  \in [h,h+3r]$ is adjusted from $h$ by an addition of at most~$3r$.



Also on $\clusteruthree$, we define the upcluster and the downcluster to be two subgraphs of $\csetuthree$. 
First recall that, on this event, the cluster $\csetuthree$ satisfies $\base(\csetuthree) = 0$.  The {\em downcluster}, denoted $\downcluster$, is the connected component containing $0 \in \Z^d$ of the subgraph  induced by $\csetuthree$ in the downspace $L_h$. 
To wit, $v \in \Z^d$ lies in the vertex set of the downcluster if 
there exists an open path in $L_h$ that starts at the origin and that ends at $v$. The edge set of the downcluster is the  union of edges belonging to all such open paths. 
 The top of the downcluster, $\Top(\downwell)$, is the set of vertices~$v$ in $\downwell$ 
 that lie in $\vertbdry_h$ (as this set is defined in Section~\ref{s.basics}).
  
The {\em upcluster}, denoted $\upcluster$, is the open connected component containing $\head(\csetuthree)$ of the subgraph induced by $\csetuthree$ in the upspace. 
Thus, an element $v \in \Z^d$ lies in the vertex set of the upcluster if and only if an open path in  $U_{\hadj}$  runs from $\head(\csetuthree)$ to $v$.
The edge set of the upcluster is the union of the edges in such paths. Note that  $\head(\upwell) = \head(\csetuthree)$. The bed of the upcluster, $\Base(\upwell)$, is the set of vertices~$v$ in $\upwell$ 
that  lie in $\vertbdry_{\hadj}$. 

The definition of the downcluster is largely consistent with the heuristics of Section~\ref{s.slideintorenewal}, but the upcluster is specified a little differently: the hyperplane that delimits the upspace (and the upcluster) on the lower side is set at a height that is $u-h$ units below that of the upcluster head. This definition is made to face the challenge of setting up the slide resample when $\vecell$ may be non-lattice and some counterpart must be found of the  space $\Z^{d-1}$  of admissible shifts of the cluster (or the upcluster) head.
By specifying the upcluster relative to its head, we obtain a definition of the upcluster that will be invariant under a suitable class of shifts.

Next, some simple height facts.
\begin{lemma}\label{l.updowncluster}
Let $H_d$ and $H_u$ denote the heights of $\downcluster$ and $\upcluster$.
Almost surely under~$\PP$ given $\clusteruthree$, $h - r \leq H_d \leq h$ and $u-h \leq H_u \leq u-h+r$. 
\end{lemma}
{\bf Proof.} In view of Lemma~\ref{l.geometrybasics}(1), we find that $h - r \leq H_d \leq h$ 
since $V(\downcluster)$ consists of vertices in $L_h$ of non-negative height, one of which  lies in $\vertbdry_h$.
In regard to $H_u$, note that elements in $V(\upcluster)$
have height at most $\hmax$ and lie in the upspace $U_{\hadj}$, with at least one  lying in $\vertbdry_h$. 
Any element of the upspace has height at least $\hadj - r$, because, in the opposing case and by the axes ordering convention, the interior of each incident edge would lie in $B_{\hadj}$.
And an element that  lies in $\vertbdry_h$ has height at most $\hadj$, because, were it to be higher, the interior of each incident edge would intersect~$F_{\hadj}$.
Thus, we confirm that  $u-h \leq H_u \leq u-h+r$.  \qed

We now present the decomposition of the cluster into three pieces that will allow us shortly to specify the midcluster $\midcluster$.
\begin{definition}\label{d.residual}
Under $\clusteruthree$, the {\em residual cluster}, denoted $\residualwell$, is in essence the residue of the cluster~$\csetuthree$ after the removal of the upcluster and the downcluster. Formally, it is the graph whose vertices are endpoints of edges in $\csetuthree$ that are edges in neither the upcluster nor the downcluster; the edge-set of $\residualwell$ is the collection of such edges. 
\end{definition}
The residual cluster is comprised of possibly several connected components, as we now describe.
\begin{lemma}\label{l.residualwell}
\leavevmode
\begin{enumerate}
\item The graph $\residualwell$ contains at least one connected component $G$ such that $V(G)$ 
intersects  $\Base(\upwell)$ and  $\Top(\downwell)$.
\item The vertex set of any connected component of $\residualwell$ intersects the union of the vertex sets of $\downwell$ and $\upwell$. Any vertex in this intersection lies in the union of  $\Top(\downwell)$ and~$\Base(\upwell)$.
\end{enumerate}
\end{lemma}
{\bf Proof: (1).} Since $\base(\csetuthree)$ and $\head(\csetuthree)$ lie in the connected subgraph $\csetuthree$, we may find a path in~$\csetuthree$ 
that runs from 
 $\base(\csetuthree)$ to $\head(\csetuthree)$. This path begins in $\downwell$ and ends in $\upwell$. 
Consider the subpath between the final vertex in $\downwell$ and the first in $\upwell$. Since these two vertices 
 lie in $\vertbdry_h$,
we may set $G$ to be the connected component of $\residualwell$ that contains this subpath.

{\bf (2).} Let $H$ be a connected component of $\residualwell$. Since the cluster $\csetuthree$ is connected, a path in~$\csetuthree$ runs from any given element of $H$ to $\base(\csetuthree)$. When this path visits $H$ for the final time, it does so at a vertex which also lies in the vertex set of either $\upwell$ or $\downwell$---in the former case, at an element of $\Base(\upwell)$; in the latter, at an element of $\Top(\downwell)$. \qed

We are ready to select and label as $\midwell$ a connected component of $\residualwell$ that connects the upcluster and the downcluster.
The residual cluster may also contain certain further fragments, which touch either $\Base(\upwell)$ or $\Top(\downwell)$ and are not central players in our analysis, that we label $\midfragment$.
\begin{definition}
\leavevmode
\begin{enumerate}
\item
Among the components $G$ specified in Lemma~\ref{l.residualwell}(1), we select one and call it the {\em midcluster} $\midwell$. For definiteness, $V(\midwell)$ contains the lexicographically minimal element among vertices in any of the concerned components. 
\item The union of connected components of $\residualwell$ excluding $\midwell$ will be denoted by $\midfragment$.
\end{enumerate}
\end{definition}


 For $A \in \N$, we identify three slabs of width~$2A$ in the neighbourhood of height $h$: a north zone, a south zone, and a tropical zone sandwiched between them. (Here, latitude may be identified with height: it acts as an alternative verbal marker for the value of the coordinate in the direction~$\vecell$.)
  
  \begin{definition}\label{d.tropics}
   The subgraphs $\northzone$, $\midzone$ and $\southzone$  are each slabs, respectively set equal to $\slab_{h+A,h+3A}$, $\slab_{h-A,h+A}$ and $\slab_{h-3A,h-A}$. 
 \end{definition}

\begin{definition}
The {\em upwell configuration space} $\ucs$ consists of  all values that can be adopted by the subgraph $\upwell$ in a percolation configuration that realizes the event~$\clusteruthree$.
\end{definition} 

 We now introduce an equivalence relation on this configuration space by identifying pairs of elements if one element of the pair may be obtained from the other by translation.
 \begin{definition}\label{d.equivalence}
 Two elements $U$ and $U'$ of $\ucs$ are {\em shift equivalent} if there exists $v \in \Z^d$ such that $U'$ is the graph obtained from $U$ by shifting by $v$ (so that the vertex and edge sets of $U'$ are obtained from the counterparts for $U$ by shifting by this vector). 
 
 We will write $[\upwell]$ for the shift equivalence class of $\upwell$. Thus,  $[\upwell]$ is a random variable that is defined whenever $\clusteruthree$ occurs.
 \end{definition}

In Section~\ref{s.slideintorenewal}, we indicated how the slide resample would be liable to act in a desired fashion if the upcluster and downcluster cardinalities are suitably limited at heights close to the resample height~$h$. We now set down conditions that express these favourable features.
\begin{definition}\label{d.fvf}
Let $R$, $\macthin$ and $\constfine$ be three positive parameters. Here we introduce several conditions of smallness for the cardinality of $\csetuthree$.
The conditions receive labels ${\bf T}$, ${\bf S}$ and ${\bf G}$ according to whether they concern the size of $\csetuthree$ in the {\em tropics} (with heights in $[h-A,h+A]$);
the {\em subtropics} (the southzone or the northzone); or {\em globally} (when the total size of $\csetuthree$ is considered).
Here are the conditions.
\begin{itemize}
\item [${\bf T}_1(R)$:] for $i \in \intint{A}$, the number of vertices in $\csetuthree$ with height in $[h+i-1,h+i]$ is at most~$i^3R$.
\item [${\bf T}_2(R)$:] for $i \in \intint{A}$, the number of vertices in $\csetuthree$ with height in $[h - i,h-i+1]$ is at most~$i^3R$.
\item  [${\bf S}_1(\macthin)$:] the number of vertices in $\csetuthree$
with height in $(i-1,i]$ is at most $\macthin$
for at least one-third of indices $i$ in 
$\llbracket h+1+A,h+3A\rrbracket$. 
\item   [${\bf S}_2(\macthin)$:]
 the number of vertices in $\csetuthree$
with height in $(i-1,i]$ is at most $\macthin$
for at least one-third of indices $i$ in 
$\llbracket h+1-3A,h-A\rrbracket$. 
\item  [${\bf G}(\constfine)$:] 
 $\vert V(\csetuthree) \vert \leq \constfine u^d$.
\end{itemize}
We further consider three adjusted conditions, marked with a prime.
\begin{itemize}
\item [${\bf T}'_1(R)$:] for $i \in \intint{A}$,  the number of vertices 
in $\upcluster$
 whose height lies in $[\hadj+i-1,\hadj+i]$ is at most $i^3R$.
\item [${\bf T}'_2(R)$:] for $i \in \intint{A}$, the number of vertices in $\downcluster$ with height in $[h - i,h-i+1]$ is at most~$i^3R$.
\item  [${\bf S}'_1(\macthin)$:]  the number of vertices in $\csetuthree$
with height in $(i-1,i]$ is at most $\macthin$
for at least one-quarter of indices $i$ in 
$\Z \cap[ \hadj +1+A, \hadj +3A ]$.
\end{itemize}
Notably, the north tropical and subtropical conditions  ${\bf T}'_1(R)$  and
 ${\bf S}'_1(\macthin)$ are expressed in terms of the adjusted equatorial height $\hadj$ rather than the ordinary equatorial height $h$.
 Equivalent elements of $\ucs$ are equal when viewed in these coordinates: this makes the coordinates useful, because it means that the primed conditions will be automatically maintained in the slide resample.
 In other regards, the changes of the primed events relative to the unprimed counterparts are weakenings of conditions: the replacement of $\csetuthree$ by its subgraphs $\upcluster$ and $\downcluster$; and of one-third by one-quarter.

(The global condition $\vert V(\csetuthree) \vert \leq \constfine u^d$ may seem absurdly weak, with $\vert V(\csetuthree) \vert \leq \constfine u$ being more natural. The weak form is adequate for our purpose and avoids some technical encumbrances.)

In terms of positive parameters $K$, $R$ and $\constfine$, we further define 
$$
\textrm{two subevents $\veryfine_h
= \veryfine_h(K,R,\constfine)$
and $\fine_h
= \fine_h(K,R,\constfine)$ 
of $\clusteruthree$} \, , 
$$
using primed (as well as some unprimed) conditions in the latter case.

\begin{figure}[htbp]
\centering
\includegraphics[width=0.6\textwidth]{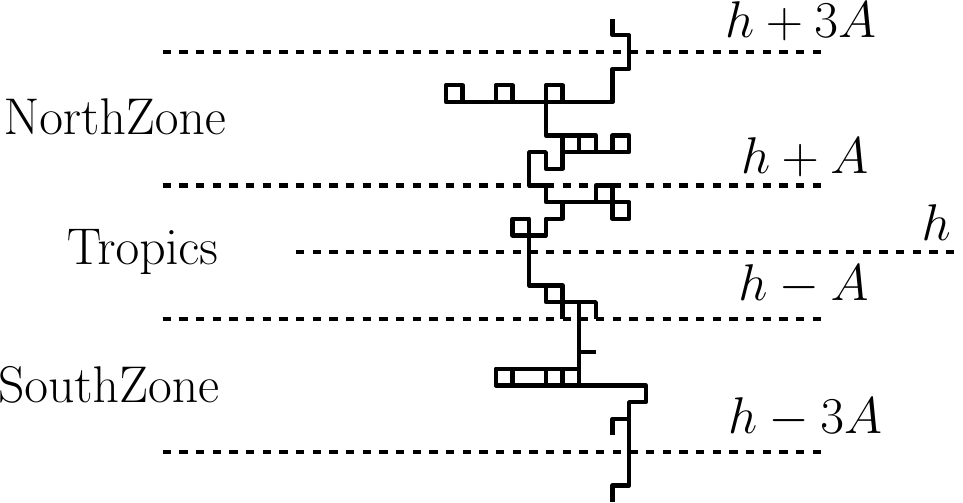}
\caption{A schematic instance of the event $\veryfine_h$ viewed in the region around the tropics. The cardinality of the intersection of the depicted cluster $\csetuthree$ with horizontal lines in $\tropics$ is bounded above quadratically in the distance from the equatorial line at height $h$, while a positive density of such cardinalities must be bounded above in the neighbouring subtropical regions.}\label{f.tropicscluster}
\end{figure}

Indeed (as Figure~\ref{f.tropicscluster} illustrates), $\veryfine_h$ occurs when the following conditions are met:
$$
  {\bf T}_1(2^{-1}\tenmac^{-3}R) \, , \,  {\bf T}_2(2^{-1}\tenmac^{-3}R) \, , \, {\bf S}_1(\macb) \, , \, {\bf S}_2(\macb) \, \, \textrm{and} \, \, {\bf G}(\constfine) \, .
  $$
And $\fine_h$ occurs when it is these conditions that are met:
$$
  {\bf T}'_1(R) \, , \, {\bf T}'_2(R) \, , \, {\bf S}'_1(\macb)  \, , \, {\bf S}_2(\macb)  \, \,  \textrm{and}  \, \, {\bf G}(\constfine) \, .
$$
\end{definition}

\begin{lemma}\label{l.veryfineimpliesfine}
When $A \geq 36$, we have that $\veryfine_h \subset \fine_h$. 
\end{lemma}
{\bf Proof.} In regard to the tropical conditions, note that the number of vertices in $\upcluster$ with height in $[\hadj +i-1,\hadj+i]$
is at most the number of vertices in $\csetuthree$ with height in $[h + \kappa-2,h+\kappa]$, where $\kappa = \lfloor \hadj - h \rfloor \in \llbracket 0, 2 \rrbracket$.
If $\veryfine_h$ occurs, then this upper bound is at most $\big( 2^3 + \tenmac^3 \big) 2^{-1}\tenmac^{-3}R \leq R$. So conditions   $\ftOne$ and $\ftTwo$  of $\fine_h$ have been verified.  In regard to the northzone condition, note that the symmetric difference between the index sets $\llbracket h+A,h+3A\rrbracket$ and $\Z \cap[ \hadj +A, \hadj +3A ]$ is at most six. So if a set of integers occupies a proportion of at least one-third of the first of these intervals, it will occupy at least one-quarter of the second, provided that $A \geq 36$.
Thus $\vfsOne$ implies $\fsOne$. Since other conditions are either shared between the events $\fine_h$ and $\veryfine_h$, or are weaker in the former case, we obtain Lemma~\ref{l.veryfineimpliesfine}.
\qed

We will obtain Theorem~\ref{t.renewal}(1) by showing that a typical cluster has $\Theta(u)$ levels $h$ for which $\veryfine_h$ occurs; and that, no matter the form of $([\upwell],\downwell)$, 
the thus typical $\fine_h$ 
often entails a simple join event, in which $\upwell$ sits directly on top of $\downwell$, ensuring renewal at height~$h$. 
The next two propositions, whose proofs we defer, show $\Theta(u)$ cardinality of very fine levels and that fine levels often present simple joins.  

\begin{proposition}\label{p.veryfinelevels}
Let $N$ denote the number of  $h \in  \N \cap [0,u]$ for which $\veryfine_h$ occurs. With the relation $\constfine = 2^{\tfrac{d-2}{d-1}} R K^{-\tfrac{2d-3}{d-1}}$ set between the parameters in the definition of the event $\veryfine_h$,
there exist positive constants $C$ and $c$ such that
$$
 \PP \Big( N < (1-a)u \, \Big\vert \, \clusteruthree  \Big) \, \leq \,   C \exp \big\{ - c R^{\tfrac{d-2}{d-1}} K^{-1} a \big\} \, .
$$
\end{proposition}

Recall that $\base(\upcluster)$ and $\head(\downcluster)$
are certain elements in the subgraphs $\Base(\upcluster)$ and $\Top(\downcluster)$. For the purpose of the ensuing definition, take $v = \head(\downcluster) + e_d$.

The event $\success$ occurs when $\clusteruthree$ does alongside the conditions that
\begin{enumerate}
\item
 $\midwell$ has merely three vertices and two edges, the pair of edges 
 having the form $[v,v+e_d]$ and $[v-e_d,v]$; and
  \item the graph $\midfragment$ is empty.
 \end{enumerate}
On $\success$, the vertex $v = \head(\downcluster) + e_d$ will be called the {\em midvertex}.


We write $\PP^{\fine_h} = \PP \big( \cdot \big\vert \, \fine_h \, \big)$. And we further 
\begin{equation}\label{e.condnot}
\textrm{denote by $\PP_{[\upwell],\downwell}^{\fine_h}$  the law  $\PP^{\fine_h}$ conditionally on  $([\upwell],\downwell)$} \, . 
\end{equation}
We have yet to define the slide resample,
but it is perhaps helpful to use it to informally locate the above notation. Indeed, the slide resample at height~$h$
concerns the just denoted law, in which conditioning occurs not only on the event $\fine_h$ but also on the random variable $([\upwell],\downwell)$. 
Note a variation here on the theme of Section~\ref{s.slideintorenewal}: the slide is holding invariant the law of $\PP$ given $\fine_h \subset \cluster_{[u,u+3r]}$, alongside the data$([\upwell],\downwell)$, so that in addition to the occurrence of  $\cluster_{[u,u+3r]}$, the properties of $\fine_h$ (which are conditionally fairly typical) must be respected when the slide resample operates.

\begin{figure}[htbp]
\centering
\includegraphics[width=0.6\textwidth]{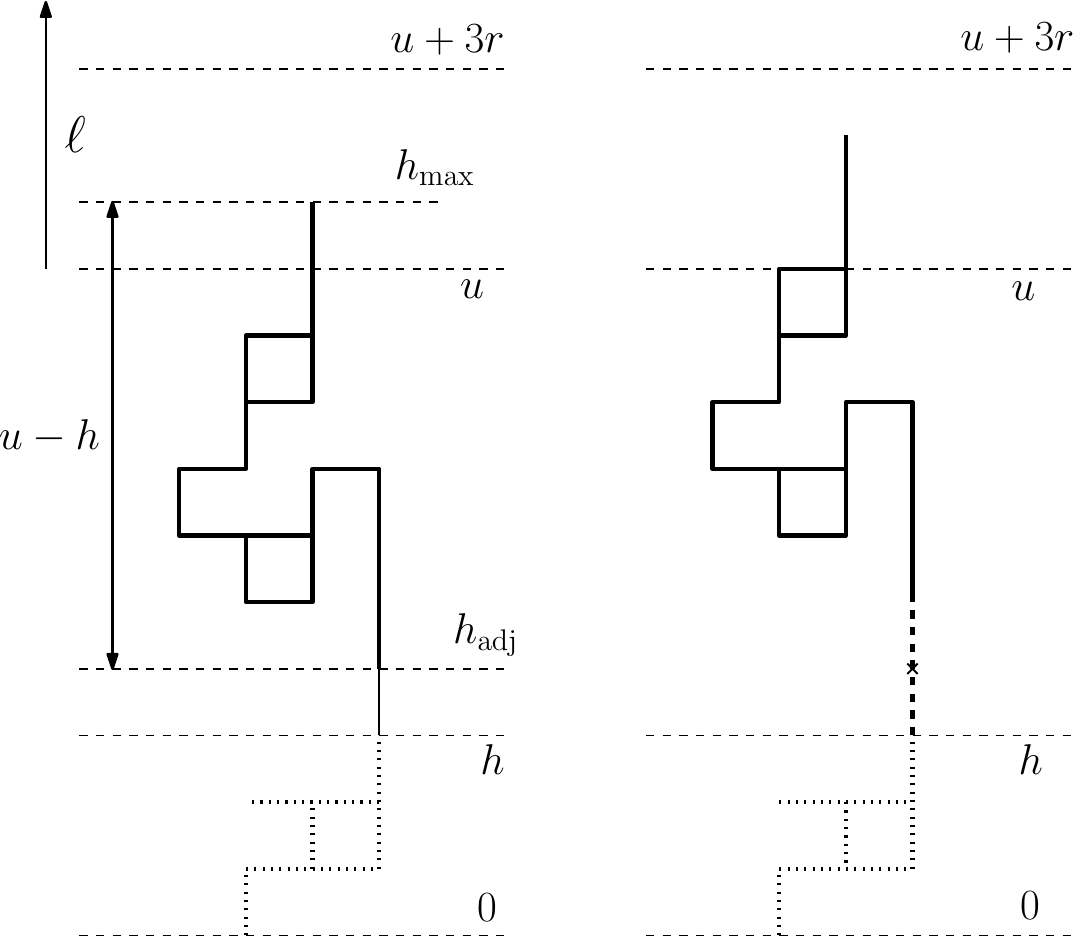}
\caption{The slide resample illustrated, with $\success$ occurring as a result. {\em Left}: 
$\upwell$ is bold, with height $u - h = h_{\rm max} - h_{\rm adj}$; $\downwell$ is bold dotted. {\em Right}: $\upwell$ moves up by one unit, to accomodate the two bold dashed edges that meet at the midvertex, marked with a cross.}\label{f.resample}
\end{figure}

The next result, which Figure~\ref{f.resample} illustrates, gives a lower bound on the probability that the slide resample realizes a simple join at height~$h$. 
\begin{proposition}\label{p.success}
For any positive $R$ and $K$, there exist $c = c(R,K) > 0$ and $C = C(K,d)$ such that the condition that $h$ and $u-h$ exceed $C \log u$ implies that   
$$
\PP^{\fine_h} \Big( \, \success \, \Big\vert \, [\upwell] , \downwell \, \Big) \, \geq \, c \, .
$$
\end{proposition}

And now an assertion, that a simple join is a renewal point, which will permit us to obtain  Theorem~\ref{t.renewal}(1) from the two preceding propositions.

 \begin{lemma}\label{l.success}
  When $\success$ occurs, the midvertex 
  is a renewal vertex of the cluster~$\csetuthree$.
  \end{lemma}
  The next claim will be used in proving this result. It is expressed using notation that develops the notion of vertex boundary $\vertbdry_h$ for $h \in \R$
  that is specified in Section~\ref{s.basics}. Indeed, we set   
   the {\em edge boundary} $\edgebdry_h$ equal to the set of elements of $E(\Z^d)$
that are incident to at least one element of $\vertbdry_h$. 
\begin{lemma}\label{l.twointernal}
Let $v \in \Z^d$ and $e \in E(\Z^d)$.
\begin{enumerate}
\item If $e \in E(U_{v + e_d})$ then $e \not\in \edgebdry_{v \cvl}$.
\item If $e \in E(L_{v - e_d})$ then $e \not\in \edgebdry_{v \cvl}$.
\end{enumerate}
\end{lemma}
 {\bf Proof.}
Let $w$ and $x$ denote the endpoints of $e$. 
By the axes ordering convention, $e_d \cdot \vecell \geq e_i \cdot \vecell \geq 0$ for $i \in \intint{d}$.
We also suppose that $w \cdot \vecell \geq x \cdot \vecell$. 
 By these assumptions, $x$ equals $w - e_i$ for some $i \in \intint{d}$.  
 
{\bf  (1).}   
   We need to show that neither $w$ nor $x$ are elements of $V(L_v)$.
Note that, since $e \in E(U_{v+e_d})$, and $w \cdot \vecell \geq x \cdot \vecell$, we have that $w \cvl \geq (v+e_d) \cvl$. Note then that 
$x \cvl \geq w \cvl - e_d \cvl \geq (v + e_d) \cvl - e_d \cvl \geq v \cvl$.

We know that $w \cvl$ either equals, or exceeds, $(v + e_d) \cvl$.
Suppose that equality holds. Then $x \cvl$ also equals $(v + e_d) \cvl$.
Indeed, in the opposite case, we would have $x \cvl < (v + e_d) \cvl$. But then the edge interior $(x,w)$
would lie in the backward space $B_{v + e_d}$ (as specified in Section~\ref{s.basics}), which would imply, contrary to our assumption, that $e \in E(L_{v+e_d})$.
We see then that, in the case that  $w \cvl = (v + e_d) \cvl$, the heights of both $w$ and $x$ exceed $v \cvl$.
This implies that neither $w$ nor $x$ lies in $V(L_v)$, so that the conclusion of Lemma~\ref{l.twointernal}(1)
holds in this case.

Suppose instead that   $w \cvl > (v + e_d) \cvl$.
From $x = w - e_i$ and $e_i \cvl \leq e_d \cvl$, we find that $x \cvl > v \cvl$. In this case, then, $w$ and $x$ both have height exceeding $v$'s. Thus, $w,x \not\in V(L_v)$,
so that Lemma~\ref{l.twointernal}(1) holds in this case also.

{\bf (2).} Here, we need to show that $w,x \not\in V(U_v)$.

Since $e \in E(L_{v - e_d})$ and the height of $w$ is at least that of $x$, we have that $w \cvl \leq (v- e_d) \cvl$.

Note that, if $z \in V(U_v)$, then $(z,z+e_j) \cap F_v \not= \emptyset$ for some $j \in\intint{d}$.
Thus, $z \cvl > v \cvl - e_d \cvl$.

For our purpose, it suffices then to show that $w \cvl$ and $x \cvl$ are at most $(v- e_d) \cvl$. Above we showed this for $w$, and by assumption $x \cvl \leq w \cvl$.
Thus the proof of Lemma~\ref{l.twointernal}(2) is complete. \qed

  {\bf Proof of Lemma~\ref{l.success}.} 
  We see that $\csetuthree = \downcluster \circ \midcluster \circ \upcluster$
  fragments at $v$ into two subgraphs whose edge sets are the union of $[v,v+e_d]$ and $E(\upcluster)$; and the union of $[v-e_d,v]$
  and $E(\downcluster)$.  
  To obtain Lemma~\ref{l.success}, it is  enough to argue that $e \in E(\upcluster) \cup  E(\downcluster)$
  implies that $e \not\in \edgebdry_{v \cvl}$.
  Given that $E(\upcluster) \subset E(U_{v+e_d})$
  and  $E(\downcluster) \subset E(L_{v-e_d})$, this statement follows from Lemma~\ref{l.twointernal}. \qed

{\bf Proof of Theorem~\ref{t.renewal}(1).} For $h \in [0,u]$, let $v = \head(\downcluster) + e_d$. Denote by  $\mathscr{R}_u(h) \subseteq \clusteruthree$ the event that $v$ is a renewal point in $\csetuthree$.
We claim that there exist positive $C$ and $c$ such that, when $h$ and $u -h$ exceed $C \log u$, 
\begin{eqnarray*}
& & \PP \big( \mathscr{R}_u(h) \big)  \, \geq \,
\PP \big( \clusteruthree , \success_h \big) \, \geq \,
\PP \big( \fine_h , \success_h \big) \\
& = & 
 \E \, \Big[ \, \PP_{[\upwell_h],\downwell_h}^{\fine_h} \big(  \success \big) \cdot {\bf 1}_{\fine_h} \, \Big] \, \geq \, 
c \, \PP \big( \fine_h \big) \, ,
\end{eqnarray*}
where we have used the notation~(\ref{e.condnot}).
Indeed, Lemma~\ref{l.success}, $\fine_h \subseteq \clusteruthree$ and Proposition~\ref{p.success} deliver the respective inequalities.  

For each $h \in [0,u]$, the associated value of  $\head(\downcluster)$ is an element of $V(L_h)$
that lies in $\vertbdry_h$. By Lemma~\ref{l.geometrybasics}(1), we see that, if two such values of $h$ differ by at least $r = e_d \cvl$, the associated vertices~$v$ cannot be equal. 
Thus, if we ascertain the occurrence of  a certain number of events $\mathscr{R}_u(h)$ for $h \in r \N \cap [0,u]$,
then we can be sure that there are at least this many elements of~$R(\csetuthree)$.

We may thus divide the bound derived above by $\PP (\clusteruthree)$  and sum the result over values $h \in r\N \cap [0,u]$ to find that
$\E \big[ \vert R(\csetuthree) \vert \big\vert \clusteruthree \big]$ is at least $c \, \E [\Nmac \vert \clusteruthree ]$, where 
 $\Nmac$ denotes the number of values of $h \in r\N \cap [0,u]$ for which $\fine_h$ occurs. Lemma~\ref{l.veryfineimpliesfine} implies that $\Nmac$ is at least the quantity $N$ specified in
Proposition~\ref{p.veryfinelevels}; thus, the latter result implies that $\E [\Nmac \vert \clusteruthree ] \geq c_1 u/r$ for some $c_1 > 0$. 
Whence, 
\begin{eqnarray*}
 & & ur^{-1} \, \PP \big(  \vert R(\csetuthree) \vert  \geq c_2 ur^{-1} \big\vert  \clusteruthree \big)  + c_2 ur^{-1} \, \PP \big(   \vert R(\csetuthree) \vert < c_2 ur^{-1} \big\vert  \clusteruthree \big) \\
 & \geq & \E \big[   \vert R(\csetuthree) \vert  \big\vert \clusteruthree \big] \, \geq \, c c_1 ur^{-1}
\end{eqnarray*}
and $\PP \big(    \vert R(\csetuthree) \vert  \geq c_2 ur^{-1} \big\vert  \clusteruthree \big)  \geq c c_1 - c_2$,
for any $c_2 > 0$.  Taking $c_2 = c c_1/2$ and relabelling $c > 0$, we obtain  Theorem~\ref{t.renewal}(1). \qed

We close this section by deriving Proposition~\ref{p.veryfinelevels}. The more substantial Proposition~\ref{p.success} will be proved in the next. 

Proposition~\ref{p.veryfinelevels} is reduced to Lemma~\ref{l.nhigh} and the latter is then proved by means of Proposition~\ref{p.dn.gen}.
\begin{lemma}\label{l.nhigh}
There  exists a positive choice of  the constant $c$  determined solely by $d \geq 3$ such that, when $\clusteruthree$ occurs, the condition
\begin{equation}\label{e.nhigh}
\big\vert \partial_{\rm ext} \overline{\cset}_u \big\vert \, \leq \,    
\min \Big\{ 
c R^{\tfrac{d-2}{d-1}} K^{-1} au , 4^{- \tfrac{d-2}{d-1}}(LK)^{\tfrac{d-2}{d-1}}u^{d-2} \Big\}
\end{equation}
ensures that $N \geq (1-a)u$.
\end{lemma}
{\bf Proof of Proposition~\ref{p.veryfinelevels}.} Set $D = c R^{\tfrac{d-2}{d-1}} K^{-1} a$.
The hypothesised condition  $\constfine = 2^{\tfrac{d-2}{d-1}} R K^{-\tfrac{2d-3}{d-1}}$  assures alongside $u \geq 1$ that 
 the minimum in~(\ref{e.nhigh})  is attained by the first term, whatever the value of $u \in [1,\infty)$. 
 We then obtain the sought result from Lemma~\ref{l.nhigh} and Proposition~\ref{p.nobigcluster} by suitably taking $C \geq D_0$ high and $c > 0$ low. \qed

{\bf Proof of Lemma~\ref{l.nhigh}.} 
First we argue that the global condition~$\fg$ in $\fine_h$, namely that $\vert V(\csetuthree) \vert \leq \constfine u^d$, is implied by  $\big\vert \partial_{\rm ext} \overline{\cset}_u \big\vert \leq 
 4^{- \tfrac{d-2}{d-1}}(LK)^{\tfrac{d-2}{d-1}}u^{d-2}$. Suppose that $\fg$ fails. Since $\csetuthree$ occupies heights in $[ 0,u ]$, 
 latitudes of elements of $\overline{\cset}_u$ are integers in $[0,u/(2Kr)]$, so they number at most $u/(2Kr) +1$, which is at most $u/(Kr)$
 since $u \geq 2Kr$. 
 There are thus more than $\constfine u^{d-1} Kr$ elements of $\overline{\cset}_u$ living at the same latitude. 
Calling this latitude $k \in \Z$, note that $\partial_{\rm ext} \overline{C}_u \cap \nkk$ separates 
$\overline{\cset}_u \cap \nkk$ from infinity in $\nkk$, so that Proposition~\ref{p.dn.gen} implies that
$$
 \big\vert \overline{\cset}_u \cap \nkk \big\vert \leq \ 2 \big\vert \nkk \cap \partial_{\rm ext} \overline{\cset}_u \big\vert^{\tfrac{d-1}{d-2}} \, ,
$$
where the notation $\nkk$ is introduced in Subsection~\ref{s.boxisoperimetry}. 
Since the left-hand side is at least $\constfine u^{d-1} Kr$,   $\partial_{\rm ext} \overline{\cset}_u$  contains more than $4^{- \tfrac{d-2}{d-1}}(LK)^{\tfrac{d-2}{d-1}}u^{d-2}$ elements. This is as we sought to argue. Hence, condition $\fg$ is satisfied, whatever the value of~$h$.

\begin{lemma}\label{l.latint}
For $q \in \R$,
let $\Gamma_q$ denote the set of latitudes of $K$ boxes that contain a vertex with $v \cvl \in [q,q+1]$. 
Then $\Gamma_q$ is an integer interval $\llbracket a,b \rrbracket$ with $a = \lfloor (q - K d^{1/2})(2K r)^{-1} \rfloor$
and $a + d^{1/2}r^{-1} \leq b \leq a + d^{1/2}r^{-1} + 1$.
\end{lemma}
{\bf Proof.} Suppose that the $K$-box $B = 2Ku + \llbracket -K,K \rrbracket^d$ contains $v$.
Then $\vert v \cvl - 2K u \cvl \vert \leq K d^{1/2}$, so that $v \cvl \in [q,q+1]$
implies that $q - K d^{1/2} \leq 2K u \cvl \leq q + K d^{1/2} +1$. 
The quantities $\lfloor (q - K d^{1/2})(2K r)^{-1} \rfloor$
and $\lfloor (q + K d^{1/2} +1)(2K r)^{-1} \rfloor$ are thus lower and upper bounds on the latitude of $B$. Lemma~\ref{l.latint} will follow once we conform that $\Gamma_q$ is an integer interval. 
If $\Gamma_q^+$ and $\Gamma_q^-$ denote the counterparts to $\Gamma_q$ in whose definition $[q,q+1]$ is respectively replaced by $[q,\infty)$ and $(-\infty,q+1]$
then $\Gamma_q$ may be written $\Gamma^+_q \cap \Gamma^-_q$ as an intersection of two semi-infinite integer intervals (that extend to infinity in opposing directions); whence we learn that $\Gamma_q$ is an integer interval indeed.  \qed

Let $F_i$ denote the set of $h \in \N \cap [0,u]$ such that there are at least $i^3 R$
elements $x \in V(\csetuthree)$ such that $x \cvl \in [h+i-1,h+i]$.
Let $G_i \subseteq F_i$
be chosen so that the distance between any pair of elements in $G_i$ exceeds $K d^{1/2} + 2$, with 
$\vert G_i \vert \geq \vert F_i \vert (K d^{1/2} + 3)^{-1}$.

Lemma~\ref{l.latint}
permits us to associate to each $j \in G_i$
a value $\ell_j \in \Gamma_{j+i-1}$ such that at least $i^3 R (d^{1/2}r^{-1} +2)^{-1}$
elements of $\overlinecsetuthree$
have latitude $\ell_j$. Elements of $\{ \Gamma_{j+i-1}: j \in G_i \}$
are at pairwise distance at least two, so that the differences between consecutive $\ell_j$ as indexed by $j \in G_i$
are also at least two. Thus, the sets $\overlinecsetuthree \cap \overline{N}_{\ell_j,\ell_j + 1}$  indexed by $j \in G_i$
are pairwise disjoint. Since $\overlinecsetuthree \cap \overline{N}_{\ell_j,\ell_j + 1}$ is separated from infinity in $\overline{N}_{\ell_j,\ell_j + 1}$ by $\partial_{\rm ext }\overlinecsetuthree$,
Proposition~\ref{p.dn.gen} implies that
$$
 \big\vert \overlinecsetuthree \cap \overline{N}_{\ell_j,\ell_j + 1} \big\vert \leq 2 \big\vert  \partial_{\rm ext }(\overlinecsetuthree) \cap \overline{N}_{\ell_j,\ell_j + 1} \big\vert^{\tfrac{d-1}{d-2}} \, .
$$
But
$$
 \big\vert \overlinecsetuthree \cap \overline{N}_{\ell_j,\ell_j + 1} \big\vert \geq \frac{i^3 R}{d^{1/2}r^{-1} + 2}
 $$
 so that 
 $$
 \big\vert  \partial_{\rm ext }(\overlinecsetuthree) \cap \overline{N}_{\ell_j,\ell_j + 1} \big\vert \geq 2^{\tfrac{2-d}{d-1}}
 \Big( \frac{i^3 R}{d^{1/2}r^{-1} + 2} \Big)^{\tfrac{d-2}{d-1}}
$$
for $j \in G_i$.
The disjointness of $\overline{N}_{\ell_j,\ell_j + 1}$
for $j \in G_i$ implies that
$$
\big\vert  \partial_{\rm ext }(\overlinecsetuthree) \big\vert \geq 2^{\tfrac{2-d}{d-1}}
 \Big( \frac{i^3 R}{d^{1/2}r^{-1} + 2} \Big)^{\tfrac{d-2}{d-1}}
 \vert G_i \vert \, .
$$
By assumption $\big\vert  \partial_{\rm ext }(\overlinecsetuthree) \big\vert \leq c R^{\tfrac{d-2}{d-1}} K^{-1}  au$. Recall that $\vert G_i \vert \geq \vert F_i \vert (K d^{1/2} + 3)^{-1}$.
Thus, 
$$
 \vert F_i \vert \leq c R^{\tfrac{d-2}{d-1}} K^{-1}  2^{\tfrac{d-2}{d-1}}  \Big( \frac{i^3 R}{d^{1/2}r^{-1} + 2} \Big)^{-\tfrac{d-2}{d-1}}
 (K d^{1/2} + 3)u \, ,
$$
where the right-hand term here is bounded above by $c \alpha i^{-3(d-2)/(d-1)} au$, with $\alpha$ a constant depending on $d$ but not on $R$ or $K$.
Thus $\sum_{i=1}^A \vert F_i \vert \leq c  \alpha \zeta(3/2)au$ for any $d \geq 3$.  
For a suitably small $d$-dependent choice of the positive constant $c$, we see that  $\sum_{i=1}^A \vert F_i \vert \leq au/4$. 
A quantity $F'_i$ may be defined as $F_i$ was, with the height interval $[h+i-1,h+i]$
replaced by $[h-i,h-i+1]$. An essentially identical argument to that just given yields the bound $\sum_{i=1}^A \vert F'_i \vert \leq au/4$.
The expression $\sum_{i=1}^A \big( \vert F_i \vert +\vert F'_i \vert \big)$ is an upper bound on the cardinality of the set of $h \in \N \cap [0,u]$
that violate at least one of the tropical conditions  ${\bf T}_j(R)$, $(j,i) \in \intint{2} \times \intint{A}$, that specify $\veryfine_h$.
Thus we see that the number of tropical condition violating values of $h$ is at most $au/2$.


Let $\chi_i$ equal the cardinality of the set of $x \in V(\csetuthree)$ for which $x \cvl \in (i-1,i]$ and let
$$
M = \Big\{ h \in \N \cap [0,u]: \chi_i > R \, \, \textrm{for at least $4A/3$ values $i \in \llbracket h+1+A,h+3A \rrbracket$} \Big\} \, .
$$
Note that $4A/3$ equals two-thirds of the cardinality of  $\llbracket h+1+A,h+3A \rrbracket$.

 Let $M$ denote the set of $h \in \N \cap [0,u]$
such that the number of $x \in \csetuthree$
with $x \cvl \in (i-1,1]$
exceeds $R$ for at least two-thirds of the indices $i$ in $\llbracket h+1+A,h+3A \rrbracket$.
For $h \in M$, write 
$$
O_h \, = \, \big[ h+A - 3Kd^{1/2}r^{-1}, h + 3A + 3Kd^{1/2}r^{-1} \big] \, .
$$
Given a finite collection $\mc{C}$ of real intervals of equal length, a pairwise disjoint subset $\mc{I}$ may be selected such that 
$\sum_{I \in \mc{I}} \vert I \vert \geq \tfrac{1}{2} \big\vert \bigcup_{I \in \mc{C}}I \big\vert$.
We may thus let $U \subseteq M$ be chosen so that $\big\{ O_u: u \in U \big\}$ is a pairwise disjoint collection that satisfies
$$
 \sum_{u \in U} \vert O_u \vert \, \geq \, \tfrac{1}{2} \cdot \bigg\vert \bigcup_{m \in M} O_m \bigg\vert \, .
$$

Let $h \in U$ and write $K_h = \big\{ i \in \llbracket h+1+A, h + 3A \rrbracket : \chi_i \geq R \big\}$. We form a subset $J_h$ of $K_h$ by first including in it the minimum value in $K_h$, and then iteratively including every $\lceil 5K d^{1/2}\rceil$\textsuperscript{th} value of $K_h$ as these are encountered in increasing order. The constructed set $J_h$ has the properties that $j \in J_h$ implies that $\chi_j > R$; that the distance between consecutive elements of $J_h$ is at least $5K d^{1/2}$; and that $\vert J_h \vert \geq \tfrac{4A}{3(5K d^{1/2} +1)}$.

Let $u \in U$ and $j \in J_u$. Since $\chi_j \geq R$, Lemma~\ref{l.latint} permits us to find $\ell_j \in \Gamma_{j-1}$
such that  $\overline\chi_j \geq R(d^{1/2}r^{-1} +2)^{-1}$, where $\overline\chi_j$ equals the number of $K$-boxes in $\overlinecsetuthree$
of latitude $\ell_j$.

We {\em claim} that there is no coincidence among $\ell_j$ as $u$ and $j$ run over $U$ and $J_u$, and indeed that the pairwise distance between these quantities is at least two.
To confirm this, we begin by arguing that when $j \in J_u$ and $j' \in J_{u'}$ for $u,u' \in U$, equality $j = j'$ is possible only if $u = u'$, and, if $j \not= j'$, then $\vert j - j' \vert \geq 5K d^{1/2}$.
Indeed, if it happens that $j,j' \in J_u$ for a common value of $u$, then, if $j \not= j'$, $\vert j - j' \vert$ is indeed at least $5K d^{1/2}$.
If $j \in J_u$ and $j \in J_{u'}$ for distinct $u$ and $u'$, then the distances $d(j,O_u^c)$ and $d(j',O_{j'}^c)$ are both at least $3K d^{1/2}r^{-1}$, so that the disjointness of $O_u$ and $O_{u'}$ implies that $\vert j - j' \vert \geq 6K d^{1/2}r^{-1}$.

To arrive at the claim, we next argue that, if $k,k' \in \Z$ differ by at least $5K d^{1/2}$, then $\Gamma_k$ and $\Gamma_{k'}$ lie at distance at least two. To see this, note that, since $d \geq 3$ and $r \leq 1$,
the difference between $k$ and $k'$ is at least $2Kr(d^{1/2}r^{-1} + 2)$; from here, Lemma~\ref{l.latint} shows the assertion just made. 

We know that $\ell_j \in \Gamma_{j-1}$ and $\ell_{j'} \in \Gamma_{j'-1}$. Taking $k =j-1$ and $k' = j'-1$, note that we have confirmed the hypothesis that $d(k,k') \geq 5K d^{1/2}$, so that the assertion in the preceding paragraph yields the sought claim. 

Next we will show that
\begin{equation}\label{e.cupo}
 \Big\vert \bigcup_{h \in M} O_h \Big\vert \, \geq \, 2M/5 \, .
\end{equation}
To see this, set $\kappa = \vert \Z \cap O_h \vert$ for any $h \in M$. 
For $x \in \Z$, let $\phi_x = \big\vert \big\{ h \in M: x \in O_h \big\} \big\vert$. Since any $x\in \Z$ belongs to at most $2\kappa+1$ intervals $O_h$ with $h \in M$,
we see that $\max \{ \phi_x : x \in\Z \} \leq 2\kappa + 1$. Note then that 
$$
 \kappa \vert M \vert = \sum_{h \in M} \vert O_h \vert = \sum_{x \in \Z} \phi_x \leq (2\kappa+1) \vert {\rm supp}(\phi) \vert \leq (2\kappa +1 ) \Big\vert \bigcup_{h \in M} O_h \Big\vert \, ,
$$
since ${\rm supp}(\phi) = \cup_{h \in M} O_h$. Thus, $\sum_{h \in M} \vert O_h \vert \geq \tfrac{\kappa}{2\kappa + 1} \vert M \vert$. 
Since
$$
 \kappa = \big\vert \Z \cap \big[ h+A - 3K d^{1/2}r^{-1}, h + 3A + 3Kd^{1/2}r^{-1} \big] \big\vert
$$
$\kappa \geq 2A - 1 + 6Kd^{1/2}r^{-1} \geq 2A$, where the latter bound depends on $d \geq 3$, $K \geq 1$ and $r \leq 1$.
We confirm that (\ref{e.cupo}) holds when $A \geq 1$.

Note that 
$$
(2A + 6K d^{1/2} r^{-1}) \vert U \vert \, = \, \sum_{u \in U} \vert O_u \vert \, \geq \, \tfrac{1}{2} \cdot \bigg\vert \bigcup_{m \in M} O_m \bigg\vert \, \geq \, M/5 \, ,
$$
where~(\ref{e.cupo}) provides the latter bound.
Proposition~\ref{p.dn.gen} implies that, for each  $u \in U$ and $j \in J_u$, 
$$
 \Big\vert  \partial_{\rm ext} \overlinecsetuthree \cap \overline{N}_{\ell_j,\ell_j +1} \Big\vert \, \geq \, 2^{\tfrac{2-d}{d-1}} \big( \tfrac{R}{d^{1/2}r^{-1}+2} \big)^{\tfrac{d-2}{d-1}} \, .
$$
Since each pair of these latitudes differ by at least two, we find that the  the number of $K$-boxes in $\partial_{\rm ext} \overlinecsetuthree$
at latitudes $\big\{ \ell_{u,j}, \ell_{u,j} + 1: u \in U, j \in J_u \big\}$ is at least $2^{\tfrac{2-d}{d-1}} \big( \tfrac{R}{d^{1/2}r^{-1}+2} \big)^{\tfrac{d-2}{d-1}} \vert U \vert \tfrac{4A}{3(5K d^{1/2} + 1)}$.
But $\big\vert \partial_{\rm ext} \overlinecsetuthree \big\vert \leq  c R^{\tfrac{d-2}{d-1}} K^{-1} au$, which implies that 
$$ 
\vert U \vert \leq  c  R^{\tfrac{d-2}{d-1}} K^{-1} 2^{\tfrac{d-2}{d-1}} \big( \tfrac{d^{1/2}r^{-1}+2}{R} \big)^{\tfrac{d-2}{d-1} }  \tfrac{3(5K d^{1/2} + 1)}{4A} au \, .
$$
 We find that $\vert M \vert \leq 10(A + 3Kd^{1/2}r^{-1}) \vert U \vert$.  Using $A \geq 3K d^{1/2}r^{-1}$, we obtain
 $$
  \vert M \vert \leq 15 r^{-1} c (2R)^{\tfrac{d-2}{d-1}} K^{-1} \big( (d^{1/2}r^{-1} + 2)R^{-1} \big)^{\tfrac{d-2}{d-1}}(5K d^{1/2} + 1)au \, ;
 $$
 whence $\vert M \vert \leq c \beta u$,  where $\beta$ is a positive constant determined solely by $d$. We may then choose $c = \beta/4$ to ensure that $\vert M \vert \leq au/4$.
 Thus, the set of $h$ for which the northzone  condition $\vfsOne$ specifying $\veryfine_h$  is violated is at most $au/4$. Similarly, the same upper bound holds in regard to the southzone condition~$\fsTwo$.

In summary, the cardinality of the set of $h \in \N \cap [0,u]$ for which $\veryfine_h$ does not occur is at most~$au$. We have  proved 
 Lemma~\ref{l.nhigh}. \qed

\subsection{Simple joins are often found at fine levels: deriving Proposition~\ref{p.success}}\label{s.simplejoin}

To prove Proposition~\ref{p.success}, we need to gain a more explicit understanding of how the slide resample operates.
A first subsection presents Lemma~\ref{l.probform}, to this effect. In this lemma appears  a denominator or partition function. This is a sum over head locations $v \in \fubuthree$ 
of what we have called connection chances: conditional probabilities 
for cluster formation given, roughly speaking, the construction of the downcluster, and of the shift of the upcluster so that its head lies at $v$.
In Subsection~\ref{s.reduction}, we reduce Proposition~\ref{p.success} to Proposition~\ref{p.finesum}, which offers an upper bound on the partition function. 
The proof of the latter result, which involves coarse-grained percolation ideas seen in Section~\ref{s.coarse}, is the subject of Subsection~\ref{s.finesum}.

\subsubsection{An explicit description of the conditional law in the slide resample}
Two definitions concerning the resample permit us to  present Lemma~\ref{l.probform}, which rigorously expresses the notion advanced in Section~\ref{s.slideintorenewal}
that any given shift value is selected under resampling with probability proportional to its connection chance. 

Recall the random shift-equivalence class $[\upwell]$ specified under $\PP$ given $\clusteruthree$ by Definition~\ref{d.equivalence}.
\begin{definition}\label{d.displace}
Let $v \in \fubuthree$ (so that $v \cvl \in [u,u+3r]$).
On the event $\clusteruthree$, we denote by $\upwell^v$ the unique element of $[\upwell]$ 
whose head equals~$v$. Note that
the equivalence class  $[\upwell]$ is equal to $\big\{ \upwellv : v \in \fubuthree \big\}$. Under the percolation law $\PP$ conditionally on the occurrence of $\clusteruthree$, the cluster $\csetuthree$ is well defined but random; thus, so is $\upwell$. We may thus define, under the event $\clusteruthree$, the random vector $V \in \fubuthree$ such that $\upwell$ is equal to $\upwell^V$. 

\end{definition}

\begin{definition}\label{d.prime}
Let $\Omega$ denote the power set of the collection of nearest-neighbour edges in~$\Z^d$. The set $\Omega$ inherits a topology by identifying this set with a unit interval in $\R$, and we set $\mc{B}$
equal to the associated Lebesgue $\sigma$-algebra on $\Omega$. 
The law $\PP$ of  bond percolation on the graph $(\Z^d,\sim)$,  in which each edge is independently declared {\em open} with probability~$p$, and {\em closed} in the opposing event, will be treated as being defined on the  measurable space $(\Omega,\mc{B})$.

 We select a second copy $(\Omega',\mc{B}',\PP')$ of this probability space $(\Omega,\mc{B},\PP)$. 
The product space  $(\Omega \times \Omega',\mc{B} \times \mc{B'},\PP \times \PP')$ describes two independent bond percolations on $\Z^d$ of parameter~$p$.
The first may be called the {\em original} percolation and the second, associated to the prime variable, the {\em auxiliary} percolation.

 
Let $\Omega''$ denote a further copy of the power set of nearest-neighbour edges in~$\Z^d$. For each $v \in \fubuthree$, we specify 
a further configuration $\omega''_v \in \Omega''$ as a function of $(\omega,\omega') \in  \Omega \times \Omega'$. 

An edge is called $\omega$-open or $\omega$-closed if it is, or is not, an element of $\omega$; naturally the same usage may be made for $\omega'$ and $\omega''_v$. 
Set $\omega''$ to be a definite choice---the identically closed configuration, say---if $\omega \not\in \clusteruthree$. Suppose now that $\omega \in \clusteruthree$.
In what follows, every random variable, such as $\upwellv$ or $\downwell$, is understood to be specified by the configuration~$\omega$. 
To specify the configuration $\omega_v''$, we classify edges in $\Z^d$.
\begin{enumerate}
\item Elements of $E(\upwellv) \cup E(\downwell)$ are called {\em discovered open}. 
\item Elements of $E(\downward_h) \cap \partial_E(\downwell)$ and $E\big( U_{\hadj} \big) \cap \partial_E(\upwellv)$ are {\em discovered closed}. 
\item Every further nearest-neighbour edge in $\Z^d$ is called {\em undiscovered}.
\end{enumerate}
Here, $\partial_E(G)$, for $G \subseteq E(\Z^d)$, denotes the set of edges in $\Z^d$ that do not belong to $G$ but each of which is incident to some element of $G$.

Each discovered edge is $\omega_v''$-open precisely when it is $\omega$-open.
 Each undiscovered edge is $\omega_v''$-open precisely when it is $\omega'$-open. 
Note that, in so specifying $\omega_v''$, discovered open edges are $\omega_v''$-open, and discovered closed edges are $\omega_v''$-closed.
\end{definition}
Let $\uclust$ and $\dclust$ denote possible values for $[\upwell]$ and $\downwell$ as $\omega$ ranges over $\clusteruthree$. Set $R^v_{\uclust,\dclust}(\omega')$ equal to the value of $\omega''_v$
given $\clusteruthree$, $[\upwell] =\uclust$ and $\downwell = \dclust$. The dependence is on $\omega'$ because it is the values of $\omega''_v$ on undiscovered edges which dictate the value of $\omega''_v$
in this circumstance.

 Set $\R^v_{\uclust,\dclust}: \mc{B} \to [0,1]$, $\R^v_{\uclust,\dclust}(A) = \PP' \big( \omega' \in \Omega' : R^v_{\uclust,\dclust}(\omega') \in A \big)$.


In the next result, we use this framework to express conditional probability given a subevent of~$\clusteruthree$. Just as random variables under $\PP$ are understood to be given by the $\omega$-configuration, so are events: for example, in the left-hand side of Lemma~\ref{l.probform}.
\begin{lemma}\label{l.probform}
Let $B,H \subseteq \Omega$ denote two $\mc{B}$-measurable events, with $\PP (B \vert \clusteruthree) > 0$. 
Let $(\uclust,\dclust)$ be a value that $([\upwell],\downwell)$ assumes with positive probability under $\PP$ given $B \cap \clusteruthree$.
Then  
\begin{eqnarray*}
& & \PP \, \Big( H \, \Big\vert \,  [\upwell] = \uclust, \downwell = \dclust , B , \clusteruthree \Big) \\
& = & \frac{\sum_{v \in \fubuthree} \R^v_{\uclust,\dclust} \big(H \cap B \cap \clusteruthree \cap \{ V = v\} \big)}{\sum_{v \in  \fubuthree} \R^v_{\uclust,\dclust} \big(B \cap \clusteruthree  \cap \{ V = v\}  \big)} \, ,
\end{eqnarray*}
where $V$ is the random vector specified under $\PP (\cdot \, \vert \, \clusteruthree)$ in Definition~\ref{d.displace}.
\end{lemma}
{\bf Proof.}
Let  $G \subseteq \Omega$ be an arbitrary $\mc{B}$-measurable event.
Take $v \in  \fubuthree$, and consider the event $\clusteruthree \cap G \cap \{ \upwell = \mathfrak{u}^v \} \cap \{ \downwell = \dclust \}$.
Does this event occur? We may find the answer in two stages. First, we check a necessary condition, an event that we label~$N$: namely, that elements of $E(\mathfrak{u}^v) \cup E(\dclust)$ are $\omega$-open
and elements of $E(\downward_u) \cap \partial_E(\dclust)$ and $E\big( U_{\hmax- (u-h)} \big) \cap \partial_E(\mathfrak{u}^v)$ are $\omega$-closed. 
If this test is passed, we move to the second stage, and realize all further edges in~$E(\Z^d)$. Check that the event $\clusteruthree \cap G$ occurs. If this happens, then $\downwell$ must equal~$\dclust$, but $\upwell$
is not yet assuredly $\mathfrak{u}^v$. We need to confirm the latter condition, and this is done by checking (at the end of the second stage) that the random vector $V$ equals $v$.

Writing $\alpha_v = \PP (N) > 0$,
note that $\alpha_v$ and $\R^v_{\uclust,\dclust} \big(G \cap \clusteruthree  \cap \{ V = v\}  \big)$ equal the respective probabilities that the checks at the first and second stages are passed.
In fact, $\alpha_v$ may be written~$\alpha$ because translation invariance of  percolation in the upspace implies that this quantity has no dependence on $v$.
What we have learnt is that, for $v \in  \fubuthree$,
\begin{equation}\label{e.resampleclaim}
\PP \, \Big( G ,  \upwell = \mathfrak{u}^v , \downwell = \dclust , \clusteruthree \Big) \, = \, \alpha \cdot \R^v_{\uclust,\dclust} \big(G \cap \clusteruthree  \cap \{ V = v\}  \big) \, .
\end{equation}
 Now set  $G = B \cap \clusteruthree$ and sum~(\ref{e.resampleclaim}) out over $v \in  \fubuthree$. The resulting value is positive by our assumptions. 
 Then take  $G$ equal to $H \cap B \cap \clusteruthree$ and sum out similarly. Dividing the latter of the resulting equalities by the former,    we obtain Lemma~\ref{l.probform}.
\qed

\subsubsection{Reducing Proposition~\ref{p.success} to a partition function upper bound}\label{s.reduction}

We reduce Proposition~\ref{p.success} to two results the first of which we then prove. 
 \begin{lemma}\label{l.qfine}
 There exists $c > 0$ such that $\sum_{v \in  \fubuthree} \R^v_{[\upwell], \downwell}(\fine, V=v) \geq c$ for $\omega \in \veryfine$. Indeed,  for such $\omega$ there exists $x \in  \fubuthree$ for which 
 $$
 \R^x_{[\upwell], \downwell}\big(\fine \cap \success \cap \{ V =x \} \big) \geq c \, .
 $$
\end{lemma}
(The value of $h$ that specifies $\fine = \fine_h$ is given in Proposition~\ref{p.success} and we will for now omit it from this notation.)

The second result is our partition function upper bound, which we prove in the next subsection.
\begin{proposition}\label{p.finesum}
Suppose  $A/K$ is bounded below by a suitably high constant multiple of $\log u$.  
Then 
there exists $C > 0$ such that $\sum_{v \in  \fubuthree} \R^v_{[\upwell],\downwell}(\fine) \leq C$ whenever $\omega \in \fine$.
\end{proposition}

 {\bf Proof of Proposition~\ref{p.success}.} Take $H = \success$ and $B = \fine_u$ in Lemma~\ref{l.probform}.
 We seek to bound the right-hand numerator below, and its denominator above.
 The assumption that $h$ and $u-h$ exceed a large constant multiple of $\log u$ permits a choice of $A$ of the same order, with $\min \{ h,u-h\} > 3A$.
 This choice of $A$ allows us to apply  Proposition~\ref{p.finesum}, to bound above the denominator. Bounding the numerator below via
 Lemma~\ref{l.qfine}, we obtain the sought assertion. \qed

\begin{lemma}\label{l.finecard}
 The event $\fine$ entails that the vertex sets $\Base(\upwell)$ and $\Top(\downwell)$ have cardinality at most~$R$.
\end{lemma}
{\bf Proof.} By Lemma~\ref{l.geometrybasics}(1), every element of $\Top(\downwell) \subset \vertbdry_h$ is a vertex in $\csetuthree$ whose height lies in $[h-r,h]$ for $r = e_d \cvl$.
Recalling that $\hadj = \hmax - (u-h)$,
note that every element of  $\Base(\upwell)$  lies in $\vertbdry_{\hadj}$, 
so that the just quoted result implies that each such element has height in $[\hadj - r,\hadj]$.
Since $r \leq 1$, 
we use properties  ${\bf T}'_1(R)$  and ${\bf T}'_2(R)$ (with $i=1$) of the event $\fine$ to 
obtain Lemma~\ref{l.finecard}. \qed

{\bf Proof of Lemma~\ref{l.qfine}.} 
Select a vector $x \in  \Z^d$ such that some vertex in $\Base(\upwellx)$ is displaced from some vertex $v \in \Top(\downwell)$ by $2 e_d$.
(If for given $x$ there is a choice of admissible $v$, we pick  $v$ from the admissible set to be the lexicographically minimal vertex among those of minimum height.) 
The height $x \cvl$ of $x$ equals the sum of the heights of $\downwell$ and $\upwell$ and the value $2 e_d \cvl$.
By Lemma~\ref{l.updowncluster}, $x \cvl  \in [u + e_d \cvl, u + 3 e_d \cvl] \subset [u,u+3r]$. We have confirmed that $x \in \fubuthree$.
 
Consider the event $E$ that
\begin{itemize}
\item the edges $[v,v+e_d]$ and $[v+e_d,v+2e_d]$ are $\omega'$-open, while all other edges incident to $v+e_d$ are $\omega'$-closed, as are all other undiscovered edges that are incident to  $\Base(\upwellu) \cup \Top(\downwell)$.
 \end{itemize}
 For each $\omega \in \fine$, 
 $$
 \R_{[\upwell],\downwell}^x ( E ) \geq p^2 (1-p)^{2(d-1)} (1-p)^{(2d-1)\cdot 2(R-1)}
  \, .
 $$
 Indeed, the right-hand term bounds the probability of $E$: $p^2$ opens the two $e_d$-oriented edges, $(1-p)^{2(d-1)}$ closes the other edges incident to the shared endpoint of this pair, and $(1-p)^{(2d-1)\cdot 2(R-1)}$ makes the other needed closures in view of Lemma~\ref{l.finecard}. 
 
 When the configuration $R^x_{[\upwell],\downwell}(\omega')$ realizes the event $E$, $\midwell$ contains precisely two edges, these of the form $[v,v+e_d]$ and $[v-e_d,v]$ for some given $v \in \Z^d$ with $v \cvl =h$; and
 the graph $\midfragment$ is empty.
 Thus,  $R^x_{[\upwell],\downwell}(\omega')$ realizes the event $\success$ in this circumstance. Moreover, we {\em claim} that, when $\omega \in \fine$ and  $R^x_{[\upwell],\downwell} (\omega') \in E$, we have that $R^x_{[\upwell],\downwell}(\omega') \in \fine \cap \{ V = u \}$.
 To confirm that $R^x_{[\upwell],\downwell}(\omega') \in \fine$,
 we need to check the tropical, subtropical and global conditions for $\fine$. To do so, 
 note that any vertex in $\csetuthree$ is an element of precisely one of the upcluster, the midcluster and  the downcluster; in the transition from $\omega$ to $R^x_{[\upwell],\downwell}(\omega')$, the upcluster may have shifted; the midcluster has collapsed, to be confined to its minimum vertex cardinality of three, with vertices $v - e_d$, $v$ and $v+e_d$; and $\downwell$ has remained unaltered. 
 The various upper bounds on vertex cardinalities in the three conditions are thus seen to pass from $\omega$ to $R^x_{[\upwell],\downwell}(\omega')$. Thus, the latter configuration indeed verifies $\fine$.  To confirm that $V =x$, it is enough to check that $\upwell$ under  $R^x_{[\upwell],\downwell}(\omega')$ is equal to $\upwellx(\omega)$. There is a unique open connected component of $\clusteruthree$ in $U_{\hmax- (u-h)}$ under this configuration, because the only element of $E(\clusteruthree)$
that our construction permits that has one endpoint in $\vertbdry_h$ and the other in  $V(U_h) \setminus \vertbdry_h$
is the edge $[v,v+e_d]$. Since this unique component equals  $\upwellx(\omega)$, we indeed find that $V =x$. This completes the proof of the claim. 
We learn that, for  $\omega \in \fine$, 
 $$
  \R_{[\upwell],\downwell}^x \big(  \fine \cap \success \big) \, \geq \, p^2 (1-p)^{2(2dR - d -R)} \, .
 $$    
The latter statement of Lemma~\ref{l.qfine} has been derived; the former is a trivial consequence. \qed


\subsubsection{Proving the partition function upper bound}\label{s.finesum}

Here we prove Proposition~\ref{p.finesum}. To argue that $\R^v_{[\upwell],\downwell}(\fine)$ is small when $v \in  \fubuthree$ is large, we will note that the tropical conditions in $\fine$
entail that $\upwell^v$ and $\downwell$ have a sparse presence near the equator. Thus, if $v$ is large, these graphs visit the tropical region at significant horizontal separation. The occurrence of $\fine \subset \clusteruthree$ under $\R^v_{[\upwell],\downwell}$ thus requires the formation of new, independent, cluster boundaries to bridge the gap, either horizontal or vertical, between $\upwell^v$ and $\downwell$.
The cost of these cluster boundaries is high, so that   $\R^v_{[\upwell],\downwell}(\fine)$  decays rapidly as the norm of $v$ rises.   

To show that such costs are indeed high, we will use the coarse-graining of percolation seen in Section~\ref{s.coarse}.
 Indeed, we let $\const  \in \N$ be a positive parameter, 
and recall from  that section that the collection of $K$-boxes tesselates~$\Z^d$, and has an index set which may be naturally identified with~$\Z^d$.
Recall that a path in $\const$-boxes makes nearest-neighbour steps given this identification;
and also that the notion of $*$-neighbour specified at the start of Section~\ref{s.coarse} passes naturally from $\Z^d$ to the set of $K$-boxes.

The notation $\overline{A}$ comes from Definition~\ref{d.plustwominustwo}.  Since $\midwell$ is nearest-neighbour connected, so is $\overline\midwell$; the latter set thus lies in
a single
 $*$-connected component of~$\overline\residualwell$. We call this component~$\mc{H}$.
 Using Definition~\ref{d.exterior}, we set $\surfacek = \partial_{\rm ext} \mc{H}$.
 \begin{lemma}\label{l.surfacek}
\leavevmode
\begin{enumerate}
\item The collection
$\surfacek$ is a nearest-neighbour connected set of $\const$-boxes.
\item Every element of $\surfacek$ is disjoint from $\residualwell$.
\item Every infinite path in $\const$-boxes lying in $\overline{S}_{0,n}$ that begins at a box that intersects $\midwell$ visits $\surfacek$.
\end{enumerate}
\end{lemma}
{\bf Proof: (1).}  This is due to~\cite[Theorem~$4$]{Timar2013}.

{\bf (2,3).} These are direct from definitions. \qed
\begin{lemma}\label{l.surfacebad}
Suppose that $\clusteruthree$ occurs. Let $B$ be a $K$-box that lies in $\surfacek$ and 
whose vertex set is disjoint from the union of the vertex sets of $\upwell$ and $\downwell$. 
Then $B$ is a bad $\const$-box. 
\end{lemma}
\noindent{\bf Proof.}
Note that $B$ is disjoint from $\csetuthree$.
Indeed, 
any vertex in $\csetuthree$ lies in at least one of $\downwell$, $\upwell$ and $\residualwell$. A vertex in $B$ may not lie in any of these graphs:
 the first and second are excluded by hypothesis; the third, by Lemma~\ref{l.surfacek}(2).

Let $D$ denote a $\const$-box $*$-neighbour of $B$ that intersects $\residualwell \subseteq \csetuthree$. Lemma~\ref{l.bdbad} now implies that $B$ is a bad $\const$-box. \qed

 It is useful to consider slight fattenings and shrinkings of collections of $K$-boxes. The next result and its proof utilize the notions of $K$-box latitude from Definition~\ref{d.boxlatitude}, and of the $K$-box   fattening $\overline{A}[2]$ and shrinking $\overline{A}[-2]$  from Definition~\ref{d.plustwominustwo}. 
\begin{lemma}\label{l.latcard}
Suppose that $\fine \subset \clusteruthree$ occurs. 
\begin{enumerate}
\item A $K$-box lying in  $\overline\downwell[2]$  has latitude at most $5d/2$.
\item  Such a box lying in 
$\overline\midzone[-2]$ has latitude at least  $- (2K)^{-1}( A +K d^{1/2})$.
\item
For $(A - 5K d^{1/2})(2Kr)^{-1}  \geq  i \geq -5d/2$, the number of $K$-boxes of latitude $-i$ lying in $\overline\downwell[2]$ is at most $5^d \big( 6K d^{1/2} + 2Kr + \vert i \vert \big)^4 R$. 
\item 
For $(A - 5K d^{1/2})(2Kr)^{-1}  \geq  j \geq -5d/2$, the number of $K$-boxes of latitude $j$ lying in $\overline\upwell[2]$ is at most  $5^d \big( 6K d^{1/2} + 2Kr + \vert j \vert \big)^4 R$.
\end{enumerate} 
\end{lemma}
{\bf Proof: (1).} Since $\downwell$ is a subgraph of $L_h$, $v \in V(\downwell)$ satisfies $v \cvl \leq h$.  
If $x$ is the central index of any element $\overline\downwell$, we thus have $x \cvl \leq h+Kd^{1/2}$; and thus, if $w$
  is the central index of any element $B \in \overline\downwell[2]$ (so that $\max_{i \in \intint{d}}\vert x_i - w_i \vert \leq 4K$), $w \cvl \leq h + K (1 + 2^{2})d^{1/2}$. We find then that 
$$
 \lat(B) = \big\lfloor \tfrac{w \cvl - h}{2K r} \big\rfloor \leq \frac{5K d^{1/2}}{2Kr} \leq 5d/2
$$
where $r = e_d \cvl \geq d^{-1/2}$ by axes ordering~(\ref{e.aoc}).

{\bf (2).} 
Let $2K w + \llbracket -K,K \rrbracket^d \in \overline{\tropics}[-2]$.
Let $B$ be the $K$-box with central vertex $2Kw - 4K e_d$. Then $B \in \overline\tropics$, so $B$ intersects $U_{h-A}$, with the argument for Lemma~\ref{l.geometrybasics}(1)
yielding that the height of every vertex in $U_{h-A}$ being at least $h-A-r$. Thus,
 $(2Kw - 4K e_d) \cvl \geq h - A -r - K d^{1/2}$.
Writing $B' = 2Kw + \llbracket -K,K \rrbracket^d$,
$$
 \lat(B') =  \big\lfloor \tfrac{w \cvl - h}{2K r} \big\rfloor \geq \frac{(2K-1)r - A - K d^{1/2}}{2Kr} \geq - \frac{A+ d^{1/2}K}{2K}
$$ 
where we used $K \geq 1/2$ and $r \leq 1$ in the latter bound.

{\bf (3).} We take $i \geq 0$; the case where $i \in \llbracket - \lfloor 5d/2 \rfloor ,-1 \rrbracket$ is similar. Let $B = 2K w + \llbracket -K,K \rrbracket^d$ be a $K$-box of latitude $-i$ lying in $\overline\downwell[2]$. Let $x \in V(\downwell) \cap V\big(\boxbig(B)\big)$.  Note that $x \cvl - h = 2Kw \cvl - h + \phi$ with $\vert \phi \vert \leq 5K d^{1/2}$.
Since $2K w \cvl \in h + 2Kr \cdot [i,i+1)$,
we find that 
$$
 x \cvl - h \in \big[ 2Kri - 5K d^{1/2} , 2K r(i+1) + 5K d^{1/2} \big) \, .
$$
We claim that $h \geq x \cvl \geq h - A$. The first bound is due to $x \in V(\downwell)$. To confirm the second, note that $\big\lfloor \tfrac{2K w \cvl - h}{2Kr} \big\rfloor = -i$
implies that $2K w \cvl -h \geq h - 2K ri$. But then $x \cvl \geq 2K w \cvl - 5Kd^{1/2}$ and 
\begin{equation}\label{e.xbound}
 x \cvl \geq h - 2Kri - 5K d^{1/2} \, ,
\end{equation}
so that the sought bound $x \cvl \geq h - A$ follows from the hypothesised condition
$i \leq \tfrac{A - 5 K d^{1/2}}{2Kr}$.
For $j \in \intint{A}$, set $\Phi_j = \big\{ x \in V(\downwell): x \cvl \in (h-j,h-j+1] \big\}$. In this way, condition  ${\bf T}'_2(R)$ in $\fine_h$ asserts that $\vert \Phi_j \vert \leq j^3 R$ for $j \in \intint{A}$.
Since $h \geq x \cvl \geq h - A$, there exists $j \in \intint{A}$ such that $x \in \Phi_j$.
Moreover, if $x \in \Phi_j$, then $x \cvl \leq h-j+1$; alongside~(\ref{e.xbound}), we find that $j \leq 5K d^{1/2} + 2K ri +1$.
Which is to say: $x$ is necessarily a member of a set $\Phi_j$ of cardinality at most $j^3 R$
for an index $j$ which is at least one and at most $5K d^{1/2} + 2Kr +i$.
Moreover, given $x$, there are at most $5^d$ compatible choices of the $K$-box $B$. We thus see that the number of choices of~$B$ is at most
$$
 5^d \sum_{j=1}^{6k d^{1/2} + 2Kr + i} j^3 R \, ,
$$
whose right-hand side is bounded above by $5^d \big( 6K d^{1/2} + 2Kr + i \big)^4 R$.
That is, we have found that
$$
 \Big\vert \overline{N}_{-i} \cap \overline{\downwell}[2] \Big\vert \, \leq \, 5^d \big( 6K d^{1/2} + 2Kr + i \big)^4 R \, .
$$

{\bf (4).} Evident changes of notation in the preceding proof yield the argument. \qed

\begin{definition}
Suppose given a finite path $P = \big( P_i: i \in \llbracket 0,k \rrbracket \big)$ whose ending point~$P_k$ coincides with the starting point~$Q_0$ of a path $Q = \big( Q_i: i \in \llbracket 0,m \rrbracket \big)$. 
Here, $k \in \N$ and $m \in \N \cup \{ \infty \}$, with $\llbracket 0, \infty \rrbracket = \N$ understood. The concatenation of $P$ and $Q$ is the sequence $P \circ Q: \llbracket 0, k+m \rrbracket \to \Z^d$
with $(P\circ Q)_i = P_i$ for $i \in \llbracket 0, k \rrbracket$ and  $(P\circ Q)_i = Q_{i-k}$ for $i \in \llbracket k, k+m \rrbracket$. 
\end{definition}

{\bf Proof of Proposition~\ref{p.finesum}.} Recall from Definition~\ref{d.fvf} that $\fine \subset \clusteruthree$ and that $V \in  \fubuthree$ is a random head location specified by Definition~\ref{d.displace} when $\clusteruthree$ occurs.  Let $v \in  \fubuthree$.
In the ensuing paragraphs, we will specify some aspects of the geometry of the cluster under the event $\fine  \cap \{ V=v \}$; we thus suppose that this event occurs. 
Recall that $\midwell$ contains an open path that connects a vertex in $\Top(\downwell)$ to a vertex in $\Base(\upwell)$. Let $Q$ denote the lexicographically minimal among the set of such open paths. Suppose that $Q$ runs from $w \in \Top(\downwell)$ to $w' \in \Base(\upwell)$. Let $Q^-$ denote the lexicographically minimal path in $\downwell$ from $\base(\csetuthree) = 0 \in \Z^d$ to $w$; and let $Q^+$ denote the similarly minimal path in $\upwell$ from $w'$ to $\head(\csetuthree) = \head(\upwell)$. Denote by $P_1$ the concatenation $Q^- \circ Q \circ Q^+$.
Note that $P_1$ is an open path in the cluster $\csetuthree$ that begins at $\base(\csetuthree)$ and ends at $\head(\csetuthree)$ and that runs successively through the down-, mid- and up-clusters. 



Prefix to $P_1$ an infinite path $P_0$ of consecutive positively oriented $e_1$-parallel edges that runs directly into $\base(\csetuthree)$, in a direction that we may call `rightwards' (and which is increasing for height because $e_1 \cvl \geq 0$ by axes ordering). Postfix such a path~$P_2$ that starts at $\head(\csetuthree)$ and that also runs directly rightwards. The resulting path $P = P_0 \circ P_1 \circ P_2$ is a bi-infinite path in $\Z^d$ that runs directly rightwards through heights at most zero into an open subpath that runs from $\base(\csetuthree)$ and takes a passage through $\midwell$ to end at $\head(\csetuthree)$, and then pursues a directly rightward trajectory through heights at least~$u$.  We now wish to specify $P[\const]$ to be the $\const$-box valued bi-infinite path that visits the $\const$-boxes consecutively visited by $P$. 
Since $\const$-boxes overlap (albeit only on their boundaries), this specification needs to be clarified. Recall that elements  $x = 2K u$, $u \in \Z^d$, of the $\const$-box index space index  $K$-boxes $x + \llbracket -K,K \rrbracket^d$. 
\begin{definition}
Let ${\rm ind} : \Z^d \to \Z^d$ be the function that sends any element of $2K u + \llbracket -K,K-1 \rrbracket^d$ to the index $2K u$. 
To the bi-infinite path $P:\Z \to \Z^d$, we may associate the mapping $\hat{P}$ from $\Z$ to the $\const$-box index space via $\hat{P}(i) = {\rm ind} \big( P(i) \big)$ for $i \in \Z$. 
Consider the bi-infinite path~$S$ in $\const$-boxes $i \to \hat{P}(i) + \llbracket -\const,\const \rrbracket^d$, $i \in \Z$. Let $R \subset \Z$ denote the set of indices $i$ at which $S(i) \not= S(i-1)$ (and note that $R$ has infinite intersection with both $\N$ and $\Z \setminus \N$).
Let $\rho: \Z \to R$ be an arbitrary increasing bijection. We then set $P[\const]$ to be the bi-infinite path  in $\const$-boxes $i \to S(\rho(i))$, $i \in \Z$.
\end{definition}
Note that indeed $P[\const]$ is a nearest-neighbour path in $\const$-boxes.
This path visits $\midwell$ in correspondence with the subpath $Q$ of $P$. It is Lemma~\ref{l.surfacek}(3) that permits the next definition.
\begin{definition}\label{d.twobs}
Let $B_1$ denote the last element of $\surfacek$ along $P[\const]$ before the first visit of $P[\const]$ to $\midwell$; and let~$B_2$ denote the first element of $\surfacek$ encountered along $P[\const]$ after the final visit of this path to $\midwell$. 
\end{definition}
Let $\Long^+ \subseteq \fine  \cap \{ V=v \}$ denote the event that   $\surfacek$ contains a $\const$-box path that starts at a box 
in $\overline{\rm Equator}[2]$
and that ends at a box in $\overline{\northzone}[2]$ with no element in the path belonging to $\overline\upwell[2]$.
Let $\Long^- \subseteq \fine  \cap \{ V=v \}$ denote the event that 
 $\surfacek$ contains a $\const$-box path  that also starts at a box in $\overline{\rm Equator}[2]$
 and that ends at a box that intersects $\overline\southzone[2]$, with no element in the path belonging to
  $\overline\downwell[2]$. 
  
  Set $\Long = \Long^+ \cup \Long^-$ and  $\Short = \fine  \cap \{ V=v \} \cap \Long^c$.
   

\begin{lemma}\label{l.short}
Suppose that $\Short$ occurs. At least one of two alternatives holds:
\begin{enumerate}
\item there exists a $K$-box in $\overline{\rm Tropics}[-2]$ that lies in both  $\overline\downwell[2]$ and $\overline\upwell[2]$; or
\item there exists a $K$-box path~$\Psi$ in  $\surfacek \cap \overline{\rm Tropics}[-2] \cap \overline\downwell[2]^c \cap \overline\upwell[2]^c$ whose first element is a nearest neighbour of an element of   $\overline\downwell[2]$ and whose final element is a nearest neighbour of an element of  $\overline\upwell[2]$. 
\end{enumerate}
\end{lemma}
Lemma~\ref{l.short} will permit us to specify a $K$-box path $\Psi$ that plays the counterpart role in rigorous analysis that the path ${\rm Connect}$ through virgin territory did in heuristic discussion in Section~\ref{s.slideintorenewal}. Indeed, 
when the latter alternative in the lemma occurs, we will set $\Psi$ equal to  some thus named path of minimum length, chosen in some definite manner, and $D[\Psi]$ and $U[\Psi]$
for the elements of  $\overline\downwell[2]$ and $\overline\upwell[2]$ that neighbour the start and end of $\Psi$.
Note also that, for the  event in question, we admit the possibility that $\Psi$ is empty. In this case, the condition amounts to the existence 
of nearest neighbour $K$-boxes  $D[\Psi]$ and $U[\Psi]$
 in   $\overline{\rm Tropics}[-2]$ such that  $D[\Psi] \in \overline\downwell[2]$ and $U[\Psi] \in \overline\upwell[2]$. In fact, even the event in the first alternative in Lemma~\ref{l.short} may be subsumed in this framework, because we may take $D[\Psi]$ and $U[\Psi]$ to be equal to the $K$-box posited in this case.

{\bf Proof of Lemma~\ref{l.short}.} Recall the $\const$-boxes $B_1$ and $B_2$ in Definition~\ref{d.twobs}.  Note that $B_1$ and $B_2$ are elements of $\surfacek$. 
 Suppose that  these two boxes are equal: $B_1 = B_2 = B$. We will argue that the first alternative presented in the lemma holds.  Since $B_1$ intersects either $P_0$ or $\downwell$, it also intersects $L_h$. Likewise, since $B_2$ intersects either $P_2$ or $\upwell$, it intersects $U_h$. Thus $B$ intersects $\equator$.  
Since $u-h$ and $h$ exceed $K$, $B$ can intersect neither $P_2$ (a path at height at least $u$) nor $P_0$ (a path at height at most zero). Thus $B=B_1 \in \overline\downwell \subset \overline\downwell[2]$ and $B = B_2 \in \overline\upwell \subset \overline\upwell[2]$.  And indeed $B$ must lie in $\overline{\rm Tropics}[-2]$, because $A - 3K > K$. Thus we confirm the first alternative in the lemma.

Suppose now that $B_1$ and $B_2$ are not equal.  
  Since $\surfacek$ is a connected set, we may find a $\const$-box path valued in $\surfacek$ that starts at $B_1$ and ends at $B_2$.   Select such a path~$\Gamma$ of minimal length, chosen in some definite manner.  We will argue that there exists a box $B_3$ on $\Gamma$ that lies in $\overline\downwell[2] \cap \overline\tropics[-2]$ such that the subpath of $\Gamma$ strictly beyond $B_3$
  makes no visit to  $\overline\downwell[2]$ and does not leave $\overline\tropics[-2]$ via its south side. To establish this, note that if $\Gamma$ begins (at $B_1$) to the south of  $\overline\tropics[-2]$, or if it leaves  $\overline\tropics[-2]$  via the south side, then it will, after reentering  $\overline\tropics[-2]$, make a visit to  $\overline\downwell[2]$ before ending at $B_2$: this is because we are supposing that $\Short$ occurs and $B_2$ lies in, or to the north of, $\overline\equator[2]$. If we take $B_3$ to be the last visit of $\Gamma$ to  $\overline\downwell[2]$, it thus has the sought property. 
 
 Now consider the reversal $\Gamma^\downarrow$ of $\Gamma$ from $B_2$, stopped at $B_3$. We assert that there is a box $B_4$ in $\Gamma^\downarrow$  that lies in $\overline\upwell[2] \cap \overline\tropics[-2]$ such that the subpath of  $\Gamma^\downarrow$ strictly beyond $B_4$
  makes no visit to  $\overline\upwell[2]$ and does not leave $\overline\tropics[-2]$ via its north side. Indeed, and similarly to above, if $\Gamma^\downarrow$ begins to the north of 
 $\overline\tropics[-2]$, or if it leaves  $\overline\tropics[-2]$  via the north side, then it will, after reentering  $\overline\tropics[-2]$, make a visit to  $\overline\upwell[2]$ before ending at $B_3$. Naturally we then take $B_4$ to be the last visit of  $\Gamma^\downarrow$  to  $\overline\upwell[2]$.
 
 It may be that $B_3 = B_4$: the first case in the lemma then holds. If not,
 we take $\Psi$ to be the subpath of $\Gamma$ strictly between $B_3$ and $B_4$.  This path remains in  $\overline\tropics[-2]$ and thus verifies the second case in the lemma. This completes the proof of Lemma~\ref{l.short}. \qed
 
 The projected distance $\dproj$ between two $\const$-boxes $B_1$ and $B_2$ is given by translating the two boxes in the $e_d$-direction so that they
 have latitude zero, and setting $\dproj(B_1,B_2)$ to be the $\ell_1$-distance between these translates.
 
Let $i,j \in \Z$ be at least  $-(d/2 +2)$ and let $k \in \N$.
We set  $\fv[i,j,k]$ to be the event that $\Short$ occurs alongside the latter alternative in Lemma~\ref{l.short}, with the $K$-box path~$\Psi$ in that alternative being such that $D[\Psi]$ has latitude $-i$; $U[\Psi]$ has latitude~$j$; and $\dproj \big( D[\Psi] , U[\Psi]\big) = k$.  The event that the first alternative in Lemma~\ref{l.short} obtains, in which $D[\Psi]$ and $U[\Psi]$ may be viewed as being equal, can in fact be included in this classification, by taking $k=0$ and $i=-j$.
\begin{lemma}\label{l.rfprob}
With $\e \in  \big(0,1/(2d-1) \big)$ denoting a parameter that may be chosen arbitrarily small by making a sufficiently high choice of $\const \in \N$, we have that, for $\omega \in \fine$,
\begin{enumerate}
\item and for $v \in \fubuthree$, $i,j \in \llbracket -  \lfloor d/2 +2 \rfloor , \lfloor A/\const \rfloor \rrbracket$ and $k \in \N$,
$$
\R^v_{[\upwell],\downwell}  \big( \fv[i,j,k] \big) \leq 
5^d \big( 6K d^{1/2} + 2Kr + \vert i \vert \big)^4   R (2d-1)^{\vert i+j \vert +k} \e^{\vert i+j \vert +k} \, ;
$$
\item for such $(i,j,k)$, the set of $v \in \fubuthree$ for which 
$\R^v_{[\upwell],\downwell}  \big(  \fv[i,j,k] \big)$ is non-zero has cardinality at most
$$
 50^d \big( 6K d^{1/2} + 2Kr + \vert i \vert \big)^4 
\big( 6K d^{1/2} + 2Kr + \vert j \vert \big)^4  \big( 2(\vert i + j \vert  +k) + 1 \big)^d R^2 K^d \, ;
$$
\item  and, for a universal constant $C > 0$,
$$
\R^v_{[\upwell],\downwell}  (\Long) \leq C 
(\constfine u^d K^{-1} + 2)^{d-1}  \big( (2d-1) \e \big)^{A/\const - 6} \, .
$$ 
\end{enumerate}
\end{lemma}
{\bf Proof: (1).} For given $\omega \in \fine$ and $v \in  \fubuthree$, the set of nearest-neighbour edges in $\Z^d$ is partitioned in Definition~\ref{d.prime} into three disjoint subsets: discovered open edges; discovered closed edges; and undiscovered edges. We have $R^v_{[\upwell],\downwell}(e) = \omega'(e)$ for any undiscovered edge~$e$.
Recall also, from the final sentence of Definition~\ref{d.prime}, that
$$
\R^v_{[\upwell],\downwell}  \big( \fv[i,j,k] \big) = \PP' \Big( \omega' \in \Omega': R^v_{[\upwell],\downwell} (\omega') \in  \fv[i,j,k]  \Big) \, .
$$
It is perhaps useful to recall that, under the probability measure $\R^v_{[\upwell],\downwell}  (\cdot)$, $\upwellv$, when it is viewed, as it will be, as a function of $\omega$, is a well-defined and non-random subset of edges, specified by 
Definition~\ref{d.displace}.  Moreover, when a sample of this probability measure realizes $V = v$, the random variable $\upwell$ defined by the sample equals $\upwellv= \upwellv(\omega)$.
In what follows, we will consider this sample when various subsets of the event $\{ V = v\}$ occur, so that $\upwell$ is well-defined and takes the just stated value.

Call a $\const$-box~$B$ {\em undiscovered} if every nearest-neighbour edge in the big box with Centre~$B$ is undiscovered for the given $\omega \in \fine$ and $v \in  \fubuthree$. 
If the configuration $R^v_{[\upwell],\downwell}$ realizes $\clusteruthree \cap \{ V = v \}$, then any box in $\surfacek$ that lies in $\overline\midzone[-2]$---recall Definition~\ref{d.tropics}---but in neither $\overline\upwell[2]$ nor $\overline\downwell[2]$ is undiscovered: indeed, in the case that $R^v_{[\upwell],\downwell}$ realizes $\clusteruthree$, any discovered edge in $\midzone$ is incident to a vertex in either $\upwell$ or $\downwell$; since no edge in a big box with Centre equal to such a  $\const$-box is thus incident, each such edge is undiscovered.

When $R^v_{[\upwell],\downwell} (\omega') \in \fv[i,j,k]$, the $\const$-box path $\Psi$ offered by Lemma~\ref{l.short} lies in $\surfacek \cap \overline\midzone[-2]$ but in neither  $\overline\upwell[2]$ nor  $\overline\downwell[2]$. Thus, the $\const$-boxes along this path are undiscovered. 
 Lemma~\ref{l.surfacebad} implies that  the $\const$-boxes along $\Psi$ are bad.  The configuration $R^v_{[\upwell],\downwell} (\omega')$ belongs to $\fv[i,j,k]$ and thus realizes $\fine$. 
 The $\const$-box $D[\Psi]$ lies in $\overline\midzone[-2]$ at latitude $-i$ and intersects $\overline\downwell[2]$. Lemma~\ref{l.latcard}(1,3) thus implies that $D[\Psi]$ may adopt one of at most $5^d \big( 6K d^{1/2} + 2Kr + \vert i \vert \big)^4 R$ locations.
By Lemma~\ref{l.percodominate}, we see that the path $\Psi$, which begins at a neighbour of $D[\Psi]$, and has length $\vert i+j\vert+k$ (the absolute value is operative when one or other of $i$ and $j$ lies in $\llbracket -\lfloor 3d/2  \rfloor,-1\rrbracket$), lies in  the set of open vertices of a site percolation on $\Z^d$ of parameter $\e > 0$, where this parameter may be chosen to be arbitrarily small provided that a sufficiently high choice of $\const \in \N$ is made. And indeed  this percolation is independent of the conditioning that specifies the law $\R^v_{[\upwell],\downwell}(\cdot)$:  this conditioning is measurable with respect to the configuration in $\downwell$, $\upwell$, $\northzone$ and $\southzone$, a region that no big box with Centre in $\Psi$ intersects, so that Lemma~\ref{l.condbadperc} delivers the claimed independence. The probability~$p(\ell,\e)$ that such a site percolation contains an open path of length~$\ell \in \N$ emanating from the origin is readily seen to satisfy
\begin{equation}\label{e.percpath}
p(\ell,\e) \leq (2d-1)^{\ell -1} \e^\ell \, .
\end{equation}
 Thus, we see that, with $\ell = \vert i+j\vert +k$,   
$$
\R^v_{[\upwell],\downwell} \big( \fv[i,j,k] \big) \leq 5^d \big( 6K d^{1/2} + 2Kr + \vert i \vert \big)^4 R
 (2d-1) (2d-1)^{\ell -1} \e^\ell \, .
$$
Using $\const \geq 1$, 
we obtain Lemma~\ref{l.rfprob}(1).

{\bf (2).} We seek an upper bound on the cardinality of the set of $v \in  \fubuthree$ for which $\omega' \in \Omega'$ exists such that  $R^v_{[\upwell],\downwell} (\omega') \in  \fv[i,j,k]$.  
It is useful to partition the set of $v \in  \fubuthree$ into $(2K)^d$ classes, with members of a given class differing by an element of $(2K)^d \cdot \Z^d$. We do so because, as $v$ varies over the elements in a given class, the collection $\overline\upwell[2]$ shifts in a straightforward fashion, by vector displacement. 

In any case, if such $\omega'$ is to exist,
then there must exist a pair $\big(D(\Psi),U(\Psi) \big)$ of $K$-boxes,  the former in $\overline\downwell[2]$ at latitude~$-i$, and the latter in $\overline\upwell[2]$ at latitude~$j$, at $K$-box nearest-neighbour distance at most $\vert i+j \vert +k+1$. The number of admissible pairs of $K$-boxes is at most  
$$
\alpha_1 =  5^d \big( 6K d^{1/2} + 2Kr + \vert i \vert \big)^4 R \cdot
5^d \big( 6K d^{1/2} + 2Kr + \vert j \vert \big)^4 R
$$ 
by Lemma~\ref{l.latcard}(1,2,3,4). When such a pair is given, the variation of $v$ within a given class yields various outcomes at which the elements of the pair lie at  $K$-box distance at most $\vert i+j \vert +k+1$. An upper bound on the number of such $v$ is  $\alpha_2 =  \big( 2(\vert i + j \vert  +k) + 1 \big)^d$, because this is the number of $K$-boxes at such a distance from a given $K$-box. The number of classes into which the set of $v$ is divided is at most $\alpha_3 = (2K)^d$. The product $\alpha_1 \alpha_2 \alpha_3$ is  the upper bound in   Lemma~\ref{l.rfprob}(2).

{\bf (3).}  Under $\fine \cap \Long$, $\Psi$ is a path in bad boxes of length at least~$A/\const - 6$ that begins at a box whose distance (in the lattice~$\Z^d$) from the vertical line through $0 \in \Z^d$ is at most $\constfine u^d + 2K$, since this box $*$-neighbours an element of $\overline{\cset}_u$, and  $\vert V(\csetuthree) \vert \leq \constfine u^d$ by condition~$\fg$.
Elements in the path $\Psi$ lie in $\overline\midzone[-2]$ and in neither $\overline\upwell[2]$ nor $\overline\downwell[2]$, so the constituent boxes are undiscovered for the given $\omega \in \fine$. 
The probability of such a path emanating from a given box is at most $\big( (2d-1)\e \big)^{A/\const - 6}$ by Lemma~\ref{l.percodominate} and the bound~(\ref{e.percpath}).   Lemma~\ref{l.rfprob}(3) follows from $\e < (2d-1)^{-1}$. \qed

 Note that $\R^v_{[\upwell],\downwell}(\fine )$ equals  $\R^v_{[\upwell],\downwell} (\fine  , V=v )$ when $v$ belongs to the set $\llbracket -\constfine u^d,\constfine u^d \rrbracket^{d-1} \times \Z$; and equals zero when $v$ is any other element of $\Z^d$. 
 This justifies the first inequality as we write
 \begin{eqnarray*}
  & & \sum_{v \in  \fubuthree}\R^v_{[\upwell],\downwell}(\fine) \\
  & \leq & \sum_{v \in  \fubuthree} \sum_{i,j,k}   \R^v_{[\upwell],\downwell}  \big( \fv[i,j,k] \big) \, + \,  \sum_{v \in \llbracket -\constfine n^d,\constfine n^d \rrbracket^{d-1} \times \Z} \R^v_{[\upwell],\downwell} (\Long) \\
  & \leq & 
250^d  R^3 \const^d
        \sum_{i,j,k}    \big( 6K d^{1/2} + 2Kr + \vert i \vert \big)^8 
\big( 6K d^{1/2} + 2Kr + \vert j \vert \big)^4  \\
& & \qquad \qquad \qquad \qquad  \qquad \times \, \big( 2(\vert i + j \vert  +k) + 1 \big)^d 
        (2d-1)^{\vert i+j \vert +k}  \e^{\vert i+j \vert +k} 
         \\
        & & \qquad \qquad \qquad 
         \, \, + \, \,
  (2\constfine u^d + 1)^{d-1}  C 
(\constfine u^d K^{-1} + 2)^{d-1}  \big( (2d-1) \e \big)^{A/\const - 6} \, ,
 \end{eqnarray*}
 where Lemma~\ref{l.rfprob}(1,2) is invoked to bound above the first term in the latter inequality, and Lemma~\ref{l.rfprob}(3) offers control on the second such term.  
The final sum is over $i,j \geq -3$ and $k \geq 0$; it converges due to $(2d-1)\e < 1$. Recall that, in Proposition~\ref{p.finesum}, we suppose that $A/K$ is bounded below by a suitably high multiple of $\log u$.  
We find then that
  $\sum_{v \in  \fubuthree}\R^v_{[\upwell],\downwell}(\fine)$ is bounded above, uniformly in $\omega \in \fine$. This completes the proof of Proposition~\ref{p.finesum}. \qed

\subsection{Controlling the damage when the single slide fails}\label{s.damagecontrol}

The event $\catastrophe$ occurs when $\residualwell$ contains a vertex whose height differs from $h$ by more than $3A$; which is to say, there exists a vertex in $\residualwell$ that lies in none of $\northzone$, $\midzone$ and $\southzone$. This tall structure prevents renewal levels at all the heights that it occupies. This section is devoted to proving the next result, which shows the rarity of this untoward outcome of slide resampling.

\begin{proposition}\label{p.cat}
For any positive $R$ and $K$, there exist $c = c(R,K) > 0$ and $C = C(K,d)$ such that the condition that $A \geq C \log u$ implies that 
$$
\PP^\fine \Big( \, \catastrophe  \, \Big\vert \, [\upwell] , \downwell \, \Big)  \, \leq \, \exp \big\{ - cA \big\} \, .
$$
\end{proposition}

The slide resample has been the centre of attention, but we have never explicitly defined it.   
In order to prove Proposition~\ref{p.cat}, it is useful to do so.  The slide operation $\slide$ 
at height~$h$
associates to each configuration $\omega \in \Omega$ a random configuration $\slide(\omega)$. Let the positive parameters $R$ and $K$ that specify the event $\fine$ be given. 
\begin{definition}\label{d.slide}
If $\omega \not\in \fine$, set $\slide(\omega) = \omega$. If $\omega \in \fine$, we set $\slide(\omega)$ to be a random element of $\Omega$ whose law equals
 $\PP^\fine \big( \cdot \big\vert  [\upwell] , \downwell \big)$.
 \end{definition}
\begin{lemma}\label{l.rvsystem}
The probability space $(\Omega,\mc{B},\PP)$ may be augmented in order to support the system of random variables $\big\{ \slide(\omega): \omega \in \Omega \big\}$ specified in Definition~\ref{d.slide}.
\end{lemma}
{\bf Proof.} The value of $\slide(\omega)$ may be chosen so that it is determined by the evaluation at $\omega$ of $([\upwell],\downwell)$. 
There are only finitely many possible values for this pair. For each, the value of the concerned  $\slide(\omega)$ may be set by sampling the probability measure in Lemma~\ref{l.probform} with $B = \fine$ and $(\uclust,\dclust) = \big([\upwell](\omega),\downwell(\omega) \big)$. \qed

\begin{definition}\label{d.crossing}
Let $\ell,\ell' \in \Z$ with $\ell \leq \ell'$. A connected subgraph~$S$ of the slab $\slab_{\ell,\ell'}$ is called {\em crossing} for this slab if $V(S)$ contains an element 
in   $\vertbdry_\ell$
 and another in $\vertbdry_{\ell'}$. 
\end{definition}

A {\em north island} is a finite open connected component in $\northzone$ that is disjoint from $\upwell$. A north island is called {\em crossing} if it is crossing for $\northzone$.
 A north island is called {\em meagre} if for at least one-quarter of  north-zone integral height indices $i \in \llbracket h+A,h+3A-1 \rrbracket$
 there at most $K\macb$ vertices in the north island with height in $[i,i+1)$.
A north island is called {\em near} if the distance from  the vertical line running through $0 \in \Z^d$ of all of its vertices is at most $\constfine u^d$.  

A {\em south island} is a finite open connected component in $\southzone$ that is disjoint from $\downwell$.
The definitions above may equally be made for south islands, with the replacements $\northzone \to \southzone$, $h+A \to h -3A$ and $h+3A \to h-A$ made.

Let  $\cmni$ denote the event that there exists a north or south island that is crossing, meagre and near.
\begin{lemma}\label{l.cmni}
There exist positive $C$ and $c$ such that, whenever the parameter $A \in \N$ appearing in Definition~\ref{d.tropics}  satisfies $A \geq C \log u$, 
$$
 \R^v_{[\upwell],\downwell} \big(  \cmni \big) \leq \exp \big\{ - c A \big\} 
$$
for any $v \in \fubuthree$ and $\omega \in \clusteruthree$.
\end{lemma}

Let $(\uclust,\dclust)$ denote a possible value of $([\upwell],\downwell)$ as $\omega$ ranges over $\clusteruthree$, and let $v \in \fubuthree$.
Let $\seed$ denote the set of $x \in \vertbdry_{h+A} \cap \llbracket -\constfine u^d ,\constfine u^d  \rrbracket^d$ such that $x$ is an endpoint of an  edge in $\northzone$ that is undiscovered given the data $(\uclust,\dclust,v)$ according to Definition~\ref{d.prime}.
The vertex set of any near crossing island in $\northzone$ intersects $\seed$. (The set $\seed$ is so named because an iterative procedure will grow a cluster from a starting point therein in the proof of the next result.)
\begin{lemma}\label{l.x}
There exists $c > 0$ such that, for $\omega \in \clusteruthree$, $x \in D$ and $v \in \fubuthree$, 
the conditional probability under $\R^v_{[\upwell],\downwell}(\cdot)$ 
given that this value is adopted
that $x$ is an element in the vertex set of a meagre crossing island is at most $e^{-cA}$.
\end{lemma}
{\bf Proof of Lemma~\ref{l.cmni}.} 
We sum the estimate in Lemma~\ref{l.x} over $x \in \seed$ to obtain the lemma with a relabelling of~$c > 0$. \qed

{\bf Proof of Lemma~\ref{l.x}.} 
It is perhaps helpful to recall the meaning of the sought assertion. We take $\omega \in \clusteruthree$ given and consider the conditional probability
\begin{eqnarray*}
& & \R^v_{[\upwell],\downwell} \big(  \cmni \big) \\
& = & \PP' \Big( \omega' \in \Omega': R^v_{[\upwell],\downwell}(\omega') \, \,  \textrm{realizes the event} \, \, \cmni \Big) \, ,
\end{eqnarray*}
where $\omega$ determines the value of $([\upwell],\downwell)$.
The sought claim is that this probability decays exponentially in the north- and south-zone width~$A$, independently of $\omega$. 
To prove it, take given $x \in \seed$. First an heuristic overview. 
Our task is to examine the north island that contains $x$---which is to say, the open connected component in $\northzone$ containing the vertex $x$ under the auxiliary percolation~$\omega'$. We will iteratively form this component by locally realizing its edges in a process that begins with edges neighbouring $x$. If the resulting component is to cross $\northzone$, the process must run for at least $2A$ steps. At each step, certain edges at the boundary of the thus realized cluster remain to be examined. If, at a given step, there are at most a given bounded number of such edges---call this the `bounded case'---there is a uniformly positive probability that the process will end at this very step, because it may be that each such edge will turn out to be closed when it is realized. On the other hand, if there are a greater number of unexamined edges at the given step---`the unbounded case'---then it is likely that a large number of open edges will enter the forming cluster when they are examined. Thus, if the bounded case is more typical, the forming north island is likely to have died out before the $(2A)$\textsuperscript{th} step is reached, so that the island cannot be crossing. And if the unbounded case is more prevalent, then typically the growing island will often welcome in many new edges and its resulting form will not be meagre. In summary, an island can be crossing and meagre only by performing the rare feat of surviving for long without thriving.

To make this scheme precise, we specify the searching procedure in terms of some notation.
Let $i \in \llbracket 0, 2A-1 \rrbracket$, and set 
$S_i$ equal to the subgraph $\slab_{h+A+i,h+A+i+1}$.
Thus, the edge set of $\northzone = \slab_{h+A,h+3A}$ is partitioned into the sets $\big\{ E(S_i) :  i \in \llbracket 0, 2A-1 \rrbracket \big\}$.

We will iteratively specify edge sets $U_j$ and $O_j$ for $j \in \llbracket 0,\termin \rrbracket$.
The value of the index $j$ will rise until its terminal value $\termin$, which is at most $2A$.
Write $\undisc \subset E(\Z^d)$ for the set of undiscovered edges in the sense of Definition~\ref{d.prime}, and $\mc{O} \subset E(\Z^d)$ for the set of open edges. 
Set 
$$
C_0 = \{ x \} \, \, \, , \, \, \, U_0 = \big\{ e \in \undisc \cap E(S_0): e \, \, \textrm{is incident to} \, \, x \big\} \, \, \,  \textrm{and} \, \, \, O_0 = U_0 \cap \mc{O} \, .
$$
For $j \in \intint{2A}$, set $\termin = j$ and terminate the procedure if $O_{j-1} = \emptyset$. In the opposing case, 
set
$C_{j-1}$ equal to the open connected component of $x$ in $\slab_{h+A,h+A+j} = \cup_{i=0}^{j-1} S_i$; then further set
$$
 U_j = \Big\{ e \in \undisc \cap E(S_j): e \, \, \textrm{is incident to} \, \, V(C_{j-1}) \Big\} \, \, \, \textrm{and} \, \, \, O_j = U_j \cap \mc{O} \, .
$$
Noting that $e \in E(\Z^d)$
belongs to $S_i$ precisely when $(e) \cap F_{h+A+i} \not= \emptyset$ and $(e) \subset B_{h+A+i+1}$, we see that $e \in S_i$ implies that the minimum height of the endpoints of $e$
lies in $[h+A+i,h+A+i+1)$. We thus find that, whenever $O_j \not= \emptyset$,
\begin{equation}\label{e.heightclaim}
\textrm{$\csetuthree$ contains at least  $\vert O_j \vert/(2d-1)$ vertices whose height lies in $[h+A+i,h+A+i+1)$} \, .
\end{equation}
Write $O = \cup_{i = 1}^\termin O_i$.
Note that $O \subseteq E(N_x)$, where $N_x$  is the north island that contains the vertex~$x$ under $R^v_{[\upwell],\downwell}$.  
We analyse how $O$ forms. Use the shorthand 
$$
\macd = \kay \, .
$$
 We say that step $j$ is {\em fecund} if $\vert U_j \vert \geq \macd$ and {\em unproductive} if  $\vert O_j \vert \leq \kaypbytwomac$. 
Let~$\termin$ denote the number of steps for which the procedure runs. Let $M$ denote the number of steps $j \in \llbracket 0,\termin \rrbracket$ that are not fecund. To express the theme of the proof sketch in these terms, a meagre crossing island emanating from $x$ must arise from a procedure with many unproductive steps. For such a procedure,  there are either many steps that fail to be fecund, or  many steps that are fecund and unproductive.  Next is a result showing that each of these eventualities is a rarity.

\begin{lemma}\label{l.geom}
\leavevmode
\begin{enumerate}
\item
The random variable $M$ is stochastically dominated by a geometric random variable~$X$ with law $\PP (X = k) = \rho^{k-1} (1 - \rho)$ for $k \geq 1$, where $\rho = 1 - (1-p)^{K-1}$. 
\item There exists $q = q(p) \in (0,1)$ such that, 
for each $k \in \N_+$,  the probability that there are at least~$k$ fecund steps, and that, among the first $k$ such, there are at least $k/2$ unproductive steps, is at most $q^k$.
\end{enumerate}
\end{lemma}
{\bf Proof: (1).} At any given step that is not fecund, the size of the unexamined set is at most $\kaymac -1$. With probability at least $(1-p)^{\kaymac -1}$, each edge will be found to be closed when it is examined. The procedure will then terminate; otherwise, it will continue, with further independent examinations of unexamined edges. Thus $M$, the number of steps that are not fecund, is stochastically dominated by a geometric random variable that is at least one and whose success probability equals $(1-p)^{\kaymac -1}$. 

{\bf (2).} The probability that
a given fecund step is unproductive 
is a constant~$c$ that is strictly less than one-half. These outcomes are independent, so that, if the list of fecund steps is recorded, and the unproductive elements are marked, the collection of marked indices is stochastically dominated by a Bernoulli-$c$ sequence truncated at a random time.  From this, Lemma~\ref{l.geom}(2) follows. \qed

We {\em claim} that, if $x$ is a vertex in a crossing island, then the procedure duration $T$ is at least~$2A-1$.
To see this, let $m_i$ denote the maximum height of a vertex in $C_i$. 
Note that $x \cvl \leq h +A$ by $x \in \vertbdry_{h+A}$ and Lemma~\ref{l.geometrybasics}(1); thus, $m_0 \leq h+A$. Note that
 $m_{i+1} - m_i \leq r$, with $r = \max_{1 \leq k \leq d} \vert e_k \cvl \vert$ so that $m_{i+1} - m_i \leq 1$. 
 If $x$ is a vertex in a crossing island, then this island intersects $\vertbdry_{h+3A}$, so that 
  $m_\termin \geq h + 3A - r$ by Lemma~\ref{l.geometrybasics}(1). We find then that $T \geq 2A/r - 1$, so that the claim follows from $r \leq 1$.

The probability that there are among the procedure's steps at least $A$ that are not fecund is at most $\PP(M \geq A) \leq \rho^{A-1}$ by Lemma~\ref{l.geom}(1). 
If $x$ is a vertex in a crossing island and there are fewer than $A$ steps that are not fecund, it must be, in view of $\termin \geq 2A -1$, that there are at least $A$ fecund steps. 
The probability that there are among these fecund steps at least $A/2$ that are unproductive is at most $q^A$ by  Lemma~\ref{l.geom}(2).
Moreover,~(\ref{e.heightclaim}) shows that, when there are at least $A/2$ steps that are not unproductive, there are at least $A/2$ values $i \in \llbracket h+A,h+3A-1 \rrbracket$
such that the crossing island of which $x$ is an element has more than $\kaypbytwomac = \kaypbytwo$ vertices whose height lies in~$[i,i+1)$. 
In this event, this crossing island is thus not meagre. Making a suitably small choice of $c > 0$ completes the proof of Lemma~\ref{l.x}. \qed

 {\bf Proof of Proposition~\ref{p.cat}.}
 We will first argue that
\begin{equation}\label{e.finecat}
\fine \cap \catastrophe \subseteq \cmni \, .
\end{equation}
Indeed, suppose that the left-hand event occurs. Since $\residualwell$ contains a vertex~$v$ with $\vert v \cvl - h \vert > 3A$,  Lemma~\ref{l.residualwell}(2) implies that  there exists an open path in $\residualwell$ that lies in one or other of $\northzone$ or $\southzone$ and 
that crosses the slab in question between its lower and upper sides. 
The open path thus lies in a crossing island. This island is a subset of $\residualwell$ and thus of the well $\csetuthree$. 
Every vertex in this island lies at distance at most $\constfine u^d$ from the vertical line through $0 \in \Z^d$, because $\vert V(\csetuthree) \vert \leq \constfine u^d$ (by condition~$\fg$ of~$\fine$) and contains $0 \in \Z^d$. Thus, the island is near. Since $\fine$ occurs, this event's subtropical conditions~$\fsOne$ and~$\fsTwo$---recall Definition~\ref{d.fvf}---ensure that the island must also be meagre. 
 We thus see that $\cmni$ occurs, so that~(\ref{e.finecat}) has been verified.

From~(\ref{e.finecat}) and Lemma~\ref{l.probform}, and using the notation~(\ref{e.condnot}), we note that, for $\omega \in \fine$,
\begin{eqnarray*}
& &  \PP_{[\upwell],\downwell}^\fine \big(  \catastrophe \big) \\
&  \leq & \frac{\sum_{v \in \llbracket -\constfine u^d, \constfine u^d \rrbracket^{d-1}} \R^v_{[\upwell],\downwell} \big(  \cmni \cap \fine \big) }{ \sum_{v \in \fubuthree} \R^v_{[\upwell],\downwell} (\fine , V =v)} \, ,
\end{eqnarray*}
where the summand is thus restricted in the right-hand numerator because summands indexed by other values of $v \in \fubuthree$ report the probability of an event that entails that $\upwell$ contains a vertex whose distance from the vertical line containing the origin exceeds $\constfine u^d$, which event is inconsistent with the global condition~$\fg$ in~$\fine$.
Proposition~\ref{p.cat}
then follows from Lemmas~\ref{l.qfine} and~\ref{l.cmni}.  \qed

\subsection{The iterated slide: fast decay for inter-renewal distance}\label{s.iterate}

Here we prove Theorem~\ref{t.renewal}(2).
To do so, we will iterate the slide resample that led to Theorem~\ref{t.renewal}(1). In analysing the steps of the iteration, we will not invoke Theorem~\ref{t.renewal}(1) but rather the underlying Propositions~\ref{p.success} and~\ref{p.cat}, which show that each step often acts favourably with a negligible probability of catastrophic failure. 

The level $h \in [ 0, u ]$ has been given, 
and we have often omitted it in describing the slide resample. 
Now, however, we will treat $h$ as variable, and indicate its role by placing it in a subscript of relevant notation, in usages such as $\fine_h$, $\upwell_h$, $\slide_h$ and $\catastrophe_h$. In the cases where the subscript role is already reserved, namely $\clusteruthree$ and $\csetuthree$, this event and random variable have no $h$-dependence, and it should be understood that these notions remain unaltered.

For convenience, recall that $\fine_h$ from Definition~\ref{d.fvf}  occurs when $\clusteruthree$ does and
 \begin{enumerate}
\item for $i \in \llbracket -A, -1 \rrbracket \cup \intint{A}$, $\big\vert \big\{ v \in V(\csetuthree): v \cvl \in [ h+i, h+i+1) \big\} \big\vert \leq i^3 R$; 
\item for at least one-quarter of $i \in I $, $\big\vert \big\{ v \in V(\csetuthree): v \cvl \in [h+i,h+i+1) \big\} \big\vert \leq K\macb$,
when $I = \llbracket h+A,h+3A-1 \rrbracket$; and likewise when $I = \llbracket h-3A,h-A-1 \rrbracket$; and
\item $\vert V(\csetuthree) \vert  \leq \constfine u^d$.
 \end{enumerate}

\begin{figure}[htbp]
\centering
\includegraphics[width=0.44\textwidth]{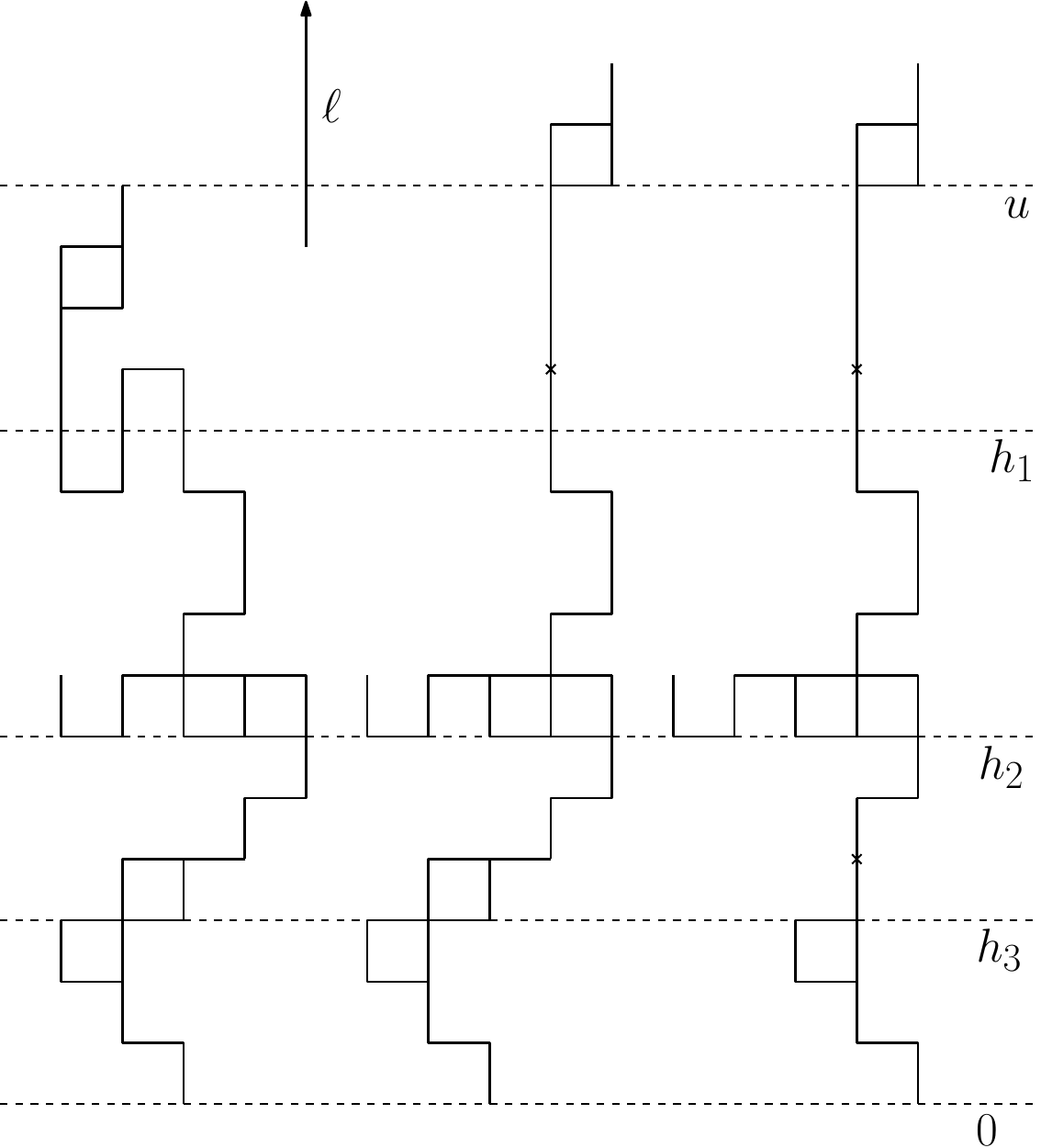}
\caption{The slide resample iterated, here with $L=3$ steps. The first and third steps are active; while the second is dormant, due to the high cardinality of the cluster around height~$h_2$.
On the left, the cluster under $S_0$ is shown; in the middle, this cluster under $S_1 = S_2$; on the right, the output cluster, under $S_3$. The active steps generate simple joins, resulting in two renewal points, marked with crosses, in the output.}\label{f.slideiterated}
\end{figure}

Let $h \in [ 0,u ]$. The random operator $\slide_h$ has been specified in Definition~\ref{d.slide}. This operator maps the event $\clusteruthree$ to itself and preserves the data $([\upwell],\downwell)$.
If $\omega$ has the law $\PP \big( \cdot \big\vert \csetuthree \cap \fine_h \big)$, then $\slide_h(\omega)$ is conditionally independently of $\omega$ given the value of $([\upwell],\downwell)$ at $\omega$.
We now offer a definition of an iterated slide resample, in which this resampling is undertaken at several heights. 
\begin{definition}
Let $\overline{h} = \big\{ h_i: i \in \intint{L} \big\}$ be a decreasing vector of elements of $[0,u]$. We now specify an iterative resampling procedure $\iter_{\overline{h}}: \clusteruthree \to \clusteruthree$, 
 $$
 \iter_{\overline{h}} \, = \, \slide_{h_L} \circ \cdots \circ \slide_{h_2} \circ \slide_{h_1} \, . 
 $$
 In this definition, the randomness used at each stage is independent of all earlier or later used randomness. 
 \end{definition}
 It is perhaps helpful to review this definition, and we set down a little notation by doing so; the reader may consult Figure~\ref{f.slideiterated}.

 The operator $\iter_{\overline{h}}$ acts on an input configuration $\omega$ that realizes $\clusteruthree$. 
 The form of its action is a composition of $L$ steps. 
 Set $S_0 = \omega$ and $S_i = \slide_h(S_{i-1})$ for $i \in \intint{L}$, evoking randomness in each step that is independent of that used at other steps. Then we take $\iter_{\overline{h}}(\omega) = S_L$. 
 The first step is illustrative of the general one.  The form of this step depends on whether $\omega$ realizes the event $\fine_{h_1}$. If it does not, the first step is dormant, so that $\omega$ remains the state of the system at the end of the first step: which is to say, $S_1 = \omega$. If $\omega \in \fine_{h_1}$, then the first-step ending-state $S_1$ is set equal to a realization of the law of the percolation measure $\PP$
 subject to the event $\clusteruthree \cap \fine_{h_1}$ occurring and to the data $\downwell_{h_1}$ and $\big[\upwell_{h_1} \big]$ coinciding with its specification under $S_0 = \omega$. 
 The height-$h$ hyperplane---the equator---moves downwards at each step, 
 visiting the successive levels in the vector $\overline{h}$, with a dormant or active step taking place at each level. The cluster $\csetuthree$ is a function of $S_i$ that is being held at equilibrium while being modified at each step. It is our aim to show that each modification introduces a renewal point at the concerned height with a uniformly positive probability, and further that earlier successes, in which renewal points have been secured higher up in the well, are unlikely to be undone by the ongoing resampling dynamic.

Let $i \in \intint{L}$.  The $i$\textsuperscript{th} step of the procedure is {\em regular} if either $S_{i-1} \not\in \fine_{h_i}$ or  $S_{i-1} \in \fine_{h_i}$ and $S_i \not\in \catastrophe_{h_i}$.
 The operator 
$\iter_{\overline{h}}$ is said to {\em act normally} if every step, indexed by $i \in \intint{L}$, is regular.
 
Let $j \in \llbracket 0, L \rrbracket$. An element $h_i$, with $i \in \intint{L}$, in the vector $\overline{h}$ is called {\em fine in $S_j$} if the configuration $S_j$ realizes the event $\fine_{h_i}$.

For $L \in \N$, let  $\big\{ h_i: i \in \intint{L} \big\}$ be a decreasing sequence of elements of $(0,u)$. We further set $h_0 = u$ and $h_{L+1} = 0$. Such a sequence is {\em well-spaced} if $h_{i+1} - h_i < -6A$ for $i \in \llbracket 0, L \rrbracket$.
\begin{proposition}\label{p.enduringveryfine}
Let $L \in \N$ and let $\overline{h} = \big\{ h_i: i \in \intint{L} \big\}$ be well-spaced. Suppose that 
$\iter_{\overline{h}}$ acts normally.
Let $j \in \llbracket 0, L \rrbracket$. If $h_j$ is very fine in $S_0$, then $h_j$ is very fine in $S_j$.
\end{proposition}
For  $j \in \llbracket 0, L \rrbracket$ and $h \in \R$, 
let $C_{j,h}$ denote the cardinality of the intersection of $\vertbdry_h$ 
and the set $V(\csetuthree)$ as specified by the configuration~$S_j$. 
 \begin{lemma}\label{l.regular} 
 Suppose that 
$\iter_{\overline{h}}$ acts normally. Let $i \in \intint{L}$ and $\ell \in \N$, with $\vert \ell \vert \leq 3A$.
Then the sequence $\llbracket 0, i-1 \rrbracket: j  \to C_{j,h_i + \ell}$ is non-increasing. 
\end{lemma}
We have seen that the value of the well $\csetuthree$ is variable during the $L$-step procedure that specifies 
 $\iter_{\overline{h}}$. We write $\csetuthree(S_i)$ to indicate the form at the end of the $i$\textsuperscript{th} step when the algorithm reports state~$S_i$. Naturally we make analogous uses such as 
 $\upwell_{h_i}(S_j)$: here, the denoted object is the upcluster specified by equatorial height~$h_i$ for the well $\csetuthree(S_j)$.
 
{\bf Proof of Lemma~\ref{l.regular}.}  Let $j \in \llbracket 0, L \rrbracket$.
At the $j$\textsuperscript{th} step of the procedure, an input configuration $S_{j-1}$ is transformed into an output $S_j$. In this change, the upcluster $\upwell_{h_j}(S_{j-1})$ shifts to become $\upwell_{h_j}(S_j)$. Indeed, we may locate a vector $V_j \in \Z^d$ such that  $\upwell_{h_j}(S_j) = \upwell_{h_j}(S_{j-1}) + V_j$. 
Let $i \in \intint{L}$ and $j \in \llbracket 1,i \rrbracket$. We {\em claim} that
\begin{itemize}
\item when $S_j \not\in \catastrophe_{h_j}$,  
every $v \in V \big( \csetuthree(S_j) \big)$ with $v \cvl < h_j - 3 A$ is a vertex in $\downwell_{h_j}(S_j)$. 
\end{itemize}
To see this, note that a vertex $v$ in $\csetuthree(S_j)$ is a vertex of at least one of  $\upwell_{h_j}(S_j)$, $\downwell_{h_j}(S_j)$ and $\residualwell_{h_j}(S_j)$. Any vertex $w$ in $\upwell_{h_j}(S_j)$ satisfies $w\cvl \geq h_j$. Because $S_j \not\in \catastrophe_{h_j}$, any vertex $w$ in $\residualwell_{h_j}(S_j)$ satisfies $w\cvl \geq h_j - 3A$. Thus we see that a vertex $v$ in $\csetuthree(S_j)$ with $v\cvl < h_j - 3 A$ is indeed a vertex in $\downwell_{h_j}(S_j)$, as we claimed. 

The case $i=1$ in Lemma~\ref{l.regular}  being trivial, we let $i \in \llbracket 2, L \rrbracket$  and $j \in \llbracket 1, i-1 \rrbracket$. Because the vector $\overline{h}$ is well-spaced,  $h_i + \ell  < h_j - 3A$  whenever $\vert \ell \vert \leq 3A$.
Since 
$\iter_{\overline{h}}$ acts normally, either $S_{j-1} \not\in \fine_{h_j}$  or $S_j \not\in \catastrophe_{h_j}$. 
In the former case, $S_j = S_{j-1}$ and thus  $C_{j,h_i+\ell} = C_{j-1,h_i+\ell}$. 
In the latter case, we invoke  the claim to find that  any vertex $v$ in $\csetuthree(S_j)$ that  lies in $\vertbdry_{h_i + \ell}$ is a vertex in $\downwell_{h_j}(S_j)$.
Since $\downwell_{h_j}(S_j)= \downwell_{h_j}(S_{j-1})$ is a subgraph of $\csetuthree(S_{j-1})$, we see that the number of vertices in $\csetuthree(S_j)$ that  lie in $\vertbdry_{h_i + \ell}$ is at most  the number of vertices in $\csetuthree(S_{j-1})$ that lie in  $\vertbdry_{h_i+\ell}$: which is to say, $C_{j,h_i+\ell} \leq C_{j-1,h_i+\ell}$. Lemma~\ref{l.regular} has been proved. \qed

{\bf Proof of Proposition~\ref{p.enduringveryfine}.}
Let $j \in \llbracket 0, i-1 \rrbracket$. It is enough to argue that, if $h_i$ is very fine in $S_j$, then it is very fine in $S_{j+1}$. 
If $j=i-1$, note that $S_i = \slide_{h_i}(S_{i-1})$ and that the definition of the random map $\slide_{h_i}$ entails the desired assertion for this value of $j$. Suppose then that $j \leq i-2$.
The event $\veryfine_{h_i}$from Definition~\ref{d.fvf}  entails tropical, subtropical and global conditions: that the first or second hold under $S_{j+1}$ when the respective assertion does  under $S_j$ follows from Lemma~\ref{l.regular}.
Note that either $S_{j+1} = S_j$, or $S_{j+1}$ realizes the event $\veryfine_{h_j}$: either way, the global condition is verified for~$S_{j+1}$ when it is for~$S_j$. \qed

Let $i \in \N$ and set $h_i = i \cdot 12 A$. Let $\slab_i$ denote $\big\{ z \in r\Z:  \vert z - h_i \vert < 3A \big\}$. Set $L = \lfloor u ( 12 A  )^{-1} \rfloor - 1$.
The set $\big\{ \slab_i: 1 \leq i \leq L \big\}$ is a collection of consecutive elements of $r \cdot \N$ that is contained in $[0,u]$. Each element has  cardinality $6Ar^{-1}$ (up to a discrepancy of two), and 
the gap between 
consecutive members,
in the sense of the cardinality of the interval in $r\Z$ lying to the right of the lower and to the left of the higher, is also $6A r^{-1}$ (up to the same discrepancy). 

Let $\overline{H} = \big( H_i: i \in \intint{L} \big)$ denote a random vector whose components are selected independently. The law of $H_i$ is uniform on $\slab_i$.

We will consider the iterative random procedure $\iter_{\overline{H}}$, where the initial condition $S_0$ has the law $\PP ( \cdot \vert \clusteruthree)$,
 the random vector $\overline{H}$ is selected independently of~$S_0$, and the randomness that is invoked during the running of the procedure is independent of this other randomness. Let~$\Q$ denote the law that governs the procedure in this form.

Let $\alpha \in (0,1)$. We say that the procedure's {\em input is  $\alpha$-satisfactory} if, in the configuration~$S_0$, there are at least $\alpha  L$ indices $i \in \intint{L}$ for which $\veryfine_{H_i}$ occurs.
\begin{proposition}\label{p.satisfactory}
There exist $\alpha \in (0,1)$ and $c > 0$ such that, under $\Q$, input is $\alpha$-satisfactory with probability at least $1 - \exp \big\{ - c   L \big\}$.
\end{proposition}
{\bf Proof.} Let $a \in (0,1/5)$. Choose $K \in \N$ to be at least $Ca^{-1}$ for the constant $C$ in Proposition~\ref{p.veryfinelevels}.
Under~$\Q$, the input configuration $S_0$ has the law $\PP \big( \cdot \big\vert \clusteruthree \big)$. 
Recall from the same proposition
that $N$ denotes the number of  $h \in \N \cap [0,u]$ for which $\veryfine_h$ occurs. 
By this result, $N \geq (1-a)u$ except with probability at most 
$\exp \big\{ - c K a u \big\}$.  
In the event $N \geq (1-a)u$, at most $au/r$ levels $i \in r\N \cap (0,u]$ fail to realize the event $\veryfine_i$. A slab $\slab_i$ may be called {\em grim} if at least one-half of its members fail in this sense.
Since the cardinality of a slab and the gap between one slab and the next both equal $6A r^{-1}$ up to an error of two, the union $\cup_{i=1}^L \slab_i$, which is a subset of $r\Z \cap [0,u]$, comprises a proportion of at least $1/2 - o(1)$ of elements of $r\N \cap [0,u]$, where $o(1) \to 0$ in high $A$.
 By choosing $A$ to be at least a suitably high constant $A_0$, we thus ensure that, when $N \geq (1-a)u$,  
 the number of grim slabs is at most $5aL$, and the number of slabs that are not grim is at least  $(1-5a) L$. If $i$ indexes a not grim slab, $\veryfine_{H_i}$ occurs with probability at least one-half, independently of other such outcomes. In the event $N \geq (1-a)u$, then, the number of indices $i \in \intint{L}$ for which $\veryfine_{H_i}$ occurs stochastically dominates a binomial random variable with parameters $(1-5a)L$ and $1/2$. Let $\alpha$ be a given positive value less than $(1-5a)/2$. Proposition~\ref{p.satisfactory} follows from binomial tail bounds. \qed

We fix a value for  $\alpha \in (0,1)$ in accordance with Proposition~\ref{p.satisfactory}.
\begin{proposition}\label{p.outputrenewal}
Consider the law $\Q$ conditionally on input being $\alpha$-satisfactory. There exist positive constants $c_0$, $c_1$ and $c_2$ such that, when $L \leq c_0 u (\log u)^{-1}$, the  probability that the output $S_ L$ contains at least $c_1 L$ renewal levels is at least $1 - \exp \big\{ -c_2  L \big\}  -  L \exp \big\{ -c_2  u/L \big\}$.
\end{proposition}
{\bf Proof.}
Note first that a small enough choice of $c_0 > 0$ yields the hypothesis $\min \{ h,u-h,A \} \geq C \log u$ that permits the application of Propositions~\ref{p.success} and~\ref{p.cat} in the upcoming argument.

Let $S \subset \intint{L}$ denote the set of indices $i$ for which $\veryfine_{H_i}$ occurs in the input configuration~$S_0$. Under the conditional law that we consider, $\vert S \vert$ is at least $\alpha L$.
Let $i \in S$. Suppose that, in the first $i-1$ steps of the procedure that specifies $\iter_{\overline{H}}$, catastrophe has not occurred: that is, for $j \in \intint{i-1}$, $S_{j-1} \not\in \catastrophe_{H_j}$.
From the proof of Proposition~\ref{p.enduringveryfine}, we learn that $\veryfine_{H_i}$ occurs in $S_{i-1}$.  By Lemma~\ref{l.veryfineimpliesfine} and Proposition~\ref{p.success}, there is probability at least some positive quantity~$c_0$ that $H_i$ is a renewal level in $S_i$.

If $S_k \not\in \catastrophe_{H_k}$ for all $k \in \intint{L}$, then, for any $i \in \intint{L}$ for which $H_i$ is a renewal level in $S_i$, the value $H_i$ is also a renewal level in the output configuration $S_L$. Indeed, in any update during the running of the procedure after the $i$\textsuperscript{th} step, the   downcluster and the residual cluster will lie below the level $H_i$. Since the upcluster is preserved at each of these later steps, the renewal status of level~$H_i$ measured relative to the top cannot be jeopardised during this later part of the procedure.

We see then that the number of renewal levels in the output of the procedure is stochastically dominated by $B {\bf 1}_\mathsf{N}$, where $B$ is a binomial random variable with parameters $\alpha L$  and $c_0$, and $\mathsf{N}$ is  the event that $S_i \not\in \catastrophe_{H_i}$ for $i \in \intint{L}$. By Proposition~\ref{p.cat}, the quantity $\Q(\mathsf{N})$ is at least $1 -  L \exp \big\{ -c_2  u/L \big\}$ for some positive value of $c_2$. Set $c_1 = \alpha c_0/2$. We have that $B \geq c_1 L$
except with probability at most $\exp \big\{ - c \alpha  L \big\}$ for some positive $c$ (that does not depend on $\alpha$). Adjusting $c_2 > 0$ to be at most $c\alpha$, we obtain Proposition~\ref{p.outputrenewal}. \qed

\begin{proposition}\label{p.renewalplenitude}
There exists  $c > 0$ such that, for $k \leq cu(\log u)^{-1}$,
the probability under $\PP$ given $\clusteruthree$ that $\csetuthree$ has fewer than $k$ renewal levels is at most $\exp \big\{ - ck \big\} +   k \exp \big\{ -c  u/k \big\}$.
\end{proposition}
{\bf Proof.}
We consider the copy of the law  $\PP \big( \cdot \big\vert \clusteruthree \big)$ offered by the distribution of the output of the iterative procedure $\iter_{\overline{H}}$. We aim to set the parameter $L$ close to the value $k$. By choosing $A = \lfloor u/(12(k+1)) \rfloor$, we see that
$$
\textrm{$L = \lfloor u (12A)^{-1} \rfloor - 1$ lies between $k$ and $k + 1 + 12(k+1)^2/u + k^3 u^{-2} O(1)$} \, .
$$
Since our assumption entails that $k \leq c u$
for a suitably small $c > 0$, we have that $L \in \llbracket k, 2k \rrbracket$.
Thus, input is $\alpha$-satisfactory with probability at least $1 - \exp \big\{ - ck \big\}$ by Proposition~\ref{p.satisfactory}.
Provided that input is $\alpha$-satisfactory, the output configuration will contain at least $c_1 k$ renewal levels with probability at least $1 - \exp \big\{ - c_2 k \big\}  -  2k \exp \big\{ - 2^{-1}c_2  u/k \big\}$ by Proposition~\ref{p.outputrenewal}. By relabelling $k$,
we obtain Proposition~\ref{p.renewalplenitude} for a suitably small choice of $c > 0$. \qed

{\bf Proof of Theorem~\ref{t.renewal}(2:Cluster).} The statement is obtained from Proposition~\ref{p.renewalplenitude} by lowering the value of $c > 0$ if need be. \qed
 
We end the section by deriving Proposition~\ref{p.stingrenewal}, which is a string counterpart to Theorem~\ref{t.renewal}(1).

 

  Let $\starcsetuthree$ denote the event $\clusteruthree$ given that $\head(\mc{C})$ and $\base(\mc{C})$ are renewal points in~$R(\mc{C})$ (and we also permit semi-open intervals $[u,u+3r)$ in using this notation). 
  Clusters with such renewal points are close cousins of strings; Corollary~\ref{c.renewal} is the counterpart of Theorem~\ref{t.renewal}
  for such clusters, and Corollary~\ref{c.renewal}(1) is the principal tool for proving Proposition~\ref{p.stingrenewal}.
  
  {\bf Proof of Corollary~\ref{c.renewal}.}  The derivation coincides with that of Theorem~\ref{t.renewal} by considering the event $\clusteruthree^*$ (or $\clusteruthreeopen^*$, for the third part) in place of
   $\clusteruthree$. \qed

   A further straightforward claim is needed.
   \begin{lemma}\label{l.stingcluster}
   We have that
   $$
   \PP \big( \cluster^*_{[u-2r,u+r)}  \big) \geq c \, \PP \big( \sting_{[u,u+r)} \big)
   $$
   \end{lemma}
   {\bf Proof.} This follows from $\PP \big( \cluster^*_{[u-2r,u+r)}  \big) \geq \PP \big( \cluster^*_{[u,u+r]}\big) \geq c \PP \big( \sting_{[u,u+r]}\big)$, the latter bound due to a bounded number of independent open/closed choices rendering a realization of  $\sting_{[u,u+r]}$ into an instance of  $\cluster^*_{[u,u+r]}$.  \qed
   
   {\bf Proof of Proposition~\ref{p.stingrenewal}.}
    Lemma~\ref{l.stingcluster}   implies that
   \begin{equation}\label{e.clineq} 
    \PP \Big( \vert R(S) \vert \geq \delta u  \Big\vert \sting_{[u,u+r)} \Big)
      \geq   c \cdot \frac{\PP \Big( \vert R(S) \vert \geq \delta u  \, , \,  \sting_{[u,u+r)} \Big)}{\PP \big(  \cluster^*_{[u-2r,u+r)}  \big)} \, .
   \end{equation}
 It is straightforward to see that
 $$
 \PP \Big( \vert R(S) \vert \geq \delta u \, , \, \sting_{[u,u+r)} \Big) \, \geq \, c \cdot \PP \Big( \vert R(S) \vert \geq \delta u  - 3 \, , \, \sting_{[u-2r,u+r)} \Big)
 $$
 where on the right-hand side, $S$ denotes the string of minimal length that manifests the event $\sting_{[u+r-3,u+r)}$.
 Indeed, the occurrence of the above right-hand event entails that on the left provided that a bounded set of open/closed decisions outside the slab implicated in the right-hand event fall the right way.
 
 By restricting to a slab, we see that the last displayed right-hand probability is at least 
 $$
  \PP \Big(  \vert R(\mc{C})\vert \geq \delta u - 3 \, , \, \cluster^*_{[u-2r,u+r)} \Big) \, .
 $$
 Returning to~(\ref{e.clineq}), we find that  
  $$
 \PP \Big( \vert R(S) \vert \geq \delta u \, \Big\vert \, \sting_{[u,u+r)} \Big) \, \geq \,   c \cdot \PP \Big( \vert R(\mc{C} \vert \geq \delta u  - 3 \, \Big\vert \, \cluster^*_{[u-2r,u+r)} \Big) \, .
 $$
Corollary~\ref{c.renewal}(3) shows that the right-hand side is at least $c^2$. This completes the proof of Proposition~\ref{p.stingrenewal}. \qed

\bibliographystyle{plain}

  \bibliography{traps}

\end{document}